\documentclass[10pt,reqno]{article}

\usepackage{amsmath, amsthm, amscd, amsfonts, amssymb, graphicx, color,mathrsfs,mathtools,enumerate}
\usepackage{lineno}
\usepackage{authblk}
\usepackage{dsfont}
\usepackage[bookmarksnumbered, colorlinks, plainpages]{hyperref}
\hypersetup{colorlinks=true,linkcolor=red, anchorcolor=green, citecolor=cyan, urlcolor=red, filecolor=magenta, pdftoolbar=true} 
\modulolinenumbers[5]

\newtheorem{thm}{Theorem}[section]
\newtheorem{lemma}[thm]{Lemma}
\newtheorem{proposition}[thm]{Proposition}
\newtheorem{corollary}[thm]{Corollary}
\newtheorem{remark}[thm]{Remark}
\newtheorem{defn}[thm]{Definition}
\newtheorem{hyp}[thm]{Hypothesis}

\theoremstyle{definition}

\theoremstyle{remark}
\newtheorem{rmk}[thm]{Remark}

\numberwithin{equation}{section}

\renewcommand{\hat}[1]{\widehat{#1}}

\newcommand{\gen}[1]{{\left\langle #1\right\rangle}}

\newcommand{\ra}{\rightarrow}

\newcommand{\N}{\mathbb{N}}

\newcommand{\R}{\mathbb{R}}

\newcommand{\elle}{\operatorname{L}}

\newcommand{\A}{{\operatorname{\mathcal{A}}}}

\newcommand{\fcon}{\operatorname{\mathscr{F}C}}

\allowdisplaybreaks

\title{Equivalence of Sobolev norms with respect to weighted Gaussian measures}

\author{{D. Addona}
\thanks{D.A. is member of G.N.A.M.P.A. of the Italian Istituto Nazionale di Alta Matematica (INdAM)}
\thanks{email: davide.addona@unipr.it}}
\affil{Department of Mathematical, Physical and Computer Sciences\\
University of Parma\\
Parco Area delle Scienze 53/a (Campus), 43124 Parma, Italy}
\date{}
\providecommand{\keywords}[1]{{\textit{Keywords}:} #1}
\providecommand{\subjclass}[1]{{\textit{SubjClass}[2010]:} #1}

\begin{document}

\frenchspacing

\maketitle 

\begin{abstract}
We consider the spaces $\elle^p(X,\nu;V)$, where $X$ is a separable Banach space, $\mu$ is a centred non-degenerate Gaussian measure, $\nu:=Ke^{-U}\mu$ with normalizing factor $K$ and $V$ is a separable Hilbert space. In this paper we prove a vector-valued Poincar\'e inequality for functions $F\in W^{1,p}(X,\nu;V)$, which allows us to show that for every $p\in(1,+\infty)$ and every $k\in\N$ the norm in $W^{k,p}(X,\nu)$ is equivalent to the graph norm of $D_H^k$ (the $k$-th Malliavin derivative) in $\elle^p(X,\nu)$. To conclude, we show exponential decay estimates for $(T^V(t))_{t\geq0}$ as $t\rightarrow+\infty$. Useful tools are the study of the asymptotic behaviour of the scalar perturbed Ornstein-Uhlenbeck $(T(t))_{t\geq0}$, and pointwise estimates for $|D_HT(t)f|_H^p$ by means both of $T(t)|D_Hf|^p_H$ and of $T(t)|f|^p$.
\end{abstract}

\vspace{5mm}
\keywords{{
Vector-valued Poincar\'e inequaility; Sobolev spaces; asymptotic behaviour; abstract Wiener spaces; vector-valued perturbed Ornstein-Uhlenbeck semigroup }}

\vspace{2mm}
\subjclass{28C20; 46E35; 47D06}

\section{Introduction}
Let $X$ be a separable Banach space, and let $(X,\mu,H)$ be an abstract Wiener space, i.e., $\mu$ is a centred non-degenerate Gaussian measure on $X$ and $H$ is the Cameron-Martin space associated to $\mu$. Further, let $\nu:=Ke^{-U}\mu$, where $U:X\rightarrow \R$ is a suitable function (see Hypothesis \ref{hyp:U}) and $K=\|e^{-U}\|_{\elle^1(X,\mu)}^{-1}$ is the normalizing constant which gives $\nu(X)=1$, and let $V$ be a separable Hilbert space. 
The aim  of this paper is to generalize to the spaces $\elle^p(X,\nu)$ and $\elle^p(X,\nu;V)$ some important results which are already known in the $\elle^p(X,\mu)$ setting.

Abstract Wiener spaces have been introduced by Gross in \cite{Gr67} to study the properties of Gaussian measures on infinite-dimensional spaces. The fundamental idea in the theory of Gaussian measures is that every centred Gaussian measure is the realization of the same "canonical'' Gaussian measure: the countable product of the standard normal Gaussian distributions on the line. This fact, and the fact that in infinite dimension does not exist an analogous of Lebesgue measure, have increased the interest around Gaussian measures in infinite dimension, both from an analytic and a probabilistic point of view (see e.g.  \cite{Bog98,Nu06,Sh04,Us95}). In an abstract Wiener space, the Gaussian measure $\mu$ factors according to an orthogonal decomposition of $H$, and this implies that many results can be obtained arguing in finite dimension and then letting the dimension to infinity. The situation completely change when $\mu$ is replaced by the measure $\nu$; in this case the finite dimensional approximations do not always work, and even when it happens computations are delicate and much more complicated.

The interest in weighted Gaussian measures $\nu$ in infinite dimension increases in last years, since in general they represent a class of infinite dimensional measures which are not decomposable along  directions of $H$. Further, by means of the Malliavin derivative $D_H$, on $\elle^p(X,\nu)$ it is possible to define a strongly continuous  semigroup $(T_p(t))_{t\geq0}$ whose infinitesimal generator $L_p$ is a perturbation of the Ornstein-Uhlenbeck operator. Features of $(T_p(t))_{t\geq0}$ and of $L_p$, such as maximal Sobolev regularity for solution to the elliptic problem $\lambda u-L_pu=f$, with $\lambda>0$ and $f\in \elle^2(X,\nu)$, smoothness properties of the semigroup $(T(t))_{t\geq0}$ and some functional inequalities are investigated both from an analytic and a probabilistic point of view also in more general setting, see \cite{AD20, ACF20, AngFerPal18, AngBigFer23, Big22,BigFer21, BigFer22, CF16,CF18,DPL14, DPL15}. We also refer to \cite{Fer19} for an in-depth analysis of the Sobolev spaces $W^{k,p}(X,\nu)$. 


The central point of our investigation is a characterization of the Sobolev spaces $W^{k,p}(X,\nu)$ with $k\in\N$ and $p\in(1,\infty)$ (see below for the definition of $W^{k,p}(X,\nu)$), obtained by means of a family of vector-valued Poincar\'e inequalities  and the Wiener chaos decomposition. To be more precise, we show that, for every $p\in[1,+\infty)$, there exists a positive constant $k_p$ such that for every $F\in W^{1,p}(X,\nu;V)$ it holds that
\begin{align}
\label{intr_poincare}
\|F-\nu(F)\|_{\elle^p(X,\nu;V)}\leq k_p\|\overline {D_{H}}F\|_{\elle^p(X,\nu;H\otimes V)},
\end{align}
where
\begin{align*}
\nu(F):=\int_XFd\nu \in V.
\end{align*}
We recall that a Poincar\'e inequality for weighted Gaussian measures has been already proved in \cite[Theorem 6.1]{FeUs00}, for smooth cylindrical real valued functions with $p=2$, in \cite[Theorem 6.2]{FeUs00} dropping the differentiability assumptions on $U$, and in \cite{AngFerPal18} on convex domains in Hilbert spaces. However, we deal with the closure of $D_H$ in $\elle^p(X,\nu)$ and therefore the belonging of $U$ to some Sobolev space is necessary in order to define the Sobolev spaces $W^{1,p}(X,\nu)$.

In the case of Gaussian measures, i.e., when $U=0$, these inequalities for vector-valued functions have been obtained in \cite{Bog98,Us95} when $p\in[1,+\infty)$ by exploiting properties of Gaussian measures when $V$ is a separable Hilbert space, 
and in \cite[Proposition 3.1]{PrVe14} for $p\in(1,+\infty)$ in the setting of Malliavin Calculus and when $V$ is a UMD space. In the last paper, the authors use a vector-valued version of Meyer's Multiplier Theorem, which does not work for $p=1$ and which cannot be adapted to our situation since it strongly relies on the decomposition of the Gaussian measure $\mu$, which is no longer available for $\nu$. We also mention that in \cite{VN15} the author provides a Poincar\'e inequality in Gaussian setting with respect to generalized directional gradients. 

In this specific case when $V=\mathcal H_k(H)$, $f\in W^{k+1,p}(X,\nu)$ and $F:=D_H^kf$, inequality \eqref{intr_poincare} reads as
\begin{align}
\label{intr_poincare_2}
\|D_H^kf-\nu(D_H^kf)\|_{\elle^p(X,\nu;\mathcal H_k(H))}\leq k_p\|D_{H}^{k+1}f\|_{\elle^p(X,\nu;\mathcal H_{k+1}(H))}.
\end{align}
Formula \eqref{intr_poincare_2} allows us to prove the characterization of the Sobolev spaces $W^{k,p}(X,\nu)$, with $k\in\N$ and $p\in(1,+\infty)$,  in terms the graph norm of $D_H^k$. The Sobolev spaces $W^{k,p}(X,\nu)$ have been introduced in \cite{Fer19} and they are defined as the domain of the closure of the operator
\begin{align*}
(D_H,\ldots,D_H^k):\fcon_b^\infty(X)\rightarrow \elle^p(X,\nu;H)\times\ldots\times \elle^p(X,\nu;\mathcal H_k(H)),
\end{align*}
in $\elle^p(X,\nu)$, with $k\in\N$ and $p\in(1,+\infty)$, endowed with the norm
\begin{align*}
\|f\|_{k,p}:=\left(\sum_{j=0}^k\|D_H^jf\|_{\elle^p(X,\nu;\mathcal H_{j}(H))}^p\right)^{1/p}, \quad f\in W^{k,p}(X,\nu).
\end{align*}
We prove that the question if $\|\cdot\|_{k,p}$ is equivalent to the graph norm 
\begin{align*}
\|f\|_{p,D_H^k}:=\|f\|_{\elle^p(X,\nu)}+\|D_H^kf\|_{\elle^p(X,\nu;\mathcal H_k(H))}, \quad f\in \fcon_b^\infty(X;V),
\end{align*}
has positive answer, and this extends the analogous result for the Gaussian measure $\mu$. The case $U=0$ and $p=1$ has been partially solved in \cite{AdMuRo21}, but it still remains an open problem. 

As a byproduct, we get the closability of the operator $D_H^k:\fcon_b^k(X)\rightarrow \elle^p(X,\nu;\mathcal H_k(H))$ in $\elle^p(X,\nu)$ when $U\in W^{k-1,q}(X,\mu)$ and a sufficient condition for the belonging of a function $f$ to $W^{k,p}(X,\nu)$. Finally, we provide an exponential decay for $(T^V(t))_{t\geq0}$.

The paper is organized as follows.

In Section \ref{sec:preliminaries} we collect the background which we need in the sequel of the paper. In particular, we present the structure of abstract Wiener spaces, listing the main features of the Cameron-Martin space, of the Ornstein-Uhlenbeck semigroup and of the Sobolev spaces $W^{k,p}(X,\mu)$. Then, we give the definition and the principal properties of the Wiener chaos decomposition which plays a crucial role in the proof of Theorems \ref{thm:equivalece_space_2} and \ref{thm:equivalence_space_k}. Later, we provide the assumptions on the function $U$, we introduce the Sobolev spaces $W^{k,p}(X,\nu)$ and $W^{1,p}(X,\nu;V)$, with $k\in\N$ and $p\in[1,+\infty)$, where $V$ is a separable Hilbert space. Besides the results in \cite{Fer19}, we show other properties of the elements of the functions belonging to these Sobolev spaces which will be used in the following. We also introduce the semigroup $(T(t))_{t\geq0}$ in $\elle^2(X,\nu)$ associated with the symmetric bilinear form
\begin{align*}
\mathcal E(u,v):=\int_X[D_Hu,D_Hv]_Hd\nu, \quad u,v\in W^{1,2}(X,\nu),
\end{align*}
and its vector-valued extension $(T^V(t))_{t\geq0}$ in $\elle^p(X,\nu;V)$, $p\in[1,+\infty)$. 

Section \ref{sec:analysis} is devoted to the study of the asymptotic behaviour of $(T(t))_{t\geq0}$ in $\elle^p(X,\nu)$ and to provide pointwise estimates for $|D_HT(t)f|_H^p$ by means both of $T(t)|D_Hf|_H^p$ and of $T(t)|f|^p$. We also generalize these results to the vector-valued semigroup $(T^V(t))_{t\geq0}$, of which we investigate both the asymptotic behaviour in $\elle^p(X,\nu;V)$ and pointwise estimates for $|\overline {D_H}T^V(t)F|_{H\otimes V}^p$ when $p\in[1,+\infty)$. The gradient estimates follow from computations similar to those in the proofs of \cite[Theorems 3.1 \& 3.3]{AngFerPal18} and are inspired by an idea of Bakry and \'Emery (see \cite{BaEm85}), while the asymptotic behaviour of $(T(t))_{t\geq0}$ is obtained with different techniques with respect to those in the quoted paper: here, we take advantage of the asymptotic estimate \eqref{asy_est_D_HT(t)} for $\|D_HT(t)f\|_{\elle^2(X,\nu;H)}$ as $t\rightarrow+\infty$, which allows us to characterize the limit of $T(t)f$ in $\elle^p(X,\nu)$ as $t\rightarrow+\infty$. The approximation results which we exploit to obtain the gradient estimates are postponed in Appendix \ref{app_A}.

Sections \ref{sec:poincare} and \ref{sec:equivalence} are the core of this paper. In the former we provide the vector-valued Poincar\'e inequalities \eqref{intr_poincare} and \eqref{intr_poincare_2}, in the latter we give the equivalence of the norms $\|\cdot\|_{k,p}$ and $\|\cdot\|_{p,D_H^k}$. We remark that the proof of \eqref{intr_poincare}, inspired by that in \cite{AdMuRo21}, relies on duality arguments: when $p\geq2$ we use the fact that, for every $F\in W^{1,p}(X,\nu)$, the function $F^*:=|F|^{p-2}F\in W^{1,p'}(X,\nu)$ and an explicit formula for $D_HF^*$ is available (see Lemma \ref{lemma:der_duale}), while when $p\in[1,2)$ we employ the duality between $\elle^p(X,\nu;H\otimes V)$ and $\elle^{p'}(X,\nu;H\otimes V)$. We stress that the constant $k_p$ which we obtain equals that in \cite{AdMuRo21} when $U=0$.

Regarding the equivalence of $\|\cdot\|_{k,p}$ and $\|\cdot\|_{p,D_H^k}$, we sketch the idea when $k=2$. In this case, it is enough to estimate $\|D_Hf\|_{\elle^p(X,\nu;H)}$ by means of $\|f\|_{p,D_H^2}$. From \eqref{intr_poincare_2} we have
\begin{align*}
\|D_Hf\|_{\elle^p(X,\nu;H)}
\leq & \|D_Hf-\nu(D_Hf)\|_{\elle^p(X,\nu;H)}+\|\nu(D_Hf)\|_H \\
\leq & k_p\|D_{H}^{2}f\|_{\elle^p(X,\nu;\mathcal H_2(H))}+\|\nu(D_Hf)\|_H.
\end{align*}
To conclude, it remains to prove that there exists a positive constant $c$ such that $\|\nu(D_Hf)\|_H\leq c\|f\|_{\elle^p(X,\nu)}$.  Similar same arguments allows us to extend the equivalence for $k\geq3$, provided further assumptions on $U$. The last result of this section is an improvement of the asymptotic behaviour of $(T^V(t))_{t\geq0}$, for which we show an exponential decay in $\elle^p(X,\nu)$ as $t\rightarrow+\infty$ for every $p\in[1,+\infty)$.

Finally, in Appendix \ref{app_A} we discuss the finite dimensional approximations which we use to prove the gradient estimates in Appendix \ref{app_B}. Many of the computations follow the lines of those in \cite{AngFerPal18}, but some slight modification is necessary since we weaken the assumptions on the weighted function $U$. Hence, we provide some details for reader's convenience.

\subsection*{Acknowledgements}
The author acknowledges that this research has financially been supported by the Programme ``FIL-Quota Incentivante'' of University of Parma and co-sponsored by Fondazione Cariparma.

\subsection{Notations}
Let $X$ be a real separable Banach space, let $X^*$ be its topological dual and let us denote by $\langle\cdot,\cdot\rangle$ and by $|\cdot|_X$ its duality and its norm, respectively. For every $k\in\N\cup\{\infty\}$ we denote by $C^k(X)$ the set of $k$-times (infinitely may times if $k=+\infty$) Fr\'echet differentiable functions $f:X\rightarrow \R$. We denote by $C_b^k(X)$ the set of functions $f\in C^k(X)$ which are bounded together with their derivatives up to order $k$ (of every order if $k=+\infty$). We denote by $\fcon_b^k(X)$ (resp. $\fcon^k(X)$) the set of functions $f\in C_b^k(X)$ (resp. $f\in C^k(X)$) such that there exist $n\in\N$, $x_1^*,\ldots,x_n^*\in X^*$ and $\varphi\in C_b^k(\R^n)$ (resp. $\varphi\in C^k(\R^n)$) such that $f(x)=\varphi(\langle x,x_1^*\rangle,\ldots,\langle x,x_n^*\rangle)$ for every $x\in X$. If $k=0$ we simply write $C(X), C_b(X), \fcon(X), \fcon_b(X)$. For every $f\in C_b(X)$ we set $\|f\|_\infty:=\sup\{|f(x)|: x\in X\}$.

Let $V$ be a separable Hilbert space with norm $|\cdot|_V$ and inner product $[\cdot,\cdot]_V$. For every $k\in \N\cup\{\infty\}$ we denote by $\fcon_b^k(X;V)$ the space of $k$-times (infinitely many times if $k=+\infty$) Fr\'echet differentiable functions $f:X\rightarrow V$ such that there exist $n\in\N$, $v_1,\ldots,v_n\in V$ and $f_1,\ldots,f_n\in \fcon_b^k(X)$ such that 
\begin{align*}
f(x)=\sum_{j=1}^nf_j(x)v_j, \quad x\in X.
\end{align*}
 For every $f\in C_b(X;V)$ we set $\|f\|_\infty:=\sup\{|f(x)|_V: x\in X\}$, where $C_b(X;V)$ denotes the space of functions $f:X\rightarrow V$ which are bounded and continuous with respect to the strong topology.

We denote by $B(X)$ the Borel $\sigma$-field of $X$. Let $\gamma$ be a Borel positive finite measure on $X$. For every $p\in[1,+\infty)$ we denote by $\elle^p(X,\gamma;V)$ the space of the (equivalence classes of) functions $F$ which are Bochner integrable endowed with the norm
\begin{align*}
\|F\|_{\elle^p(X,\gamma;V)}:=\left(\int_X|F|^p_Vd\gamma\right)^{1/p}, \quad F\in \elle^p(X,\gamma;V).
\end{align*}

Let $K$ be a separable Hilbert space with norm $|\cdot|_K$ and inner product $[\cdot,\cdot]_K$. We denote by $V\times K$ the set $\{(v,k):v\in V,\ k\in K\}$ and by $V\otimes K$ the tensor product of $V$ and $K$, i.e., the closure of span$\{v\otimes k:v\in V, \ k\in K\}$ with respect to the inner product
\begin{align*}
[v_1\otimes k_1,v_2\otimes k_2]_{V\otimes K}:=[v_1,v_2]_V[k_1,k_2]_K, \quad v_1,v_2\in V, \ k_1,k_2\in K.
\end{align*}
We set $\mathcal L(V;K)$ the space of bounded linear operators from $V$ to $K$. If $V=K$ we simply write $\mathcal L(V)$. We have $V\otimes K\subseteq \mathcal L(V;K)$ or $V\otimes K\subseteq \mathcal L(K;V)$ by setting $(v\otimes k)w:=[v,w]_Vk$ for every $v,w\in V$ and every $k\in K$, or $(v\otimes k)h:=[k,h]_Kv$ for every $v\in V$ and every $k,h\in K$, respectively.
For every $k\in\N$ we set
\begin{align*}
V^{\otimes k}:=\stackrel{\displaystyle k{{\rm -times}}}{\overbrace{V\otimes \cdots\otimes V}}
\end{align*}
As above, we can see $V^{\otimes 2}$ as a subset of $\mathcal L(V)$ by setting $(v\otimes w)(h)=[v,h]_Vw$ for every $v,w,h\in V$.

We say that a symmetric non-negative operator $A\in\mathcal L(V)$ is a Trace class operator if there exists an orthonormal basis $\{v_n:n\in\N\}$ of $V$ such that
\begin{align*}
{\rm Tr}[A]_V:=\sum_{n\in\N}[Av_n,v_n]_V<+\infty.
\end{align*}
It is possible to prove that the value of ${\rm Tr}[A]_V$ does not depend on the choice of the basis $\{v_n:n\in\N\}$. We denote by $\mathcal L_{(1)}^+(V)$ the space of Trace class operators endowed with the norm
\begin{align*}
\|A\|_{\mathcal L_{(1)}^+(V)}:={\rm Tr}[A]_V.
\end{align*}

For every $k\in\N$ we denote by $\mathcal H_k(V)$ the space of $k$-linear Hilbert-Schmidt operators on $V$, i.e., the space of the operators $A:V^k\rightarrow \R$ such that there exists an orthonormal basis $\{v_n:n\in\N\}$ of $V$  which gives
\begin{align*}
\|\A\|_{\mathcal H_k(V)}^2:=\sum_{i_1,\ldots,i_k=1}^\infty|\A(v_{i_1},\ldots,v_{i_k})|^2<+\infty.
\end{align*}
Again, the value of $\|\A\|_{\mathcal H_k(V)}$ does not depend on the basis $\{v_n:n\in\N\}$. The space $\mathcal H_k(V)$ with inner product
\begin{align*}
\langle A,B\rangle_{\mathcal H_k(V)}
:=\sum_{i_1,\ldots,i_k=1}^\infty A(v_{i_1},\ldots,v_{i_k})B(v_{i_1},\ldots,v_{i_k}),
\end{align*}
is a Hilbert space. Thanks to the Riesz representation theorem it is possible to identify $V$ with $\mathcal H_1(V)$ by setting, for every $h\in V$, $h(v):=[h,v]_V$, for every $v\in V$. For every $A\in\mathcal H_2(V)$ we have
\begin{align*}
\|A\|_{\mathcal H_2(V)}^2=\sum_{n\in\N}[Av_n,Av_n]_{V},
\end{align*}
where $\{v_n:n\in\N\}$ is any orthonormal basis of $V$. Further, $A^*A$ is a trace class operator and
\begin{align}
\label{traccia_hilbert_sc}
{\rm Tr}[A^*A]_V=\|A\|_{\mathcal H_2(V)}^2=\sum_{n\in\N}|Av_n|_V^2, \quad {\textrm{$\{v_n:n\in\N\}$ orthonormal basis of $V$}}.
\end{align}
We notice that for every $w_1,\ldots,w_k\in V$ the element $(w_1\otimes \cdots \otimes w_k)\in V^{\otimes k}$ belongs to $\mathcal H_k(V)$: indeed, for every orthonormal basis $\{v_n:n\in\N\}$ of $V$ we have
\begin{align}
\label{intro_k_prod}
\sum_{i_1,\ldots,i_k=1}^\infty|(w_1\otimes \cdots \otimes w_k)(v_{i_1},\ldots,v_{i_k})|^2
= & \sum_{i_1,\ldots,i_k=1}^\infty[w_1,v_{i_1}]_V^2[w_2,v_{i_2}]_V^2\cdots [w_k,v_{i_k}]_V^2 
= \prod_{j=1}^k|w_j|_V^2.
\end{align}

The writing $(a_n)$ denotes the sequence $\{a_n:n\in\N\}$, and $(a_{k_n})$ denotes the subsequence $\{a_{k_n}:n\in\N\}$ of $(a_n)$. For every $B\in B(X)$ we denote by $\chi_B$ the characteristic function of the set $B$, i.e., $\chi_B(x)=1$ if $x\in B$ and $\chi_B(x)=0$ if $x\in X\setminus B$.

For every $f:X\rightarrow \R$ we set $f^+:=\max\{f,0\}$ and $f^-:=\max\{-f,0\}$.

\section{Preliminaries}
\label{sec:preliminaries}
\subsection{The abstract Wiener space}\label{subsect_Hinfinity}

Let $X$ be a real separable Banach space and let $\mu$ be a centered non-degenerate Gaussian measure on $X$. We follow \cite[Chapter 2]{Bog98} to construct the Cameron--Martin space $H$ associated to $\mu$. This construction will give us the abstract Wiener space $(X,\mu,H)$ which will be the primary setting of our studies.
From the Fernique's theorem \cite[Theorem 2.8.5]{Bog98}, it follows that $X^*\subseteq L^2(X,\mu)$, and we denote by $j:X^*\rightarrow L^2(X,\mu)$ the injection of $X^*$ in $L^2(X,\mu)$, namely for every $x^*\in X^*$ we have
\begin{align*}
(j(x^*))(x):=\gen{x,x^*}, \quad x\in X.
\end{align*}
We remark that by \cite[Theorem 2.2.4]{Bog98} there exists a nonnegative symmetric linear bounded operator $Q:X^*\rightarrow X$ such that for every $x^*_1,x^*_2\in X^*$ we have
\begin{align}
\label{caratt_cova_inf}
\langle Q x^*_1,x^*_2\rangle=\int_X j(x^*_1)j(x^*_2)d\mu.
\end{align}
We denote by $X_{\mu}^*$ the closure of $j(X^*)$ in $L^2(X,\mu)$ and we define $R:X^*_{\mu}\rightarrow (X^*)'$ by
\begin{align}
\label{op_R}
R(f)(x^*):=\int_Xfj(x^*)d\mu, \quad f\in X_{\mu}^*, \ x^*\in X^*.
\end{align}
It is well-known that $R(X^*_{\mu})\subseteq X$. This means that for every $f\in X^*_{\mu}$ there exists $\mathcal{R}(f)\in X$ such that for every $x^*\in X^*$
\[\gen{\mathcal{R}(f),x^*}=\int_X fj(x^*) d\mu.\]
In particular the operator $\mathcal{R}:R(X^*_\mu)\ra X$ is the adjoint of $j$. 

We are now able to define the Cameron--Martin space $H$ (see \cite[Section 2.2]{Bog98}).

\begin{defn}
The Cameron--Martin space associated to $\mu$ is $H :=\left\{h\in X:|h|_{H}<+\infty\right\}$,
where 
\[|h|_{H} :=\sup\left\{\langle h,x^*\rangle\,\middle|\,x^*\in X^*\text{ such that } \|j(x^*)\|_{L^2(X,\mu)}\leq 1\right\}.\]
\end{defn} 

From \cite[Lemma 2.4.1]{Bog98} it follows that $h\in H$ if and only if there exists $\hat h\in X_{\mu}^*$ such that $\mathcal{R}(\hat h)=h$. Furthermore $H$ is a Hilbert space if endowed with the inner product
\begin{align}\label{prod_scal_Hinf}
[h,k]_{H}=\int_X\hat h\hat kd\mu, \quad h,k\in H.
\end{align}
We stress that for every $x^*\in X^*$, from \eqref{caratt_cova_inf} and \eqref{op_R} we have $Q x^*\in H$ and that $\mathcal{R}(j(x^*))=Q x^*$, i.e., $\widehat {Q x^*}=j(x^*)$. 
Further, from \eqref{prod_scal_Hinf} we deduce that
\begin{align}\label{inner_prod_inf_dual}
\langle Q x^*_1,x^*_2\rangle=[Q x^*_1,Q x^*_2]_{H}, \quad x^*_1,x^*_2\in X^*.
\end{align}
We provide the following characterization of $H$.

\begin{lemma}\label{caratt_cameron_martin_inf}
$H$ is the RKHS associated to $Q$ in $X$, i.e., $H=\overline{QX^*}^{|\cdot|_H}$, where $[Qx^*,Qy^*]_H:=\langle Qx^*,y^*\rangle$.
\end{lemma}

\begin{proof}
The proof is quite simple but we provide it for the convenience of the reader. Let $h\in H$. Then, there exists $\hat h \in X_{\mu}^*$ such that $\mathcal{R}(\hat h)=h$. In particular, there exists $(x^*_n)\subseteq X^*$ such that the sequence $(j(x^*_n))$ converges to $\hat h$ in $L^2(X,\mu)$. We claim that the sequence $(Q x^*_n)$ converges to $h$ in $H$. Indeed, by \eqref{prod_scal_Hinf} and recalling that $\widehat {Q x^*_n}=j(x^*_n)$ for every $n\in\N$ and \eqref{prod_scal_Hinf}, it follows that
\begin{align*}
\lim_{n\ra+\infty}|Q x^*_n-h|_{H}^2= \lim_{n\ra+\infty}\int_X|j(x^*_n)-\hat h|^2d\mu=0.
\end{align*}
This means that $\displaystyle H\subseteq\overline{Q X^*}^{|\cdot|_{H}}$. The converse inclusion follows by analogous arguments.
\end{proof}

Let us denote by $i$ the injection $i: H\ra X$, and let $i^*:X^*\rightarrow H$ be the adjoint operator of $i$ (here we have identified $H^*$ with $H$ by means of the Riesz representation theorem), then $Q=i\circ  i^*$. In particular, for every $x^*_1,x^*_2\in X^*$, by \eqref{caratt_cova_inf} and \eqref{prod_scal_Hinf} we have
\begin{align}\label{prop_covariance_op}
 [i^*x^*_1,i^*x^*_2]_{H}
 =\langle (i \circ i^*) x^*_1,x^*_2\rangle
=\langle Q x^*_1,x^*_2\rangle
=\int_X j(x^*_1)j(x^*_2)d\mu.
\end{align}

Finally we introduce the Ornstein-Uhlenbeck semigroup $(P(t))_{t\geq0}$: for every $f\in C_b(X)$ we set
\begin{align*}
P(t)f(x)=\int_Xf(e^{-t}x+\sqrt{1-e^{-2t}}y)\mu(dy), \quad x\in X, \ t\geq0.
\end{align*}
$(P(t))_{t\geq0}$ extends to a strongly continuous semigroup of contractions on $\elle^p(X,\mu)$ for every $p\in[1,+\infty)$ (see \cite[Theorem 2.9.1]{Bog98}). We suggest \cite[Section 2.9]{Bog98} for an in-depth study of $((P(t))_{t\geq0}$.

\subsection{The \texorpdfstring{$H$}{H}-gradient and the Sobolev spaces \texorpdfstring{$W^{k,p}(X,\mu)$}{W Xmu}}

In this subsection we define the Sobolev spaces $W^{k,p}(X,\mu)$. For an extensive treatment of these spaces we refer to \cite[Chapter 5]{Bog98}.
We consider the operator $D_{H}$ defined on $C^1(X)$ by
\begin{align*}
D_{H}F(x):=i^* DF(x), \ F\in C^1(X), \ x\in X.
\end{align*}
On the space $\fcon^1(X)$ the operator $D_{H}$ acts as follows:
\begin{align}
\label{i_inf-diff}
D_{H}f(x)=\sum_{j=1}^n\frac{\partial\varphi}{\partial\xi_j}(\langle x,x_1^*\rangle,\ldots,\langle x,x_n^*\rangle)i^*x^*_j, \quad x\in X,
\end{align}
where $f(x)=\varphi(\langle x,x_1^*\rangle,\ldots,\langle x,x_n^*\rangle)$ for some $n\in\N$, $\varphi\in C^1(\R^n)$ and $x_1^*,\ldots,x_n^*\in X^*$. We say that a function $f$ is $H$-differentiable at $x\in X$ if there exists $k\in H$ such that
\begin{align}
\label{H-diff}
f(x+h)=f(x)+[k,h]_H+o(|h|_H), \quad |h|_H\rightarrow0.
\end{align}
In this case we set $D_Hf(x):=k$. It is not hard to see that if $f\in\fcon_b^\infty(X)$ then the definitions of $D_H$ given in \eqref{i_inf-diff} and \eqref{H-diff} coincide.

We define the Sobolev spaces $W^{1,p}(X,\mu)$ as the domain of the closure of the operator $D_H:\fcon_b^\infty(X)\rightarrow \elle^p(X,\mu;H)$ in $\elle^p(X,\mu)$, for every $p\in[1,+\infty)$, and we denote by $D_H^\mu$ its closure. These spaces are Banach spaces if endowed with the norm
\begin{align*}
\|f\|_{W^{1,p}(X,\mu)}:=\left(\|f\|_{\elle^p(X,\mu)}^p+\|D^\mu_Hf\|_{\elle^p(X,\mu;H)}^p\right)^{1/p}, \quad f\in W^{1,p}(X,\mu),
\end{align*}
and the space $W^{1,2}(X,\mu)$ is a Hilbert space if endowed with the scalar product
\begin{align*}
\langle f,g\rangle_{W^{1,2}(X,\mu)}
:=\langle f,g\rangle_{\elle^2(X,\mu)}+\langle D^\mu_Hf,D_H^\mu g\rangle_{\elle^2(X,\mu;H)}, \quad f,g\in W^{1,2}(X,\mu).
\end{align*}
For every $k\in\N$ we define the operator $D_{H}^k$ which acts on elements of $\fcon^k(X)$ as follows:
\begin{align}
\label{k-H-der}
D_{H}^kf(x)=\sum_{j_1,\ldots,j_k=1}^n\frac{\partial^k\varphi}{\partial\xi_{j_1}\cdots\partial \xi_{j_k}}(\langle x,x_1^*\rangle,\ldots,\langle x,x_n^*\rangle)i^*x^*_{j_1}\otimes\cdots\otimes i^* x^*_{j_k}, \quad x\in X,
\end{align}
where $f(x)=\varphi(\langle x,x_1^*\rangle,\ldots,\langle x,x_n^*\rangle)$ for some $n\in\N$, $\varphi\in C^k(\R^n)$ and $x_1^*,\ldots,x_n^*\in X^*$. 
The Sobolev spaces $W^{k,p}(X,\mu)$ are the domain of the closure of the operator $(D_H,\ldots,D_H^k):\fcon_b^\infty\rightarrow \elle^p(X,\mu;H)\times\cdots\times \elle^p(X,\mu;\mathcal H_k(H))$ in $\elle^p(X,\mu)$, for every $p\in[1,+\infty)$, and we denote by $D_H^{\mu,k}$ its closure. These spaces are Banach spaces if endowed with the norm
\begin{align*}
\|f\|_{W^{k,p}(X,\mu)}:=\left(\|f\|_{\elle^p(X,\mu)}^p+\sum_{j=1}^k\|D^{\mu,j}_Hf\|_{\elle^p(X,\mu;\mathcal H_j(H))}^p\right)^{1/p}, \quad f\in W^{k,p}(X,\mu),
\end{align*}
and $W^{k,2}(X,\mu)$ is a Hilbert space if endowed with the scalar product
\begin{align*}
\langle f,g\rangle_{W^{k,2}(X,\mu)}
:=\langle f,g\rangle_{\elle^2(X,\mu)}+\sum_{j=1}^k\langle D^{\mu,j}_Hf,D_H^{\mu,j}g \rangle_{\elle^2(X,\mu;\mathcal H_j(H))}, \quad f,g\in W^{k,2}(X,\mu).
\end{align*}

\subsection{Wiener chaos decomposition}
\label{sub:wiener_chaos}
In this subsection we briefly present the Wiener chaos decomposition. For a systematic study of the Wiener chaos decomposition in the setting of Malliavin calculus we refer to \cite{Nu06}, while for the main results of the Wiener chaos decomposition in our setting we refer to \cite[Section 2.9]{Bog98}. For every $n\in\N$ let $H_n$ be the $n$-th Hermite polynomial, which is defined by
\begin{align*}
H_n(\xi):=\frac{(-1)^n}{n!}e^{\xi^2/2}\frac{d^n}{d\xi^n}\left(e^{-\xi^2/2}\right), \quad \xi\in\R.
\end{align*}
We notice that $H_1(\xi)=\xi$ and $H_2(\xi)=\frac12(\xi^2-1)$. We denote by $\Lambda$ the set of multiindices $\alpha=(\alpha_1,\ldots,\alpha_n,\ldots)$ such that all the terms, except a finite number of them, vanish. For every $\alpha\in \Lambda$ we set 
\begin{align*}
\alpha!:=\prod_{i=1}^\infty\alpha_i, \quad |\alpha|:=\sum_{i=1}^\infty\alpha_i.
\end{align*}
For every multiindex $\alpha$ we define the generalized Hermite polynomial $H_\alpha$ by
\begin{align*}
H_\alpha(\xi):=\prod_{i=1}^\infty H_{\alpha_i}(\xi_i), \quad \xi=(\xi_1,\ldots,\xi_n)\in\R^n, \ n=|\alpha|.
\end{align*}
Let us fix an orthonormal basis $\{h_n=i^*(x^*_n):\ x^*_n\in X^*, \ n\in\N\}$ of $H$ (this basis exists since $j(X^*)$ is dense in $X^*_\mu$), and for every $\alpha\in\Lambda$ let us define
\begin{align*}
\Phi_\alpha(x):=\sqrt{\alpha!}\prod_{i=1}^\infty H_{\alpha_i}(\langle x,x^*_i\rangle), \quad x\in X.
\end{align*}
It is possible to prove (see \cite[Proposition 1.1.1]{Nu06}) that the set $\{\Phi_\alpha:\alpha\in\Lambda\}$ is a complete orthonormal system in $\elle^2(X,\mu)$. Further, if for every $n\in\N\cup\{0\}$ we the set
\begin{align*}
E_n:={\rm  span}\{\Phi_\alpha:\alpha\in\Lambda, \ |\alpha|=n\},
\end{align*}
then $E_0=\R$, $E_1=X_\mu^*$, $E_n$ and $E_m$ are orthogonal subspaces of $\elle^2(X,\mu)$ when $n\neq m$ and
\begin{align*}
\elle^2(X,\mu)
:=\bigoplus_{n\in\N\cup\{0\}}E_n.
\end{align*}
Let us denote by $I_n$ the projection on $E_n$, $n\in\N\cup\{0\}$. Then, for every $f\in\elle^2(X,\mu)$ we have
\begin{align*}
f=\sum_{n=0}^\infty I_nf, \quad I_nf=\sum_{\alpha\in\Lambda, |\alpha|=n}\left(\int_Xf\Phi_\alpha d\mu\right)\Phi_\alpha.
\end{align*}
\begin{remark}
\label{rmk:chaos_ind_base}
The properties of $\Phi_\alpha$ and $E_n$ are independent of the choice of the orthonormal basis $\{h_n:n\in\N\}$ of $H$, i.e., if $\{k_n:k_n=i^*(y_n^*), \ y^*_n\in X^*, \ n\in\N\}$ is another orthonormal basis of $H$ and we set 
\begin{align*}
\Phi'_\alpha(x):=\sqrt{\alpha!}\prod_{i=1}^\infty H_{\alpha_i}(\langle x,y^*_i\rangle), \quad x\in X, \ \alpha\in\Lambda, \quad E'_n:={\rm  span}\{\Phi'_\alpha:\alpha\in\Lambda, \ |\alpha|=n\}, 
\end{align*}
then $E'_n$ and $E'_m$ are orthogonal subspaces of $\elle^2(X,\mu)$ when $n\neq m$,
\begin{align*}
\elle^2(X,\mu)
:=\bigoplus_{n\in\N\cup\{0\}}E'_n,
\end{align*}
and for every $f\in\elle^2(X,\mu)$ we have
\begin{align*}
f=\sum_{n=0}^\infty I'_nf, \quad I'_nf=\sum_{\alpha\in\Lambda, |\alpha|=n}\left(\int_Xf\Phi'_\alpha d\mu\right)\Phi'_\alpha.
\end{align*}
\end{remark}

\begin{remark}
\label{rmk:dec_wiener_chaos}
For every multiindex $\alpha$ with $|\alpha|=1$ it follows that $\Phi_\alpha(x)=\hat{h_n}(x)$ for every $x\in X$, where $\alpha_n=1$, with $n\in\N$, is the unique component of $\alpha$ different from $0$. Further, for every multiindex $\alpha$ with $|\alpha|=2$, we have $\Phi_\alpha(x)=\frac{1}{\sqrt2}((\hat{h_n}(x))^2-1)$ for every $x\in X$, if $\alpha_n=2$, with $n\in\N$, is the unique component of $\alpha$ different from $0$, while $\Phi_\alpha(x)=\hat{h_n}(x)\hat{h_m}(x)$ for every $x\in X$, if $\alpha_n=\alpha_m=1$, with $n,m\in\N$, $n\neq m$, are the unique components of $\alpha$ which do not vanish.
\end{remark}
The Ornstein-Uhlenbeck semigroup $(P(t))_{t\geq0}$ behaves good on $E_n$: indeed, (see \cite[Theorem 2.9.2]{Bog98}) for every $t\geq0$ and every $f\in\elle^2(X,\mu)$ we have
\begin{align*}
P(t)f=\sum_{n=0}^\infty e^{-nt}I_nf.
\end{align*}
In particular, for every $f\in E_n$ it follows that $T(t)f=e^{-nt}f$, for every $n\in\N$ and every $t\geq0$. 

We prove the following useful lemma which is a consequence of Nelson's hypercontractivity theorem for the Ornstein-Uhlenbeck semigroup $(P(t))_{t\geq0}$ (see \cite[Theorem 5.5.3]{Bog98}). 

\begin{lemma}
\label{lemma:L2-Lp-norm}
Let  $p\in(1,+\infty)$. The projection $I_k$ is a bounded linear operator from $\elle^p(X,\mu)$ onto $\elle^2(X,\mu)$, and
\begin{align*}
\|I_kf\|_{\elle^2(X,\mu)}\leq & \frac{1}{(p-1)^{k}}\|f\|_{\elle^p(X,\mu)}, \quad f\in\elle^p(X,\mu), \quad p\in(1,2), \\
\|I_kf\|_{\elle^2(X,\mu)}\leq & (p-1)^{k/2}\|f\|_{\elle^p(X,\mu)}, \quad f\in\elle^p(X,\mu), \quad p\in[2,+\infty), 
\end{align*} 
\end{lemma}
\begin{proof}
The case $p=2$ follows from the fact that $I_k$ is a projection on $\elle^2(X,\mu)$. If $p\in(2,+\infty)$ the statement is easy to prove. Indeed, let $p\in(2,+\infty)$ and let $f\in \elle^p(X,\mu)$. Then, $f\in\elle^2(X,\mu)$ and from \cite[Chapter IV, Theorem 1]{Us95} we have
\begin{align*}
 \|I_kf\|_{\elle^p(X,\mu)}\leq  (p-1)^{k/2}\|f\|_{\elle^p(X,\mu)}.
\end{align*}
Further,
\begin{align*}
\|I_kf\|_{\elle^2(X,\mu)}\leq \|I_kf\|_{\elle^p(X,\mu)},
\end{align*}
and combining the above inequalities we have the thesis.

Let $p\in(1,2)$, let $k\in\N$ and let $f\in \elle^p(X,\mu)$. We consider a sequence $(f_n)\subseteq C_b(X)\subseteq \elle^2(X,\mu)$ which converges to $f$ in $\elle^p(X,\mu)$ as $n\rightarrow+\infty$. We recall that for every element $g\in E_k$ we have $P(t)g=e^{-kt}g$. Then,
\begin{align*}
\|I_k f_n\|_{\elle^2(X,\mu)}
=e^{kt}\|P(t)I_k f_n\|_{\elle^2(X,\mu)}, \quad n\in\N.
\end{align*}
Further, the hypercontractivity property of $(P(t))_{t\geq0}$ (see \cite[Theorem 5.5.3]{Bog98}) gives
\begin{align*}
\|P(t)I_k f_n\|_{\elle^2(X,\mu)}
\leq \|I_kf_n\|_{\elle^p(X;\mu)}, \quad n\in\N,
\end{align*}
with $t=-\ln(p-1)^{1/2}$. It follows that
\begin{align*}
\|I_k f_n\|_{\elle^2(X,\mu)}
\leq e^{kt} \|I_kf_n\|_{\elle^p(X;\mu)}=\frac{1}{(p-1)^{k/2}} \|I_kf_n\|_{\elle^p(X;\mu)}, \quad n\in\N.
\end{align*}
\cite[Chapter IV, Theorem 1]{Us95} implies that
\begin{align*}
\|I_kf_n\|_{\elle^p(X,\mu)}
\leq \frac{1}{(p-1)^{k/2}}\|f_n\|_{\elle^p(X,\mu)}.
\end{align*}
Therefore, $(I_kf_n)$ is a Cauchy sequence in $\elle^2(X,\mu)$. If we denote by $I_kf$ the limit of $(I_kf_n)$ in $\elle^2(X,\mu)$ as $n\rightarrow+\infty$, it follows that $I_kf$ is well defined, it belongs to $\elle^2(X,\mu)$ and
\begin{align*}
\|I_kf\|_{\elle^2(X,\mu)}
\leq \frac{1}{(p-1)^{k}}\|f\|_{\elle^p(X,\mu)}.
\end{align*}
\end{proof}

\subsection{The Sobolev spaces \texorpdfstring{$W^{k,p}(X,\nu)$}{W1pXm} and the \texorpdfstring{$H$}{H}-divergence operator}

Let us provide the assumptions on the function $U$, which will be used to define the weighted Gaussian measure $\nu$.
\begin{hyp}
\label{hyp:U}
$U$ is a $H$-convex $|\cdot|_X$-lower semicontinuous function which belongs to $W^{1,p}(X,\mu)$ for every $p\in(1,+\infty)$. \end{hyp}

\begin{rmk}
\begin{enumerate}
\item From \cite[Theorem 4.3]{FeUs00} it follows that $U$ admits a modification $U'$ which is a Borel measurable convex function on $X$. Since $U=U'$ for $\mu$-a.e., in the following we simply denote by $U$ this modification, and therefore we can assume that $U$ is a convex function.
\item It is possible to weaken the assumptions on $U$ requiring that it belongs to $W^{1,p}(X,\mu)$ for some $p\in(1,+\infty)$. This implies that the results of this paper hold for exponents $q$ smaller than a fixed exponent $\overline p\in(1,+\infty)$ which depends on $p$.  
\end{enumerate}
\end{rmk}

Hereafter, we always assume Hypothesis \ref{hyp:U}. The convexity of $U$ implies that $U$ is bounded from below by a linear functional. From Fernique's theorem (see \cite[Theorem 2.8.5]{Bog98}) we infer that the function $e^{-U}$ belongs to $\elle^p(X,\mu)$ for every $p\in[1,+\infty)$. We introduce the probability measure
\begin{align}
\nu:=Ke^{-U}\mu, \quad K=\|e^{-U}\|^{-1}_{\elle^1(X,\mu)},
\end{align}
where $K$ is the normalizing factor. Hypothesis \ref{hyp:U} implies that $|D^\mu_HU|_H\in \elle^p(X,\nu)$ for every $p\in[1,+\infty)$. Indeed, for every $p\in[1,+\infty)$ we have
\begin{align*}
\int_X|D_H^\mu U|^p_Hd\nu
= & \int_X|D_H^\mu U|^p_H e^{-U}d\mu
\leq \|D_H^\mu U\|_{\elle^{pq}(X,\mu)}\|e^{-U}\|_{\elle^{q'}(X,\mu)}<+\infty,
\end{align*}
for every $q\in(1,+\infty)$. Analogously, we can prove that $\hat h\in\elle^q(X,\nu)$ for every $q\in[1,+\infty)$.

We define the Sobolev spaces $W^{1,p}(X,\nu)$ as in \cite{Fer19}. The integration by parts formula (see \cite[Lemma 3.1]{Fer19})
\begin{align}
\label{int_by_parts}
\int_X[D_{H}f,h]_{H}d\nu
= \int_Xf(\hat h+[D_{H}^\mu U,h]_{H})d\nu,
\end{align}
implies that the operator $D_{H}:\fcon_b^\infty(X)\rightarrow \elle^p(X,\nu;H)$ is closable for $p\in[1,+\infty)$. 
\begin{proposition}
\label{prop:closability_gradient}
The operator $D_{H}:\fcon_b^\infty(X)\rightarrow \elle^p(X,\nu;H)$ is closable in $\elle^p(X,\nu)$ for every $p\in[1,+\infty)$. We denote by $D_{H}$ its closure and by $W^{1,p}(X,\nu)$ the domain of its closure. $W^{1,p}(X,\nu)$ is a Banach space if endowed with the norm
\begin{align*}
\|f\|_{1,p}:=\left(\|f\|_{\elle^p(X,\nu)}^p+\|D_{H}f\|_{\elle^p(X,\nu;H)}^p\right)^{1/p}, \quad f\in W^{1,p}(X,\nu),
\end{align*}
$W^{1,2}(X,\nu)$ is a Hilbert space if endowed with inner product
\begin{align*}
\langle f,g\rangle_{W^{1,2}(X,\nu)}
= \langle f,g\rangle_{\elle^2(X,\nu)}
+ \langle D_{H}f,D_{H} g\rangle_{\elle^2(X,\nu;H)}, \quad f,g\in W^{1,2}(X,\nu).
\end{align*}
\end{proposition}
\begin{remark}
We have introduced the notation $D_H^\mu$ to underline the difference between the closure of $D_H$ in $\elle^p(X,\mu)$ and in $\elle^p(X,\nu)$; the importance of this different notation will arise in the proof of Proposition \ref{prop:gradiente_nullo_funz_costante}, where we use a property of $D_H^\mu$. This means that the notation $D_Hf$ represents both the action of the closure of $D_H$ in $\elle^p(X,\nu)$ for $f\in W^{1,p}(X,\nu)$, and the action of the operator $D_H$ defined in \eqref{i_inf-diff} on smooth functions $f\in C^1(X)$.
\end{remark}
\begin{proof}
If $p>1$ the statement has been already proved in \cite[Proposition 3.2]{Fer19}. Let us consider $p=1$ and let $\{h_m:m\in\N\}$ be an orthonormal basis of $H$. In this case we cannot directly use the method applied in the proof of \cite[Proposition 3.2]{Fer19}, since $\hat h_n-[D_H^\mu U,h_n]_{H}$ does not belong to $\elle^\infty(X,\nu)$. We introduce the function $\theta\in C_b^2(\R)$ such that $\theta(0)=0$ and $\theta'(0)\neq0$. Let $(f_n)\subseteq \fcon_b^\infty(X)$ be such that $f_n\rightarrow 0$ in $\elle^1(X,\nu)$ and $D_{H}f_n\rightarrow G$ in $\elle^1(X,\nu;H)$ as $n\rightarrow+\infty$. Then, for every $\varphi\in C_b^1(X)$ and every $m\in\N$ we have
\begin{align}
\label{p=1_inizio}
\int_X[D_{H}(\theta\circ f_n),h_m]_{H}\varphi d\nu
= \int_X(\theta'\circ f_n)[D_{H}f_n,h_m]_{H}\varphi d\nu.
\end{align}
Letting $n\rightarrow+\infty$, the right-hand side converges to $\int_X\theta'(0)[G,h_m]_{H}\varphi d\nu$. Indeed, 
\begin{align*}
& \left|\int_X(\theta'\circ f_n)[D_{H}f_n,h_m]_{H}\varphi d\nu-\int_X\theta'(0)[G,h_m]_{H}\varphi d\nu \right| \\
\leq & \left|\int_X(\theta'\circ f_n)[D_{H}f_n,h_m]_{H}\varphi d\nu-\int_X(\theta'\circ f_n)[G,h_m]_{H}\varphi d\nu \right| \\
& + \left|\int_X(\theta'\circ f_n)[G,h_m]_{H}\varphi d\nu-\int_X\theta'(0)[G,h_m]_{H}\varphi d\nu \right| \\
\leq & \int_X|(\theta'\circ f_n)\varphi ||[D_{H}f_n,h_m]_{H}-[G,h_m]_{H}|d\nu
+ \int_X|(\theta'\circ f_n)-\theta'(0)||[G,h_m]_{H}\varphi |d\nu.
\end{align*}
Since $\theta'$ and $\varphi$ are bounded, the first integral vanishes as $n\rightarrow+\infty$. Further, there exists a subsequence $(f_{k_n})\subseteq (f_n)$ such that $f_{k_n}(x)\rightarrow 0$ for $\nu$-a.e. $x\in X$. Moreover,
\begin{align*}
|(\theta'\circ f_n)-\theta'(0)||[G,h_m]_{H}\varphi |
\leq 2\|\theta'\|_\infty\|\varphi\|_\infty |G|_H.
\end{align*} 
By the dominated convergence theorem we infer that
\begin{align}
\label{p=1_convergenza_secondo_integrale}
 \int_X|(\theta'\circ f_{k_n})-\theta'(0)||[G,h_m]_{H}\varphi |d\nu\rightarrow 0, \quad n\rightarrow+\infty.
\end{align}
In particular, we have prove that every subsequence $(f_{k_n})\subseteq (f_n)$ admits a subsequence $(f_{k_{n_s}})\subseteq (f_{k_n})$ such that \eqref{p=1_convergenza_secondo_integrale}, with $(f_{k_n})$ replaced by $(f_{k_{n_s}})$ holds true. Hence,
\begin{align*}
 \int_X|(\theta'\circ f_{n})-\theta'(0)||[G,h_m]_{H}\varphi |d\nu\rightarrow 0, \quad n\rightarrow+\infty.
\end{align*}
Let us apply \eqref{int_by_parts} to the left-hand side of \eqref{p=1_inizio}. We get
\begin{align*}
\int_X[D_{H}(\theta\circ f_n),h_m]_{H}\varphi d\nu
=\int_X (\theta\circ f_n)(\hat{h_m}\varphi+[D_{H}^\mu U,h_m]_{H}\varphi -[D_{H}\varphi,h_m]_{H})d\nu.
\end{align*}
As above, let $(f_{k_n})\subseteq (f_n)$ be such that $f_{k_n}(x)\rightarrow 0$ for $\nu$-a.e. $x\in X$. Since
\begin{align*}
& \left|(\theta\circ f_n)(\hat{h_m}\varphi+[D_{H}^\mu U,h_m]_{H}\varphi -[D_{H}\varphi,h_m]_{H})\right| \\
\leq &\|\theta\|_\infty(|\hat{h_m}|\|\varphi\|_\infty+|D_H^\mu U|_{H}\|\varphi\|_\infty+\|D_{H}\varphi\|_\infty)\in\elle^1(X,\nu),
\end{align*}
by the dominated convergence theorem we infer that
\begin{align}
\label{p=1_conv_terzo_int}
=\int_X (\theta\circ f_{k_n})(\hat{h_m}\varphi+[D_H^\mu U,h_m]_{H}\varphi -[D_{H}\varphi,h_m]_{H})d\nu\rightarrow0,\quad n\rightarrow+\infty.
\end{align}
In particular, we have prove that every subsequence $(f_{k_n})\subseteq (f_n)$ admits a subsequence $(f_{k_{n_s}})\subseteq (f_{k_n})$ such that \eqref{p=1_conv_terzo_int} holds true with $(f_{k_n})$ replaced by $(f_{k_{n_s}})$. Hence,
\begin{align*}
\int_X (\theta\circ f_{n})(\hat{h_m}\varphi+[D_H^\mu U,h_m]_{H}\varphi -[D_{H}\varphi,h_m]_{H})d\nu\rightarrow0,\quad n\rightarrow+\infty.
\end{align*}
Therefore, we get
\begin{align*}
\theta'(0)\int_X[G,h_m]_{H}\varphi d\nu=0,
\end{align*}
for every $\varphi\in C_b^1(X)$. By recalling that $\theta'(0)\neq0$, it follows that $[G,h_m]_{H}=0$ for $\nu$-a.e. in $X$ for every $m\in\N$, which gives $G=0$ for $\nu$-a.e. in $X$.

The second part of the statement follows from standard arguments.
\end{proof}

\begin{remark}
\label{rmk:U_sob_space_nu}
If $(U_n)\subseteq \fcon_b^\infty(X)$ converges to $U$ in $W^{1,p}(X,\mu)$ for some $p>1$, then for every $q\in[1,p)$ we get
\begin{align*}
& \int_X|U_n-U|^qd\nu
= \leq \|U_n-U\|_{\elle^p(X,\mu)}^q\|e^{-U}\|_{\elle^{p/(p-q)}}\rightarrow 0, \\
& \int_X|D_HU_n-D_H^\mu U|_H^qd\nu
= \leq \|D_HU_n-D_H^\mu U\|_{\elle^p(X,\mu;H)}^q\|e^{-U}\|_{\elle^{p/(p-q)}}\rightarrow 0,
\end{align*}
as $n\to +\infty$. This means that $U\in W^{1,p}(X,\nu)$ for every $p\in[1,+\infty)$ and $D_HU=D_H^\mu U$. 
\end{remark}

The following technical lemma will be used in the sequel. 
\begin{lemma}
\label{lem:der_modulo}
Let $p\in[1,+\infty)$ and let $f\in W^{1,p}(X,\nu)$. Then, $|f|, f^+,f^-\in W^{1,p}(X,\nu)$ and $D_{H}|f|={\rm  sgn}(f)D_{H}f$, $D_{H}f^+=D_{H}f\chi_{\{f>0\}}$ and $D_{H}f^-=D_{H}f\chi_{\{f<0\}}$. Further, $D_Hf=0$ for $\nu$-a.e. in $f^{-1}(0)$. The same results hold true if we replace $\nu$ with $\mu$ and $D_{H}$ with $D^\mu_{H}$.
\end{lemma}
\begin{proof}
The proof can be obtain by repeating verbatim that of \cite[Lemma 2.7]{DPL14}.
\end{proof}

Thanks to Lemma \ref{lem:der_modulo} we are able to prove that a bounded function $g\in W^{1,p}(X,\nu)$ with $D_{H}g=0$ is constant for $\nu$-a.e. in $X$.
\begin{proposition}
\label{prop:gradiente_nullo_funz_costante}
Let $g\in W^{1,p}(X,\nu)$, with $p>1$, be a bounded function such that $D_{H}g=0$. Then, $g$ equals a constant for $\nu$-a.e. in $X$.
\end{proposition}
\begin{proof}
We will show that $g\in W^{1,1}(X,\mu)$ and $D_{H}^\mu g=0$. From \cite[Example 5.4.16]{Bog98} and the equivalence of the measures  $\mu$ and $\nu$ the thesis follows. Let $(g_n)\subseteq \fcon_b^1(X)$ be such that $g_n\rightarrow g$ in $W^{1,p}(X,\nu)$ as $n\rightarrow+\infty$. It follows that
\begin{align*}
D_{H}^\mu(g_ne^{-U})
= & e^{-U}(D_{H}^\mu g_n)-g_ne^{-U}(D_H^\mu U)
= e^{-U}(D_{H} g_n)-g_ne^{-U}(D_H^\mu U) \\
& \rightarrow  e^{-U}(D_{H} g)-ge^{-U}(D_H^\mu U)=-ge^{-U}(D_H^\mu U), \quad n\rightarrow+\infty,
\end{align*}
in $\elle^q(X,\mu)$ for every $1<q<p$. Hence, $ge^{-U}\in W^{1,q}(X,\mu)$ and $D_{H}^\mu(ge^{-U})= -ge^{-U}(D_H^\mu U)$. Let us set $U_n:=U\wedge n$. For every $n\in\N$ the function $U_n$ is bounded, from Lemma \ref{lem:der_modulo} we infer that $e^{U_n}\in W^{1,r}(X,\mu)$ for every $r\in[1,+\infty)$ and
$D_{H}^\mu (e^{U_n})=e^{U_n}(D_H^\mu U)\chi_{\{U<n\}}$. If we consider $r=q'$ the conjugate exponent of $q$, we have $ge^{-U}e^{U_n}=g\chi_{\{U\leq n\}}+ge^{n-U}\chi_{\{U> n\}}\in W^{1,1}(X,\mu)$ and
\begin{align*}
D_{H}^\mu(ge^{-U}e^{U_n})
= & e^{U_n}(D_{H}^\mu(ge^{-U}))+ge^{-U}(D_{H}^\mu e^{U_n}) \\
= & -ge^{U_n-U}(D_H^\mu U)+ge^{U_n-U}(D_H^\mu U)\chi_{\{U<n\}} \\
= & -ge^{U_n-U}(D_H^\mu U)\chi_{\{U>n\}}\\
=& -ge^{n-U}(D_H^\mu U)\chi_{\{U>n\}}.
\end{align*}
We notice that $(ge^{-U}e^{U_n})=(g\chi_{\{U\leq n\}}+ge^{n-U}\chi_{\{U> n\}})$ is a sequence of bounded functions which pointwise converges to $g$ as $n\rightarrow+\infty$, and $ge^{n-U}(D_H^\mu U)\chi_{\{U>n\}}\rightarrow0$ in $\elle^r(X,\mu;H)$ for every $r\in[1,+\infty)$. Hence, $g\in W^{1,1}(X,\mu)$ and $D_{H}^\mu g=0$, which gives that $g$ is constant $\mu$-a.e. in $X$.
\end{proof}

Now we define the $\nu$-divergence operator as the adjoint of $D_H$ in $\elle^2(X,\nu)$.
\begin{defn}
\label{def:H-divergence}
We set
\begin{align}
D(D_H^*)& :=\left\{f\in \elle^2(X,\nu;H):\exists g\in \elle^2(X,\nu), \int_X[f,D_Hv]_Hd\nu
=-\int_X gvd\nu, \ \forall v\in\fcon_b^1(X)\right\}, \notag\\
D_H^*f& :=g.
\label{divergence_operator}
\end{align}
We say that the operator $D_H^*:D(D_H^*)\subseteq \elle^2(X,\nu;H)\rightarrow \elle^2(X,\nu)$ is the $\nu$-divergence operator.
\end{defn}

Since $D(D_H)=W^{1,2}(X,\nu)$ is densely defined in $\elle^2(X,\nu)$ and $D_H$ is a closed operator in $\elle^2(X,\nu)$, from \cite[Theorem 13.12]{Ru91} it follows that $D(D_H^*)$ is densely defined in $\elle^2(X,\nu;H)$ and $(D_H^*)^*=D_H$.

Let us consider again the operator $D^k_H$ introduced in \eqref{k-H-der}. The next proposition states that the operator $(D_{H},D_{H}^2,\ldots,D_{H}^k)$ is closable in $\elle^p(X,\nu)$.

\begin{proposition}
\label{prop:closability_k_gradient}
Let $k\geq3$, $k\in\N$. The operator $(D_{H},\ldots,D^k_{H}):\fcon_b^\infty(X)\rightarrow \elle^p(X,\nu;H)\times\ldots\times  \elle^p(X,\nu;\mathcal H_k(H))$ is closable in $\elle^p(X,\nu)$ for every $p\in[1,+\infty)$. We denote by $(D_{H},\ldots,D_{H}^k)$ its closure and by $W^{k,p}(X,\nu)$ the domain of its closure. $W^{k,p}(X,\nu)$ is a Banach space if endowed with the norm
\begin{align*}
\|f\|_{k,p}:=\left(\|f\|_{k-1,p}^p+\|D^k_{H}f\|_{\elle^p(X,\nu;\mathcal H_k(H))}^p\right)^{1/p}, \quad f\in W^{k,p}(X,\nu).
\end{align*}
Finally, $W^{k,2}(X,\nu)$ is a Hilbert space if endowed with inner product
\begin{align*}
\langle f,g\rangle_{W^{k,2}(X,\nu)}
= \langle f,g\rangle_{W^{k-1,2}(X,\nu)}
+ \langle D^k_{H}f,D^k_{H} g\rangle_{\elle^2(X,\nu;\mathcal H_k(H))}, \quad f,g\in W^{k,2}(X,\nu).
\end{align*}
\end{proposition}
\begin{proof}
Let us sketch the proof when $k=2$, since the other cases can be analogously treated. Let $(F_n)\subseteq \fcon_b^\infty(X)$ be such that $F_n\rightarrow 0$ in $L^p(X,\nu)$ and $(D_HF_n,D_H^2F_n)\rightarrow (G_1,G_2)$ in $L^p(X,\nu;H)\times L^p(X,\nu;\mathcal H_2(H))$ as $n\rightarrow+\infty$. From Proposition \ref{prop:closability_gradient} we already know that $G_1=0$ for $\nu$-a.e. in $X$. Let $\{h_n:n\in\N\}$ be an orthonormal basis of $H$ and let $g\in C_b^1(X)$. For every $j,m\in\N$ we have
\begin{align*}
\int_Xg[D^2_HF_n h_m,h_j]_H d\nu
= & \int_X\left(\hat h_j g+[D_HU,h_j]g-[D_Hg,h_j]_H\right)[D_HF_n,h_m]_Hd\nu\rightarrow 0, \quad n\rightarrow+\infty.
\end{align*} 
The conclusion follows from standard arguments.
\end{proof}
\begin{remark}
Let us apply \eqref{int_by_parts} to the function $D^k_Hf$ with $f\in \fcon_b^\infty(X)$ as in \eqref{k-H-der}. Let $h_1,\ldots,h_k\in H$. We have
\begin{align*}
\left[D_H\left(D_H^{k-1}f(h_1,\ldots,h_{k-1})\right)(x),h_k\right]_H
=&  \left[D_H^kf(x)(h_1,\ldots,h_{k-1}),h_k\right]_H \\
= & D_H^kf(x)(h_1,\ldots,h_k), \quad x\in X,
\end{align*} 
where
\begin{align*}
& D_H^kf(x)(h_1,\ldots,h_{k-1}) \\
& =\sum_{i_1,\ldots,i_k=1}^n\frac{\partial^k\varphi}{\partial\xi_{i_1}\cdots \partial \xi_{i_k}}(\langle x,x_1^*\rangle,\ldots,\langle x,x_n^*\rangle)[i_\infty^*x^*_{i_1},h_1]\cdots[i_{\infty}^*x^*_{i_{k-1}},h_{k-1}]_Hi^*_\infty x^*_{i_k}, \quad x\in X,
\end{align*}
is understood as an element of $H$. Formula \eqref{int_by_parts} gives
\begin{align}
\int_XD_H^kf(x)(h_1,\ldots,h_k)d\nu
= & \int_X\left[D_H\left(D_H^{k-1}f(h_1,\ldots,h_{k-1})\right),h_k\right]_Hd\nu \notag \\
= & \int_XD_H^{k-1}f(h_1,\ldots,h_{k-1})(\hat{h_k}+[D_HU,h_k]_H)d\nu.
\label{int_parti_k}
\end{align}
\end{remark}

\begin{remark}
One may ask if the operator $D_H^k:\fcon_b^\infty(X)\subseteq \elle^p(X,\nu)\to \elle^p(X,\nu;\mathcal H_k(H))$ is closable in $\elle^p(X,\nu)$ for $k\geq 2$. The answer is positive, providing we consider additional assumptions on $U$. Let us consider a sequence $(f_n)\subseteq \fcon_b^\infty(X)$ such that $f_n\rightarrow 0$ in $\elle^p(X,\nu)$ and $D_H^2f_n\to \Phi$ in $\elle^p(X,\nu;\mathcal H_k(H))$ as $n\rightarrow+\infty$. For every $h_1,h_2\in H$ and every $g\in \fcon_b^\infty(X)$ we get
\begin{align*}
\int_X[D^2_Hf_n h_1,h_2]_H gd\nu
= & \int_X[D_Hf_n,h_1](g\hat h_2+g[D_HU,h_2]_H-[D_Hg,h_2]_H)d\nu.
\end{align*}
Hence, we need a second integration by parts to exploit the convergence to $0$ of $(f_n)$. This means that the function $U$ should be smoother, namely it should belong to some Sobolev space $W^{2,q}(X,\nu)$ for some $q\in(1,+\infty)$. Since this condition will be assumed in Section \ref{sec:equivalence}, we postponed therein the proof of closability of $D_H^k$ by means of a different technique. 
\end{remark}

Now we define the $H$-gradient on vector-valued functions. Let $V$ be a separable Hilbert space. We introduce the operator $\overline{D_{H}}$ defined on smooth functions $F\in\fcon_b^\infty(X;V)$ as follows:
\begin{align*}
\overline{D_{H}}F(x)=\sum_{i=1}^mD_{H}f_i(x)\otimes v_i\in H\otimes V, \quad x\in X,
\end{align*}
where
\begin{align*}
F(x)=\sum_{i=1}^m f_i(x)v_i,\quad f_i\in \fcon_b^\infty(X), \ v_i\in V, \ i=1,\ldots,m.
\end{align*}

The arguments of Proposition \ref{prop:closability_gradient}, adapted to vector-valued functions, allow us to prove the closability of $\overline {D_{H}}$ in $\elle^p(X,\nu;V)$ with $p\in[1,+\infty)$.
\begin{proposition}
\label{prop:closability_vector_valued_gradient}
For every $p\in[1,+\infty)$ the operator $\overline{D_{H}}:\fcon_b^\infty(X;V)\rightarrow \elle^p(X,\nu; H\otimes V)$ is closable in $\elle^p(X,\nu;V)$. We still denote by $\overline{D_{H}}$ its closure and by $W^{1,p}(X,\nu;V)$ the domain of its closure. $W^{1,p}(X,\nu;V)$ is a Banach space if endowed with the norm
\begin{align*}
\|F\|_{1,p,V}:=\left(\|F\|_{\elle^p(X,\nu;V)}^p+\|\overline{D_{H}}F\|_{\elle^p(X,\nu; H \otimes V)}^p\right)^{1/p}, \quad F\in W^{1,p}(X,\nu;V).
\end{align*}
Finally, $W^{1,2}(X,\nu;V)$ is a Hilbert space if endowed with the inner product
\begin{align*}
\langle F,G\rangle_{W^{1,2}(X,\nu;V)}
= \langle F,G\rangle_{\elle^2(X,\nu;V)}
+ \langle \overline {D_{H}}F,\overline{D_{H}}G\rangle_{\elle^2(X,\nu;H\otimes V)}, \quad F,G\in W^{1,2}(X,\nu;V).
\end{align*}
\end{proposition}

We conclude this subsection by showing that for every $F\in W^{1,p}(X,\nu;V)$ with $p\in[2,+\infty)$, the dual function $F^*:=|F|^{p-2}_VF$ belongs to $W^{1,p'}(X,\nu;V)$, where $p'=\frac{p}{p-1}$.

\begin{lemma}
\label{lemma:der_duale}
Let $F\in W^{1,p}(X,\nu;V)$ with $p\in[2,+\infty)$. Then, the function $F^*:=|F|^{p-2}_VF\in W^{1,p'}(X,\nu;V)$ and
\begin{align}
\label{der_duale}
\overline {D_{H}}F^*
= (p-2)|F|^{p-4}_V[\overline {D_{H}}F,F]_V\otimes F+|F|^{p-2}_V\overline {D_{H}}F.
\end{align}
Here, $[\overline {D_H}F,F]_V$ is understood as an element of $H$.
\end{lemma}
\begin{proof}
The fact that $F^*\in \elle^{p'}(X,\nu;V)$ is trivial. Let us consider a sequence $(F_n)\subseteq \fcon_b^\infty(X;V)$ such that $F_n\rightarrow F$ in $W^{1,p}(X,\nu;V)$ as $n\rightarrow+\infty$. For every $n\in\N$ we set $F_n^*:=|F_n|^{p-2}_VF_n$. 
Hence, there exists a subsequence $(F_{k_n})\subseteq (F_n)$ such that
\begin{align}
\label{conv_n_k_dual_funz}
F^*_{k_n}(x)\rightarrow F^*(x), \ F_{k_n}(x)\rightarrow F(x), \ \overline{D_H}F_{k_n}(x)\rightarrow \overline{D_H}F(x), \quad \nu\textup{-a.e. }x\in X.
\end{align}
Further,
\begin{align}
\label{der_duale_appr}
\overline {D_{H}}F_{k_n}^*
= (p-2)|F_{k_n}|^{p-4}_V[\overline {D_{H}}F_{k_n},F_{k_n}]_V\otimes F_{k_n}+|F_{k_n}|^{p-2}_V\overline {D_{H}}F_{k_n}, \quad n\in\N,
\end{align}
where the writing $[\overline {D_{H}}F_{k_n},F_{k_n}]_V$ is meant as an element of $H$. From \eqref{conv_n_k_dual_funz} it follows that 
\begin{align*}
\overline{D_H}F^*_{k_n}(x)\rightarrow
(p-2)|F(x)|^{p-4}_V[\overline {D_{H}}F(x),F(x)]_V\otimes F(x)+|F(x)|^{p-2}_V\overline {D_{H}}F(x)=:\Phi(x), \ \  \nu{\textup{-a.e. }}x\in X,
\end{align*}
as $n\rightarrow+\infty$. We claim that $F^*_{k_n}\rightarrow F^*$ in $\elle^{p'}(X,\nu;V)$ and that $\overline{D_H}F^*_{k_n}\rightarrow \Phi$ in $\elle^{p'}(X,\nu;H\otimes V)$. If the claim is true, the fact that $\overline{D_H}$ is a closed operator in $\elle^{p'}(X,\nu;V)$ implies that $F^*\in W^{1,p'}(X,\nu;V)$ and $\overline {D_H}F^*=\Phi$, which gives the thesis. Hence, it remain to prove the claim. Egoroff's Theorem (see \cite[Theorem 2.2.1]{Bog07}, whose proof can be easily generalized to vector-valued functions) implies that for every $\varepsilon>0$ there exists a Borel set $X_\varepsilon\subseteq X$ such that $\nu( X_\varepsilon)<\varepsilon$ and $F^*_{k_n}\rightarrow F^*$ uniformly on $X\setminus X_\varepsilon$. Let $\varepsilon>0$. For every $n\in\N$ we have
\begin{align}
\label{funz_duale_1}
\int_X|F^*_{k_n}-F^*|_V^{p'}d\nu
= \int_{X\setminus X_\varepsilon}|F^*_{k_n}-F^*|_V^{p'}d\nu+\int_{X_{\varepsilon}}|F^*_{k_n}-F^*|_V^{p'}d\nu.
\end{align}
Since $F^*_{k_n}\rightarrow F^*$ uniformly on $X\setminus X_\varepsilon$, it follows that
\begin{align}
\label{funz_duale_2}
 \int_{X\setminus X_\varepsilon}|F^*_{k_n}-F^*|_V^{p'}d\nu\rightarrow0, \quad n\rightarrow+\infty.
\end{align}
Further, from the definition of $F^*_{k_n}$ and of $F^*$ it follows that
\begin{align}
\int_{X_{\varepsilon}}|F^*_{k_n}-F^*|_V^{p'}d\nu
= & \|\chi_{X_\varepsilon}(F^*_{k_n}-F^*)\|_{\elle^{p'}(X,\nu;V)}^{p'} \notag\\
\leq & \left(\|\chi_{X_\varepsilon}F^*_{k_n}\|_{\elle^{p'}(X,\nu;V)}+\|\chi_{X_\varepsilon}F^*\|_{\elle^{p'}(X,\nu;V)}\right)^{p'} \notag\\
= & \left(\|\chi_{X_\varepsilon}F_{k_n}\|_{\elle^{p}(X,\nu;V)}^{p+1}+\|\chi_{X_\varepsilon}F\|_{\elle^{p}(X,\nu;V)}^{p-1}\right)^{p'} .
\label{funz_duale_3}
\end{align}
From \eqref{funz_duale_1}, \eqref{funz_duale_2} and \eqref{funz_duale_3}, and recalling that $F_{k_n}\rightarrow F$ in $\elle^p(X,\nu;V)$ as $n\rightarrow+\infty$, it follows that
\begin{align*}
\limsup_{n\rightarrow+\infty}
\int_X|F^*_{k_n}-F^*|_V^{p'}d\nu\leq 2\|\chi_{X_\varepsilon}F\|_{\elle^{p}(X,\nu;V)}^{p'},
\end{align*}
for every $\varepsilon>0$. Since $\nu(X_\varepsilon)<\varepsilon$ for every $\varepsilon>0$, we conclude that $F^*_{k_n}\rightarrow F^*$ in $\elle^{p'}(X,\nu;V)$ as $n\rightarrow+\infty$.

Let us consider $\overline{D_H}F^*_{k_n}$. Again from Egoroff's Theorem, for every $\varepsilon>0$ there exists a Borel set $X_\varepsilon\subseteq X$ such that $\nu( X_\varepsilon)<\varepsilon$ and 
$\overline{D_H}F^*_{k_n}\rightarrow \Phi$ uniformly on $X\setminus X_\varepsilon$. Let us fix $\varepsilon>0$. For every $n\in\N$ we have
\begin{align*}
\int_X|\overline{D_H}F^*_{k_n}-\Phi|^{p'}_{H\otimes V}d\nu
= & \int_{X\setminus X_\varepsilon}|\overline{D_H}F^*_{k_n}-\Phi|^{p'}_{H\otimes V}d\nu
+\int_{X_\varepsilon}|\overline{D_H}F^*_{k_n}-\Phi|^{p'}_{H\otimes V}d\nu
=:I_1^n+I_2^n, \ \  n\in\N.
\end{align*}
Since $\overline{D_H}F^*_{k_n}\rightarrow \Phi$ uniformly on $X\setminus X_\varepsilon$, it follows that $I_1^n\rightarrow 0$ as $n\rightarrow+\infty$. Let us estimate $I_2^n$. We have
\begin{align}
I_2^n
= & \|\chi_{X_\varepsilon}(\overline{D_H}F^*_{k_n}-\Phi)\|_{\elle^{p'}(X,\nu;H\otimes V)}^{p'}\notag \\
\leq & \left(\|\chi_{X_\varepsilon}\overline{D_H}F^*_{k_n}\|_{\elle^{p'}(X,\nu;H\otimes V)}+\|\chi_{X_\varepsilon}\Phi\|_{\elle^{p'}(X,\nu;H\otimes V)}\right)^{p'} \notag\\
=: & (J_1^n+J_2)^{p'}, \quad n\in\N.
\label{der_duale_primo}
\end{align}
Taking \eqref{der_duale_appr} into account we infer that
\begin{align}
(J_1^n)^{p'}
\leq & (p-1)^{p'}\int_X\chi_{X_\varepsilon}|F_{k_n}|^{p'(p-2)}_V|\overline{D_H}F_{k_n}|_{H\otimes V}^{p'}d\nu \notag \\
\leq & (p-1)^{p'}\left(\int_X\chi_{X_\varepsilon}|F_{k_n}|_V^{p}d\nu\right)^{\frac{p-2}{p-1}}\left(\int_X\chi_{X_\varepsilon}|\overline{D_{H}}F_{k_n}|^{p}_{H\otimes V}d\nu\right)^{\frac{1}{p-1}} \notag \\
= & (p-1)^{p'}\|\chi_{X_\varepsilon}F_{k_n}\|_{\elle^p(X,\nu;V)}^{p(p-2)/(p-1)}\cdot
 \|\chi_{X_\varepsilon}\overline{D_H}F_{k_n}\|_{\elle^p(X,\nu;H\otimes V)}^{p/(p-1)},
\label{der_duale_1}
\end{align}
where in the second inequality we have applied the H\"older's inequality with $q=p-1$ and $q'=\frac{p-1}{p-2}$, and we have used the fact that $p'=\frac{p}{p-1}$. 
Recalling that $F_n\rightarrow F$ in $W^{1,p}(X,\nu;V)$ as $n\rightarrow+\infty$, we deduce that
\begin{align}
\label{stima_der_duale_1}
\limsup_{n\rightarrow +\infty}J^n_1
\leq (p-1)\|\chi_{X_\varepsilon}F\|_{\elle^p(X,\nu;V)}^{p-2}
\cdot  \|\chi_{X_\varepsilon}\overline{D_H}F\|_{\elle^p(X,\nu;H\otimes V)}.
\end{align}
As far as $J_2$ is considered, arguing as in \eqref{der_duale_1} with $F_{k_n}$ replaced by $F$ it follows that
\begin{align}
\label{stima_der_duale_2}
J_2
\leq (p-1)\|\chi_{X_\varepsilon}F\|_{\elle^p(X,\nu;V)}^{p-2}\cdot
 \|\chi_{X_\varepsilon}\overline{D_H}F\|_{\elle^p(X,\nu;H\otimes V)}.
\end{align}
From \eqref{der_duale_primo}, \eqref{stima_der_duale_1} and \eqref{stima_der_duale_2}, it follows that
\begin{align*}
\limsup_{n\rightarrow+\infty}\int_X|\overline {D_H}F^*_{k_n}-\Phi|_{H\otimes V}^{p'}d\nu
\leq (2p-2)^{p'}\|\chi_{X_\varepsilon}F\|_{\elle^p(X,\nu;V)}^{p'(p-2)}\cdot
 \|\chi_{X_\varepsilon}\overline{D_H}F\|_{\elle^p(X,\nu;H\otimes V)}^{p'},
\end{align*}
for every $\varepsilon>0$. Since $\nu(X_\varepsilon)< \varepsilon$ for every $\varepsilon>0$, we get
\begin{align*}
\limsup_{n\rightarrow+\infty}\int_X|\overline {D_H}F^*_{k_n}-\Phi|_{H\otimes V}^{p'}d\nu=0,
\end{align*}
which proves the claim.
\end{proof}

\subsection{The perturbed Ornstein-Uhlenbeck operator in \texorpdfstring{$\elle^p(X,\nu)$}{LpXm}}
We introduce the symmetric bilinear form
\begin{align*}
\mathcal E(u,v):=\int_X[D_H^\mu U,D_{H}v]_{H}d\nu, \quad u,v\in W^{1,2}(X,\nu),
\end{align*}
with domain $\mathcal D=W^{1,2}(X,\nu)$. From Proposition \ref{prop:closability_gradient} it follows that $\mathcal E$ is a symmetric bilinear form which satisfies the strong sector condition, hence it is closed and coercive. \cite[Theorem 2.8 \& Corollary 2.10]{MR92} imply that the operator $(L_2,D(L_2))$ defined as
\begin{align*}
D(L_2) & :=\left\{u\in W^{1,2}(X,\nu):\exists g\in \elle^2(X,\nu), \ \mathcal E(u,v)=-\int_Xgvd\nu, \ \forall v\in\fcon_b^\infty(X)\right\} ,\vspace{2mm}\\
L_2u& :=g,
\end{align*}
is the infinitesimal generator of an analytic symmetric strongly continuous semigroup of contractions $(T_2(t))_{t\geq0}$ on $\elle^2(X,\nu)$. 

\begin{remark}
The integration by parts formula \eqref{int_by_parts} allows us to provide an explicit expression to $L$ applied to smooth functions.
Indeed, we have
$\fcon_b^2(X)\subseteq D(L_2)$ and for every $u\in \fcon_b^2(X)$ it follows that
\begin{align}
\label{esplicit_L2}
(L_2u)(x)={\rm Tr}[D^2_{H}u(x)]_{H}-\langle x,Du(x)\rangle-[D_H^\mu U(x),D_H u(x)]_{H}, \quad x\in X.
\end{align}    
\end{remark}


We provide some properties of $(T_2(t))_{t\geq0}$ which arise from the theory of Dirichlet forms. 
For reader's convenience, we recall the definition of Dirichlet forms, Dirichlet operators and sub-Markovian semigroups, and their main properties (see e.g. \cite[Chapter 1, Definitions 4.1 \& 4.5, Proposition 4.3 \& Theorem 4.4]{MR92}).
\begin{defn}
\label{defn:dirichlet_form}
Let $(E,B,\mu)$ be a measure space and let $\mathscr H:=\elle^2(E,\mu)$ be a Hilbert space.
\begin{itemize}
\item [(i)] A symmetric closed coercive form $(\mathcal E,D(\mathcal E))$ on $\mathscr H$ is called a Dirichlet form if for every $u\in D(\mathcal E)$ one has $u^+\wedge 1\in D(E)$ and $\mathcal E(u^+\wedge 1,u^+\wedge 1)\leq \mathcal E(u,u)$.
\item[(ii)] A semigroup $(S(t))_{t\geq0}$ on $\mathscr H$ is called sub-Markovian if for every $t\geq0$ and every $f\in \mathscr H$ with $0\leq f\leq 1$ $\mu$-a.e. in $E$, we have $0\leq S(t)f\leq 1$ $\mu$-a.e. in $E$.
\item[(iii)] A closed linear densely defined operator $A$ on $\mathscr H$ is called Dirichlet operator on $\mathscr H$ if
\begin{align*}
\int_E Au(u-1)^+d\mu\leq 0, \quad u\in D(A).
\end{align*}
\end{itemize}
\end{defn}

\begin{proposition}
\label{prop:dirichlet-submark}
Let $(\mathcal E,D(\mathcal E))$ be a symmetric closed coercive Dirichlet form on $\elle^2(E,\mu)$, let $A$ be the operator associated to $\mathcal E$ and let $(S(t))_{t\geq0}$ be the strongly continuous contraction semigroup on $\elle^2(E,\mu)$ generated by $A$. Then, the following are equivalent:
\begin{itemize}
\item [(i)] $(\mathcal E,D(\mathcal E))$ is a Dirichlet form on $\elle^2(E,\mu)$.
\item[(ii)] $(S(t))_{t\geq0}$ is a sub-Markovian semigroup on $\elle^2(E,\mu)$.
\item[(iii)] $\mathcal A$ is a Dirichlet operator on $\elle^2(E,\mu)$.
\end{itemize}
\end{proposition}

\begin{proposition}\label{prop:quasi_regular}
The bilinear form $(\mathcal E, W^{1,2}(X,\nu))$ is s Dirichlet form on $\elle^2(X,\nu)$. Then:
\begin{enumerate}[(i)]
\item $(T_2(t))_{t\geq0}$ is non-negative, i.e., for every $t>0$ we have $T(t)f\geq0$ $\nu$-a.e. in $X$ for every $f\in\elle^2(X,\nu)$ such that $f\geq0$ $\nu$-a.e. in $X$. Further, if $f\in \elle^2(X,\nu)$ and has positive infimum, i.e., there exists a positive constant $c$ such that $f\geq c$ $\nu$-a.e. in $X$, then $T_2(t)f\geq c$ $\nu$-a.e. in $X$ for every $t>0$.
\item For every $f\in C_b(X)$ we have $|T_2(t)f|\leq T_2(t)|f|$ $\nu$-a.e. in $X$ and $\|T_2(t)f\|_\infty\leq \|f\|_\infty$ for every $t\geq0$.
\item for every $1\leq p\leq q\leq +\infty$ we have
\begin{align}
\label{jensen_semigroup}
\big(T_2(t)(|f|^p)\big)^{1/p}\leq\big (T_2(t)(|f|^q)\big)^{1/q}, \quad \nu\textup{-a.e. in }X, \ f\in C_b(X), \  t\geq0.
\end{align}
\item For every $p\in(1,+\infty)$ we have
\begin{align}
\label{holder_semigroup}
|T_2(t)(fg)|\leq \big(T_2(t)(|f|^p)\big)^{1/p} \big(T_2(t)(|g|^{p'})\big)^{1/p'}, \quad  \nu\textup{-a.e. in }X, \   f,g\in C_b(X), \  t\geq0,
\end{align}
where $p'$ is the conjugate exponent of $p$, i.e., $p'=\frac{p}{p-1}$.
\item For every $p\in[1,+\infty)$ we have
\begin{align}
\label{mink_inequality_somma}
\big(T_2(t)(|f+g|^p)\big)^{1/p}
\leq \big(T_2(t)(|f|^p)\big)^{1/p}+\big(T_2(t)|g|^p)\big)^{1/p}, \quad  \nu\textup{-a.e. in }X, \   f,g\in C_b(X), \  t\geq0.
\end{align}
\end{enumerate}
\end{proposition}
\begin{proof}
At first we show that $(\mathcal E,W^{1,2}(X,\nu))$ is a Dirichlet form. From Definition \ref{defn:dirichlet_form} it is enough to show that for every $u\in W^{1,2}(X,\nu)$ we have $u^+\wedge 1\in W^{1,2}(X,\nu)$ and
\begin{align*}
\mathcal E(u^+\wedge 1,u^+\wedge 1)\leq \mathcal E(u,u), \quad  u\in W^{1,2}(X,\nu).
\end{align*}
From Lemma \ref{lem:der_modulo} we infer that $u^+\wedge 1\in W^{1,2}(X,\nu)$ for every $u\in W^{1,2}(X,\nu)$ and $D_{H}(u^+\wedge 1)=D_H^\mu u\chi_{\{u\in(0,1)\}}$. Then,
\begin{align*}
\mathcal E(u^+\wedge 1,u^+\wedge 1)
= \int_{\{u\in(0,1)\}}|D_H^\mu u|^2_{H}d\nu
\leq \int_X|D_H^\mu u|^2_{H}d\nu
=\mathcal E(u,u),
\end{align*}
which implies that $(\mathcal E,W^{1,2}(X,\nu))$ is a Dirichlet form. From Proposition \ref{prop:dirichlet-submark} it follows that $(T_2(t))_{t\geq0}$ is a sub-Markovian semigroup on $\elle^2(X,\nu)$.

\vspace{2mm}
Let us prove $(i)$. Let $f\in L^\infty(X,\nu)$ be such that $f\geq0$ $\nu$-a.e. in $X$, and let $t>0$. Then, $0\leq f\|f\|_\infty^{-1}\leq 1$ $\nu$-a.e. in $X$, which gives
\begin{align*}
0\leq T_2(t)(f\|f\|^{-1}_{\infty})\leq 1, \quad \nu\textup{-a.e. in } X,
\end{align*}
and this implies that $0\leq T_2(t)f\leq \|f\|_\infty$ $\nu$-a.e. in $X$. If $f\in\elle^2(X,\nu)$ satisfies $f\geq0$ $\nu$-a.e. in $X$, we consider the sequence $(f_n:=f\wedge n)\subseteq \elle^\infty(X,\nu)$. Since $f_n\geq0$ $\nu$-a.e. in $X$, it follows that $T_2(t)f_n\geq0$ $\nu$-a.e. in $X$. Moreover, $(T_2(t)f_n)$ converges to $T_2(t)f$ in $\elle^2(X,\nu)$, and, up to a subsequence, pointwise $\nu$-a.e. in $X$. hence, $T_2(t)f\geq0$ $\nu$-a.e. in $X$. To prove the second part, we notice that if $f\geq c>0$ $\nu$-a.e. in $X$ for some constant $c$, then the function $g:=f-c\geq0$ $\nu$-a.e. in $X$. By recalling that $T_2(t)a=a$ for every $a\in\R$, it follows that
\begin{align*}
T_2(t)f-c
= & T_2(t)f-T_2(t)c
= T_2(t)(f-c)\geq0,
\end{align*}
which gives the thesis.

\vspace{2mm}
Now we prove $(ii)$. Let $f\in C_b(X)$ and let us consider $T_2(t)f^+$ and $T_2(t)f^-$. From $(i)$ it follows that $T_2(t)f^+,T_2(t)f^-\geq0$ $\nu$-a.e. in $X$. Then,
\begin{align*}
|T_2(t)f|
= |T_2(t)(f^+-f^-)|
\leq |T_2(t)f^+|+|T_2(t)f^-|
= T_2(t)(f^++f^-)
=T_2(t)|f|.
\end{align*}
In the proof of $(i)$ we have shown that $0\leq T_2(t)|f|\leq \|f\|_\infty$. Hence, we get $|T_2(t)f|\leq \|f\|_\infty$, which gives the second part of $(ii)$.

To prove $(iii)$ and $(iv)$ it is enough to repeat the computations in \cite[Lemma 2.1]{Sh97} which hold for self-adjoint sub-Markovian semigroups.

Finally, $(v)$ follows arguing as in the proof of the classical Minkowski inequality.
\end{proof}

$(T_2(t))_{t\geq0}$ extends to a positive  strongly continuous semigroup of contractions $(T_p(t))_{t\geq0}$ on $\elle^p(X,\nu)$ for every $p\in[1,+\infty)$. We state this result in the following proposition. The proof of this fact can be obtain in \cite[Theorem 1.4.1]{Dav89} (see also \cite[Proposition 3.7]{AD20}).
\begin{proposition}
$(T_2(t))_{t\geq0}$ extends to a positive  strongly continuous semigroups of contractions $(T_p(t))_{t\geq0}$ on $\elle^p(X,\nu)$ with infinitesimal generator $L_p$ for every $p\in[1,+\infty)$. These semigroups are consistent in the sense that if $f\in \elle^p(X,\nu)$ for some $p\in[1,+\infty)$, then $T_p(t)f=T_q(t)f$ for every $q\in[1,p]$ and every $t\geq0$.
\end{proposition}

\begin{remark}
Where there is no danger of confusion we omit the subscript $p$ and we simply denote by $(T(t))_{t\geq0}$ and by $L$ the semigroup and its infinitesimal generator on $\elle^p(X,\nu)$, respectively.
\end{remark}

When $p\in[2,+\infty)$ we can associate a bilinear form to $L_p$.
\begin{corollary}
\label{coro:bilinear_p}
For every $p\in[2,+\infty)$ we have
\begin{align*}
\int_XL_pfgd\nu=-\int_X[D_{H}f, D_{H}g]_{H}d\nu,
\end{align*}
for every $f\in D(L_p)$ and every $g\in W^{1,p'}(X,\nu)$.
\end{corollary}
\begin{proof}
Since $D(L_p)\subseteq D(L_2)$, the thesis follows from the definition of $L_2$ and from the density of $\fcon_b^\infty(X)$ in $W^{1,p'}(X,\nu)$.
\end{proof}

In the following proposition we state a Minkowski's integral inequality for the semigroup $(T(t))_{t\geq0}$. The proof is postponed in the appendix for reader's convenience. This result will be crucial to extend gradient estimates for the scalar semigroup $(T(t))_{t\geq0}$ to the vector-valued one.
\begin{proposition}
\label{prop:mink_ineq}
Let $(M,\mathcal M,\gamma)$ be a measure space and let $\gamma$ be a non-negative $\sigma$-finite measure. Then, for every $q\in[1,+\infty)$ and for every $F:X\times M\rightarrow \R$ measurable function with respect to the $\sigma$-field $B(X)\otimes \mathcal M$ we have
\begin{align}
\label{minkowski_smgr}
\left(\int_M\left|T(t)(F(\cdot,y))(x)\right|^q\gamma(dy)\right)^{1/q}
\leq T(t)\left(\int_M|F(\cdot,y)|^q\gamma(dy)\right)^{1/q}(x), \quad t\geq 0, \ x\in X.
\end{align}
We admit that both the sides in \eqref{minkowski_smgr} are $\infty$. Here, the notation $T(t)(G(\cdot,y))$, $y\in M$, means that we are applying the operator $T(t)$ to the function $x\mapsto G(x,y)$, where $G:X\times M\rightarrow \R$ is measurable with respect to the $\sigma$-field $B(X)\otimes \mathcal M$.
\end{proposition}

\subsection{The perturbed vector-valued Ornstein-Uhlenbeck semigroup}
Let $V$ be a separable Hilbert space. It is well-known that $\fcon_b^\infty(X;V)$ is dense in $\elle^p(X,\nu;V)$. We consider the vector-valued semigroup $(T_p^V(t))_{t\geq0}$ on $\elle^p(X,\nu;V)$, extension of $(T_p(t))_{t\geq0}$ (for vector-valued extensions of positive operators see \cite[Subsection 4.5.3]{Gr14}), defined on $F\in \elle^p(X,\nu;V)$ with finite range by
\begin{align*}
T_p^V(t)F:=\sum_{i=1}^n T_p(t)f_iv_i, \quad F:=\sum_{i=1}^n f_iv_i, \quad f_i\in\elle^p(X,\nu), \ v_i\in V, \ i=1,\ldots,n.
\end{align*}
It turns out that $(T_p^V(t))_{t\geq0}$ is a strongly continuous semigroup of contractions on $\elle^p(X,\nu;V)$. We denote by $L_p^V$ its infinitesimal generator, which acts on functions $F\in D(L_p^V)$ with finite range as
\begin{align*}
L_p^VF=\sum_{i=1}^nL_pf_iv_i, \quad F:=\sum_{i=1}^n f_iv_i, \quad f_i\in D(L_p), \ v_i\in V, \ i=1,\ldots,n,
\end{align*}
where we have supposed, without loss of generality, than $\{v_1,\ldots,v_n\}$ is an orthonormal system in $V$. Where there is no danger of confusion we omit the subscript $p$ and we simply denote by $(T^V(t))_{t\geq0}$ and by $L^V$ the semigroup and its infinitesimal generator on $\elle^p(X,\nu;V)$, respectively. We generalize some properties of the scalar semigroup and of its infinitesimal generator to $(T^V(t))_{t\geq0}$ and of $L^V$ in the following lemma.
\begin{lemma}
\label{lem:prop_vettoriali}
Let $(T^V(t))_{t\geq0}$ and $L^V$ be defined as above. Then:
\begin{enumerate}[(i)]
\item $(T^V(t))_{t\geq0}$ is self-adjoint in $\elle^2(X,\nu;V)$, i.e.,
\begin{align*}
\int_X[T^V(t)F,G]_Vd\nu=\int_X[F,T^V(t)G]_Vd\nu, \quad F\in \elle^p(X,\nu;V), \quad G\in \elle^{p'}(X,\nu;V).
\end{align*}
\item Let $p\in[2,+\infty)$. For every $F\in D(L_p^V)$ and $G\in W^{1,p'}(X,\nu;V)$ we have
\begin{align}
\label{vector_ipp}
\int_X[L_p^VF,G]_V d\nu
=\int_X[\overline{D_{H}}F, \overline{D_{H}}G]_{H\otimes V}d\nu.
\end{align}
\end{enumerate}
\end{lemma}
\begin{proof}
$(i)$ By density we can limit ourselves to prove the statement for $F,G\in \fcon_b^\infty(X;V)$. We fix
 \begin{align*}
F:=\sum_{i=1}^nf_iv_i, \quad G:=\sum_{j=1}^n g_jw_j, \quad f_i,g_j\in \fcon_b^\infty(X), \ v_i,w_j\in V, \ i,j=1,\ldots,n.
\end{align*}
Then,
\begin{align*}
\int_X[T^V(t)F,G]_Vd\nu
= &\sum_{i,j=1}^n[v_i,w_j]_V\int_XT(t)f_ig_jd\nu
= \sum_{i,j=1}^n[v_i,w_j]_V\int_Xf_iT(t)g_jd\nu\\
= & \int_X[F,T^V(t)G]_Vd\nu
\end{align*}
where we have used the property of symmetry of the scalar semigroup $(T(t))_{t\geq0}$.

\vspace{2mm}
$(ii)$ By approximation we can limit ourselves to prove the statement for functions with finite range. Let us fix
 \begin{align*}
F:=\sum_{i=1}^nf_iv_i, \quad G:=\sum_{j=1}^n g_jw_j, \quad f_i\in D(L_p), \ g_j\in W^{1,p'}(X,\nu), \ v_i,w_j\in V, \ i,j=1,\ldots,n,
\end{align*}
where we have supposed, without loss of generality, than $\{v_1,\ldots,v_n\}$ and $\{w_1,\ldots,w_n\}$ are orthonormal systems in $V$. From Corollary \ref{coro:bilinear_p} we infer that
\begin{align*}
\int_X[L_p^VF,G]_Vd\nu
= & \sum_{i,j=1}^n[v_i,w_j]_{V}\int_XL_pf_ig_jd\nu
= -\sum_{i,j=1}^n[v_i,w_j]_{V}\int_X[D_{H}f_i,D_{H}g_j]_{H}d\nu\\
= & -\int_X[\overline {D_{H}}F,\overline {D_{H}}G]_{H\otimes V}d\nu.
\end{align*}
\end{proof}

\section{Analysis of \texorpdfstring{$(T(t))_{t\geq0}$}{T(t)t} and of \texorpdfstring{$(T^V(t))_{t\geq0}$}{TV(t)t}}
\label{sec:analysis}
\subsection{Asymptotic behaviour and pointwise gradient estimates to \texorpdfstring{$(T(t))_{\geq0}$}{T(t)t}}
In this subsection we show that for every $p\in[1,+\infty)$ the semigroup $(T(t))_{t\geq0}$ is ergodic in $\elle^p(X,\nu)$, i.e., for every  $f\in \elle^p(X,\nu)$ we have $T(t)f\rightarrow \nu(f)$ in $\elle^p(X,\nu)$ as $t\rightarrow+\infty$, where
\begin{align}
\label{media_f}
\nu(f):=\int_Xfd\nu, \quad f\in \elle^p(X,\nu).
\end{align}
To prove this fact we use the following intermediate result, whose proof is inspired by those in \cite[Corollary 3.6 \& Proposition 3.7]{DPG01} in finite dimension.
\begin{lemma}
\label{lemma:conv_infinito_L2}
Let $f\in \elle^2(X,\nu)$. Then,
\begin{align}
\label{asy_est_D_HT(t)}
\lim_{t\rightarrow+\infty}\|D_{H}T(t)f\|_{\elle^2(X,\nu;H)}^2=0.
\end{align}
\end{lemma}
\begin{proof}
We claim that for every  $f\in D(L_2)$, the function $t\mapsto \chi_f(t):=\|D_{H}T(t)f\|_{\elle^2(X,\nu;H)}^2\in L^1(0,+\infty)$. To prove the claim, let us fix $f\in D(L_2)$ and $t>0$. We have
\begin{align*}
\frac{d}{dt}\int_X|T(t)f|^2d\nu
= & 2\int_X(T(t)f)(L_2T(t)f)d\nu
=-2\int_X|D_H T(t)f|^2_Hd\nu.
\end{align*}
Integrating between $0$ and $t$ we get
\begin{align*}
\|T(t)f\|_{\elle^2(X,\nu)}^2-\|f\|_{\elle^2(X,\nu)}^2
= -2\int_0^t\|D_{H}T(s)f\|_{\elle^2(X,\nu;H)}^2ds, \quad t\geq0.
\end{align*}
This implies that
\begin{align*}
\|T(t)f\|_{\elle^2(X,\nu)}^2+2\int_0^t\|D_{H}T(s)f\|_{\elle^2(X,\nu;H)}^2ds
\leq \|f\|_{\elle^2(X,\nu)}^2, \quad t\geq0,
\end{align*}
and the claim is so proved. Let us consider $f\in D(L_2^2)$, where $L_2^2$ is the square power of the operator $L_2$. This means that both $\chi_f$ and $\chi_{Lf}$ belong to $\elle^1(0,+\infty)$. Since
\begin{align*}
\left|\frac{d}{dt}\chi_f(t)\right|=2\left|\int_X[D_{H}T(t)L_2f,D_{H}T(t)f]_{H}d\nu\right|
\leq \chi_f(t)+\chi_{Lf}(t),
\end{align*}
it follows that both $\chi_f$ and $\chi_f'$ belong to $\elle^1(0,+\infty)$. This implies that $\chi_f\in W^{1,1}(0,+\infty)$, and therefore
\begin{align*}
\lim_{t\rightarrow+\infty}\chi_f(t)=0, \quad f\in D(L_2^2).
\end{align*}
Since $(T(t))_{t\geq0}$ is an analytic semigroup in $\elle^2(X,\nu)$, it follows that $T(1)f\in D(L_2^n)$ for every  $n\in\N$. Then,
\begin{align*}
\lim_{t\rightarrow+\infty}\|D_{H}T(t)f\|_{\elle^2(X,\nu;H)}^2
= \lim_{t\rightarrow+\infty}\|D_{H}T(t-1)T(1)f\|_{\elle^2(X,\nu;H)}^2=0, \quad f\in \elle^2(X,\nu).
\end{align*}
\end{proof}

\begin{proposition}
\label{prop:comp_asintotico_fz}
For every  $p\in[1,+\infty)$ and every  $f\in \elle^p(X,\nu)$ we have
\begin{align}
\label{convergenza_Lp_asymp}
\lim_{t\rightarrow+\infty}\|T(t)f-\nu(f)\|_{\elle^p(X,\nu)}=0,
\end{align}
where $\nu(f)$ has been defined in \eqref{media_f}.
\end{proposition}
\begin{proof}
Let us split the proof into three steps. 
In the former we show that for every  $f\in C_b(X)$ the function $T(t)f$ weakly converges to $\nu(f)$ in $\elle^2(X,\nu)$ as $t\rightarrow+\infty$, in the second we prove \eqref{convergenza_Lp_asymp} for $f\in C_b(X)$, in the latter we conclude.

\vspace{2mm}
{\bf{STEP $1$}}. Let $f\in C_b(X)$. Since $(T(t))_{t\geq0}$ is a semigroup of contractions in $\elle^2(X,\nu)$, it follows that there exists a sequence $(t_n)$ diverging to $+\infty$ and a function $g\in\elle^2(X,\nu)$ such that $T(t_n)f\rightarrow g$ weakly in $\elle^2(X,\nu)$ as $n\rightarrow+\infty$. Further, $g$ is bounded and $\|g\|_\infty\leq \|f\|_\infty$. Indeed, for every  positive, bounded and continuous function $v$, from Proposition \ref{prop:quasi_regular}$(ii)$ we have
\begin{align*}
\left |\int_X T(t_n)fvd\nu\right|\leq \|f\|_\infty\int_Xvd\nu.
\end{align*}
Letting $n\rightarrow+\infty$ we get
\begin{align*}
\left |\int_X gvd\nu\right|\leq \|f\|_\infty\int_Xvd\nu.
\end{align*}
The arbitrariness of $v$ implies that $-\|f\|_\infty\leq g\leq \|f\|_\infty$ for $\nu$-a.e. in $X$.

Let us consider $u\in D(D_{H}^*)$. We have
\begin{align*}
\int_XgD_{H}^*ud\nu
= & \lim_{n\rightarrow+\infty}\int_XT(t_n)f  D_{H}^*ud\nu
= \lim_{n\rightarrow+\infty}\int_X[D_{H}T(t_n)f, u]_{H}d\nu=0,
\end{align*}
where the last equality follows from Lemma \ref{lemma:conv_infinito_L2}. Therefore, $g\in D((D_{H}^*)^*)=W^{1,2}(X,\nu)$ and $D_{H} g=0$. From Proposition \ref{prop:gradiente_nullo_funz_costante} it follows that $g$ is constant $\nu$-a.e. in $X$. Finally, we have
\begin{align*}
g= & \int_Xgd\nu=\lim_{n\rightarrow+\infty}\int_X T(t_n)fd\nu=\lim_{n\rightarrow+\infty}\int_Xfd\nu=\nu(f),
\end{align*}
where the third equality follows from the fact that $\nu$ is an invariant measure for $(T(t))_{t\geq0}$. In particular, the above arguments show that for every  sequence $(t_n)\subseteq (0,+\infty)$ diverging to $+\infty$ as $n\rightarrow+\infty$ there exists a subsequence $(t_{k_n})\subseteq (t_n)$ such that $T(t_{k_n})f\rightarrow \nu(f)$ weakly in $\elle^2(X,\nu)$ as $n\rightarrow+\infty$. Hence, $T(t)f\rightarrow \nu(f)$ weakly in $\elle^2(X,\nu)$ as $t\rightarrow+\infty$.

\vspace{2mm}
{\bf {STEP $2$}}.
From Step $1$ we know that $T(t)f\rightarrow \nu(f)$ weakly in $\elle^2(X,\nu)$ for every  $f\in C_b(X)$. Then,
\begin{align}
\label{conv_L2_smgr}
\|T(t)f\|_{\elle^2(X,\nu)}^2
=\int_X(T(t)f)(T(t)f)d\nu
= \int_X(T(2t)f) fd\nu\rightarrow \int_X\nu(f)fd\nu
= \|\nu(f)\|_{\elle^2(X,\nu)}^2,
\end{align}
as $t\rightarrow+\infty$. Here, we have used the self-adjointness of $(T(t))_{t\geq0}$ in $\elle^2(X,\nu)$ and the semigroup property of $(T(t))_{t\geq0}$. \eqref{conv_L2_smgr} implies that $T(t)f\rightarrow \nu(f)$ in $\elle^2(X,\nu)$ as $t\rightarrow+\infty$. Therefore, for every  sequence $(t_n)\subseteq (0,+\infty)$ diverging to $+\infty$ as $n\rightarrow+\infty$ there exists a subsequence $(t_{k_n})\subseteq (t_n)$ such that $T(t_{k_n})f(x)\rightarrow \nu(f)$ for $\nu$-a.e. $x\in X$. By the dominated convergence theorem it follows that $T(t_{k_n})f\rightarrow \nu(f)$ in $\elle^p(X,\nu)$ as $n\rightarrow+\infty$, for every  $p\in[1,+\infty)$. This means that $T(t)f\rightarrow \nu(f)$ in $\elle^p(X,\nu)$ as $t\rightarrow+\infty$ for every  $p\in[1,+\infty)$.

\vspace{2mm}
{\bf STEP $3$}. Let $f\in\elle^p(X,\nu)$ and let $(f_n)\subseteq C_b(X)$ be such that $f_n\rightarrow f$ in $\elle^p(X,\nu)$ as $n\rightarrow+\infty$. From H\"older's inequality we get
\begin{align*}
\|\nu(f_n)-\nu(f)\|_{\elle^p(X,\nu)} 
\leq \|f_n-f\|_{\elle^p(X,\nu)}, \quad n\in\N.
\end{align*}
Hence,
\begin{align}
\|T(t)f-\nu(f)\|_{\elle^p(X,\nu)}
\leq & \|T(t)f-T(t)f_n\|_{\elle^p(X,\nu)}
+\|T(t)f_n-\nu(f_n)\|_{\elle^p(X,\nu)} \notag\\
& +\|\nu(f_n)-\nu(f)\|_{\elle^p(X,\nu)} \notag \\
\leq & 2\|f-f_n\|_{\elle^p(X,\nu)}+\|T(t)f_n-\nu(f_n)\|_{\elle^p(X,\nu)},
\label{stima_T(t)f-f}
\end{align}
where we have used the fact that $(T(t))_{t\geq0}$ is a semigroup of contractions in $\elle^p(X,\nu)$. Let us fix $\varepsilon>0$. Then, there exists $\overline n\in\N$ such that $\|f-f_{\overline n}\|_{\elle^p(X,\nu)}\leq \varepsilon/4$. From Step $2$ there exists $\overline t>0$ such that $\|T(t)f_{\overline n}-\nu(f_{\overline n})\|_{\elle^p(X,\nu)}\leq \varepsilon/2$ for every  $t>\overline t$. This means that, by choosing $n=\overline n$ in \eqref{stima_T(t)f-f}, for every  $\varepsilon>0$ there exists $\overline t>0$ such that $\|T(t)f-\nu(f)\|_{\elle^p(X,\nu)}\leq \varepsilon$ for every  $t>\overline t$. This gives the thesis.
\end{proof}

Here we show two different estimates of the $D_HT(t)f$. In Proposition \ref{prop:pointwise_stime_grad} we estimate $|D_{H}T(t)f|_H^p$ by means of $T(t)|D_{H}f|_H^p$, in Proposition \ref{prop:stime_gradiente_funzione} we estimate $|D_{H}T(t)f|_H^p$ by means of $T(t)|f|^p$. The proofs are inspired by those of \cite [Theorems 3.1 \& 3.3]{AngFerPal18}, do not present every  new idea and are quite technical, hence we postponed them in Appendix \ref{app_B}. 
\begin{proposition}
\label{prop:pointwise_stime_grad}
For every  $p\in[1,+\infty)$ and every  $f\in W^{1,p}(X,\nu)$ we have
\begin{align}
\label{stima_grad_grad_tesi}
|D_{H}T(t)f(x)|^p_{H}
\leq e^{-pt}(T^H(t)|D_{H}f|^p_{H})(x), \quad t\geq0, \ \nu{\textup{-a.e. }}x\in X.
\end{align}
\end{proposition}

\begin{proposition}
\label{prop:stime_gradiente_funzione}
Let $p\in(1,+\infty)$. Then, there exists a positive constant $c_p$, which only depends on $p$, blows up as $p\rightarrow1^+$ and equals $2^{-p/2}$ for $[2,+\infty)$, such that for every  $f\in\elle^p(X,\nu)$ we have 
\begin{align}
\label{stima_grad_funz_prima_parte_finale}
|D_{H}T(t)f(x)|^p_{H}\leq c_pt^{-p/2}T(t)|f|^p(x), \quad t>0, \ \nu\textup{-a.e.} \ x\in X.
\end{align}
\end{proposition}

\subsection{Asymptotic behaviour and gradient estimates to the vector-valued semigroup \texorpdfstring{$(T^V(t))_{t\geq0}$}{TV(t)t}}
Let $V$ be a separable Hilbert space. In this subsection we show some properties that the $V$-valued semigroup $(T^V(t))_{t\geq0}$ inherits from the scalar semigroup $(T(t))_{t\geq0}$. We collect such results in a unique proposition, whose proof is postponed in Appendix \ref{app_B}.

\begin{proposition} 
\label{prop:vett_smgr_prop}
\begin{enumerate}[(i)]
\item 
For every $p\in[1,+\infty)$ we have
\begin{align}
\label{convergenza_as_vett}
\nu(F):=\int_XFd\nu=\elle^p-\lim_{t\rightarrow+\infty}T^V(t)F, \quad F\in \elle^p(X,\nu;V).
\end{align}
\item 
For every $p\in[1,+\infty)$, every $t\geq0$ and every $F\in W^{1,p}(X,\nu;V)$ we have 
\begin{align}
\label{vector_stime_grad_grad}
|\overline{D_{H}}T^V(t)F(x)|_{H\otimes V}^p
\leq e^{-pt}T(t)\left(|\overline{D_{H}}F|_{H\otimes V}^p\right)(x),\qquad \nu\textup{-a.e.} \ x\in X.
\end{align}
As a byproduct,
\begin{align}
\label{vector_stime_grad_grad_2}
\int_X|\overline {D_{H}}T^V(t)F|_{H\otimes V}^pd\nu
\leq {e^{-pt}}\int_X|\overline {D_{H}}F|_{H\otimes V}^pd\nu, \quad t>0, \ F\in  W^{1,p}(X,\nu;V).
\end{align}
\item 
For every $p\in(1,+\infty)$, every $t>0$ and every $F\in \elle^p(X,\nu)$ we have
\begin{align}
\label{vector_stima_grad_funz}
|\overline{D_{H}}T^V(t)F(x)|_{H\otimes V}^p
\leq \frac{c_p}{t^{p/2}}T(t)\left(|F|_{V}^p\right)(x), \qquad  \nu\textup{-a.e.} \ x\in X,
\end{align}
where $c_p$ is the constant in \eqref{stima_grad_funz_prima_parte_finale}.
As a byproduct,
\begin{align}
\label{int_vector_stima_grad_funz}
\int_X|\overline {D_{H}}T^V(t)F|_{H\otimes V}^pd\nu
\leq \frac{c_p}{t^{p/2}}\int_X|F|_{ V}^pd\nu, \quad t>0, \ F\in \elle^p(X,\nu;V).
\end{align}
\item 
Let $p\in(1,+\infty)$. Then, 
\begin{align}
\label{stima_grad_funz_mista}
\|\overline {D_{H}}T^V(t)F\|_{\elle^p(X,\nu;H\otimes V)}
\leq (c_p)^{1/p}\max\{t^{-1/2},1\}\min\{1,e^{-t+1}\}\|F\|_{\elle^p(X,\nu;V)}, \quad t>0,
\end{align}
for every $F\in \elle^p(X,\nu;V)$, where $c_p$ is the constant in \eqref{stima_grad_funz_prima_parte_finale}. In particular, if $p\in[2,+\infty)$ we infer that
\begin{align}
\label{stima_grad_funz_mista_2+infty}
\|\overline {D_{H}}T^V(t)F\|_{\elle^p(X,\nu;H\otimes V)}
\leq \frac{\max\{t^{-1/2},1\}}{\sqrt2}\min\{1,e^{-t+1}\}\|F\|_{\elle^p(X,\nu;V)}, \quad t>0,
\end{align}
for every $F\in \elle^p(X,\nu;V)$.
\end{enumerate}
\end{proposition}

\section{Vector-valued Poincar\'e inequality}
\label{sec:poincare}

\subsection{Poincar\'e inequality}
Thanks to Proposition \ref{prop:vett_smgr_prop} we are able to prove a vector-valued Poincar\'e inequality with explicit constant.
\begin{thm}
\label{thm:vector_poincare}
Let $V$ be a separable Hilbert space and let $p\in[1,+\infty)$. Then, there exists a positive constant $k_p$ defined by
\begin{align}
\label{kp_constant}
k_p=
\begin{cases}
\sqrt{(p-1)}, & p\in[2,+\infty), \\
\displaystyle \frac{3}{\sqrt2},  & p\in[1,2),
\end{cases}
\end{align}
such that for every $F\in W^{1,p}(X,\nu;V)$ we have
\begin{align}
\label{vector_poincare}
\|F-\nu(F)\|_{\elle^p(X,\nu;V)}\leq k_p\|\overline {D_{H}}F\|_{\elle^p(X,\nu;H\otimes V)},
\end{align}
where $\nu(F)$ has been defined in \eqref{convergenza_as_vett}.
\end{thm}
\begin{remark}
We notice that, for $p\geq2$, the constant $k_p$ equals the constant in \cite{AdMuRo21} when $U=0$.    
\end{remark}
\begin{proof}
We split the proof into two parts. In the former we prove the thesis for $p\in[2,+\infty)$, in the latter we consider the remaining cases $p\in[1,2)$. As usual, by density it is enough to prove the statement for $F\in\fcon_b^\infty(X;V)$.

\vspace{2mm}
{\bf STEP $1$}. Let $p\in[2,+\infty)$ and let $F\in \fcon_b^\infty(X;V)$ 
with 
\begin{align*}
F:=\sum_{i=1}^mf_iv_i, \quad f_i\in \fcon_b^\infty(X), \ v_i\in V, \ i=1,\ldots,m,
\end{align*}
where $\{v_1,\ldots,v_m\}$ are orthonormal vectors in $V$. Then,
\begin{align*}
T^V(t)F=\sum_{i=1}^mT(t)f_iv_i,
\end{align*}
belongs to $D(L^V_p)$ and $\displaystyle \frac{d}{dt}T^V(t)F=L_p^VT^V(t)F$ for every $t\geq0$. We define $G:=F-\nu(F)\in\fcon_b^\infty(X;V)$. From Lemma \ref{lemma:der_duale} we infer that $G^*:=|G|_V^{p-2}G\in W^{1,p'}(X,\nu;V)$ and
\begin{align}
\label{der_G_vett_poincare}
\overline {D_{H}}G^*
= (p-2)|G|^{p-4}_V[\overline {D_{H}}G,G]_V\otimes G+|G|^{p-2}_V\overline {D_{H}}G,
\end{align}
where $[\overline {D_{H}}G,G]_V$ is meant as an element of $H$. Formula \eqref{convergenza_as_vett} gives
\begin{align}
\int_X|G|_V^pd\nu
= & \int_X[G,G^*]_Vd\nu
= \int_X[F-\nu(F),G^*]_Vd\nu \notag\\
= & \lim_{t\rightarrow+\infty}\int_X[T^V(0)F-T^V(t)F,G^*]_Vd\nu.
\label{poincare_vett_1}
\end{align}
Let us consider the argument of the limit. We get
\begin{align}
\int_X[T^V(0)F-T^V(t)F,G^*]_Vd\nu
= & -\int_X\left(\int_0^t\left[\frac{d}{ds}T^V(s)F,G^*\right]_Vds\right)d\nu \notag\\
= & -\int_X\left(\int_0^t\left[L_p^VT^V(s)F,G^*\right]_Vds\right)d\nu \notag \\
= &- \int_0^t\left(\int_X\left[L_p^VT^V(s)F,G^*\right]_Vd\nu \right)ds \notag \\
= & \int_0^t\left(\int_X\left[\overline{D_{H}}T^V(s)F,\overline{D_{H}}G^*\right]_{H\otimes V}d\nu \right)ds,
\label{poincare_vett_2}
\end{align}
where we have used the Fubini's theorem and Lemma \ref{lem:prop_vettoriali}$(ii)$. From \eqref{der_G_vett_poincare} we have
\begin{align}
\left|\int_X\left[\overline{D_{H}}T^V(s)F,\overline{D_{H}}G^*\right]_{H\otimes V}d\nu\right|
\leq & \int_X|\overline{D_{H}}T^V(s)F|_{H\otimes V}|\overline{D_{H}}G^*|_{H\otimes V}d\nu \notag \\
\leq  & (p-1)\int_X |\overline{D_{H}}T^V(s)F|_{H\otimes V}|\overline{D_{H}}F|_{H\otimes V}|G|_V^{p-2}d\nu, 
\label{poincare_vett_3}
\end{align}
since $\overline {D_{H}}G=\overline{D_{H}}F$. Let us consider $p>2$. From the generalized H\"older's inequality for three functions with exponents $p,p$ and $\frac{p}{p-2}$ and \eqref{vector_stime_grad_grad_2} it follows that
\begin{align}
\int_X&  |\overline{D_{H}}T^V(s)F|_{H\otimes V}|\overline{D_{H}}F|_{H\otimes V}|G|_V^{p-2}d\nu \notag  \\
\leq & \left(\int_X |\overline{D_{H}}T^V(s)F|_{H\otimes V}^{p}d\nu\right)^{1/p}
\left(\int_X |\overline{D_{H}}F|_{H\otimes V}^{p}d\nu\right)^{1/p}
\left(\int_X|G|_V^{p}d\nu\right)^{1-2/p} \notag \\
\leq & e^{-s} \left(\int_X |\overline{D_{H}}F|_{H\otimes V}^{p}d\nu\right)^{2/p}
\left(\int_X|G|_V^{p}d\nu\right)^{1-2/p}.
\label{poincare_vett_4}
\end{align}
Putting together \eqref{poincare_vett_1}-\eqref{poincare_vett_4} we infer that
\begin{align*}
\int_X|G|^pd\nu
\leq & (p-1)\left(\int_0^\infty e^{-s}ds\right) \left(\int_X|\overline{D_{H}}F|_V^{p}d\nu\right)^{2/p}\left(\int_X|G|_V^{p}d\nu\right)^{1-2/p} \\
= & (p-1) \left(\int_X|\overline{D_{H}}F|_V^{p}d\nu\right)^{2/p}\left(\int_X|G|_V^{p}d\nu\right)^{1-2/p}. 
\end{align*}
Dividing both the sides by $\left(\int_X|G|_V^{p}d\nu\right)^{1-2/p}$ we get
\begin{align*}
\left(\int_X|G|^pd\nu\right)^{2/p}
\leq  (p-1)\left(\int_X|\overline{D_{H}}F|_V^{p}d\nu\right)^{2/p},
\end{align*}
which gives the thesis with $k_p=\left(p-1\right)^{1/2}$. If $p=2$, then $G^*=G$ and $\overline {D_H}G^*=\overline{D_H}F$, and from \eqref{vector_stime_grad_grad_2} we get
\begin{align}
\left|\int_X\left[\overline{D_{H}}T^V(s)F,\overline{D_{H}}G^*\right]_{H\otimes V}d\nu\right|
= & \left|\int_X\left[\overline{D_{H}}T^V(s)F,\overline{D_{H}}F\right]_{H\otimes V}d\nu\right| \notag \\
\leq &  \left(\int_X |\overline{D_{H}}T^V(s)F|_{H\otimes V}^{2}d\nu\right)^{1/2}
\left(\int_X |\overline{D_{H}}F|_{H\otimes V}^{2}d\nu\right)^{1/2} \notag \\
\leq & e^{-s} \left(\int_X |\overline{D_{H}}F|_{H\otimes V}^{2}d\nu\right),
\label{vector_poincare_p=2}
\end{align}
for every $s>0$. By collecting \eqref{poincare_vett_1}, \eqref{poincare_vett_2} and \eqref{vector_poincare_p=2} we get the thesis with $c_2=1=(2-1)^{1/2}$.

\vspace{2mm}
{\bf STEP $2$}.
Let $p\in(1,2)$ and let $F\in\fcon_b^\infty(X;V)$. Then, we have
\begin{align}
\label{poincare_vett_stima_duale_1}
\|F-\nu(F)\|_{\elle^p(X,\nu;V)}
= & \sup_{
\begin{array}{l}
G\in\fcon_b^\infty(X;V), \\ 
\|G\|_{\elle^{p'}(X,\nu;V)}\leq 1
\end{array}}\int_X[F-\nu(F),G]_Vd\nu.
\end{align}
Let $G\in \fcon_b^\infty(X;V)$ with $\|G\|_{\elle^{p'}(X,\nu;V)}\leq 1$. Then,
\begin{align}
\label{poincare_vett_stima_duale_2}
\int_X[F-\nu(F),G]_Vd\nu
=\lim_{t\rightarrow +\infty}\int_X[T^V(0)F-T^V(t)F,G]_Vd\nu
= \lim_{t\rightarrow +\infty}\int_X[F,T^V(0)G-T^V(t)G]_Vd\nu,
\end{align}
where we have applied Lemma \ref{lem:prop_vettoriali}$(i)$ and \eqref{convergenza_as_vett}. Arguing as in Step $1$ we infer that
\begin{align}
\int_X[F,T^V(0)G-T^V(t)G]_Vd\nu
= & \int_0^t\left(\int_X[\overline {D_{H}}F,\overline{D_{H}}T^V(s)G]_{H\otimes V}d\nu\right)ds \notag \\
\leq &  \int_0^t\left(\int_X|\overline {D_{H}}F|_{H\otimes V}|\overline{D_{H}}T^V(s)G|_{H\otimes V}d\nu\right)ds \notag\\
\leq &  \left(\int_X|\overline {D_{H}}F|_{H\otimes V}^pd\nu\right)^{1/p}\!\!
\int_0^t\left(\int_X|\overline{D_{H}}T^V(s)G|^{p'}_{H\otimes V}d\nu\right)^{1/p'}\!\!ds,
\label{poincare_vett_stima_duale_3}
\end{align}
with $p'>2$. By applying \eqref{stima_grad_funz_mista_2+infty} to $\int_X|\overline{D_{H}}T^V(s)G|^{p'}_{H\otimes V}d\nu$ and recalling that $\|G\|_{\elle^{p'}(X,\nu;V)}\leq 1$ we infer that
\begin{align}
\left(\int_X|\overline{D_{H}}T^V(s)G|^{p'}_{H\otimes V}d\nu\right)^{1/p'}
\leq & \frac{\max\{s^{-1/2},1\}}{\sqrt2}\min\{1,e^{-s+1}\}\|G\|_{\elle^{p'}(X,\nu;V)} \notag\\
\leq& \frac{\max\{s^{-1/2},1\}}{\sqrt2}\min\{1,e^{-s+1}\},
\label{poincare_vett_stima_duale_4}
\end{align}
which gives
\begin{align}
\int_0^t\left(\int_X|\overline{D_{H}}T^V(s)G|^{p'}_{H\otimes V}d\nu\right)^{1/p'}\!\!ds
\leq &\frac{1}{\sqrt2} \left(\int_0^1s^{-1/2}ds+\int_1^te^{-s+1}ds\right)=\frac{3-e^{-t+1}}{\sqrt2},
\label{poincare_vett_stima_duale_5}
\end{align}
for every $t>1$. Collecting \eqref{poincare_vett_stima_duale_2}-\eqref{poincare_vett_stima_duale_5} we get
\begin{align*}
\int_X[F-\nu(F),G]_Vd\nu
\leq \frac{3}{\sqrt2}\|\overline{D_{H}}F\|_{\elle^p(X,\nu;H\otimes V)},
\end{align*}
for every $G\in \fcon_b^\infty(X;V)$ with $\|G\|_{\elle^{p'}(X,\nu;V)}\leq 1$. From \eqref{poincare_vett_stima_duale_1} we get
\begin{align}
\label{poincare_vett_duale_finale}
\|F-\nu(F)\|_{\elle^p(X,\nu;V)}
\leq \frac{3}{\sqrt2}\|\overline{D_{H}}F\|_{\elle^p(X,\nu;H\otimes V)},
\end{align}
for every $p\in(1,2)$ and every $F\in \fcon_b^\infty(X;V)$. Since the constant $k_p=\frac{3}{\sqrt2}$ does not depend on $p$, letting $p\rightarrow 1^+$ we get \eqref{poincare_vett_duale_finale} also for $p=1$.
\end{proof}

As a byproduct of Theorem \ref{thm:vector_poincare} we get the following result.
\begin{corollary}
Let $p\in[1,+\infty)$ and let $F\in W^{k+1,p}(X,\nu;V)$ with $k\in\N$. Then,
\begin{align}
\label{poincare_vett_derivate}
\|\overline D_{H}^kF-\nu(\overline D_{H}^kF)\|_{\elle^p(X,\nu;H^{\otimes k}\otimes V)}
\leq k_p\|\overline D_{H}^{k+1}F\|_{\elle^p(X,\nu; H^{\otimes k+1}\otimes V)},
\end{align}
where $k_p$ is the positive constant in \eqref{kp_constant}.
\end{corollary}
\begin{proof}
The thesis follows by applying Theorem \ref{thm:vector_poincare} with $V=H^{\otimes k}\otimes V$ and $F=\overline D_{H}^kF$.
\end{proof}

\section{Equivalent norms on \texorpdfstring{$W^{k,p}(X,\nu)$}{WkpXm}}
\label{sec:equivalence}
The combination of Wiener chaos decomposition (see Subsection \ref{sub:wiener_chaos}) and of vector-valued Poincar\'e inequality \eqref{poincare_vett_derivate} allows us to prove that $D_H^k:\fcon_b^\infty(X)\subseteq \elle^p(X,\nu)\to \elle^p(X,\nu;\mathcal H_k(H))$  is closable in $\elle^p(X,\nu)$ and the norm $\|\cdot\|_{k,p}$ is equivalent to the graph norm of $D_{H}^k$ in $\elle^p(X,\nu)$ for every $p\in(1,+\infty)$ and every $k\in\N$, where we still denote by $D_H^k$ the closure of $D_H^k$ in $\elle^p(X,\nu)$. We split this section into two parts: in the former we consider the case $k=2$, which holds true when $U$ only satisfies Hypothesis \ref{hyp:U}, in the latter we consider the cases $k\geq 3$, when we need additional assumptions on $U$. For every $p\in(1,+\infty)$ and every $k\in\N$ we introduce the norm
\begin{align*}
\|f\|_{p,D_{H}^k}:=\left(\|f\|_{\elle^p(X,\nu)}^p+\|D_{H}^kf\|_{\elle^p(X,\nu;\mathcal H_k(H))}^p\right)^{1/p}, \quad f\in\fcon_b^\infty(X).
\end{align*}

\subsection{The case \texorpdfstring{$k=2$}{k=2}}
\label{sub:eq_norm_=2}
\begin{thm}
\label{thm:equivalece_space_2}
Let $p\in(1,+\infty)$. Then, there exists two positive constants $K_p$ and $\widetilde K_p$ such that $K_p\|f\|_{p,D_{H}^2}\leq \|f\|_{2,p}\leq \widetilde K_p\|f\|_{p,D^2_H}$ for every $f\in \fcon_{b}^\infty(X)$.
\end{thm}
\begin{proof}
From the definition it follows that $\|f\|_{p,D_{H}^2}\leq 2^{1-1/p}\|f\|_{2,p}$ for every $f\in \fcon_{b}^\infty(X)$. 

Let us prove the second inequality. We show that for every $p\in(1,+\infty)$ there exist a positive constant $\widetilde C_p$, which only depends on $p$ and $U$, such that for every $f\in \fcon_{b}^\infty(X)$ we have
\begin{align}
\label{controllo_1_derivative}
\|D_{H}f\|_{\elle^p(X,\nu;\mathcal H_{2}(H))}
\leq k_{p}\|D_{H}^2f\|_{\elle^p(X,\nu;\mathcal H_{2}(H))}
+\widetilde C_p\|f\|_{\elle^p(X,\nu)},
\end{align}
and $k_p$ is the constant in \eqref{poincare_vett_derivate}. Let $p\in(1,+\infty)$ and let $f\in \fcon_{b}^{\infty}(X)$. Then, 
\begin{align*}
\|D_{H}f\|_{\elle^p(X,\nu;H)}
\leq & \|D_{H}f-\nu(D_{H}f)\|_{\elle^p(X,\nu;H)}
+\|\nu(D_{H}f)\|_{H} \\
\leq & k_p\|D^2_{H}f\|_{\elle^p(X,\nu;\mathcal H_2(H))}
+ \left(\sum_{n\in\N}\left|\int_X[D_{H}f,e_n]_{H}d\nu\right|^2\right)^{1/2},
\end{align*}
where in the last inequality we have taken advantage from \eqref{poincare_vett_derivate} and $\{e_n:n\in\N\}$ is any orthonormal basis of $H$. Let us estimate the second addend above. Integrating by parts we get
\begin{align}
\notag
\sum_{n\in\N}\left|\int_X[D_{H}f,e_n]_{H}d\nu\right|^2
= & \sum_{n\in\N}\left|\int_Xf(\hat{e_n}+[D_H^\mu U,e_n]_{H})d\nu\right|^2 \\
\leq & 2\sum_{n\in\N}\left|\int_Xe^{-U}f\hat{e_n}d\mu\right|^2
+2\sum_{n\in\N}\left|\int_Xf[D_H^\mu U,e_n]_{H}d\nu\right|^2.
\label{stima_1_stima_D_D^2}
\end{align}
From Remarks \ref{rmk:chaos_ind_base} and \ref{rmk:dec_wiener_chaos} it follows that $\{\hat e_n:n\in\N\}$ is an orthonormal basis of the first Wiener chaos $E_1$ and that
\begin{align*}
\sum_{n\in\N}\left|\int_Xe^{-U}f\hat{e_n}d\mu\right|^2
= \|I_1(e^{-U}f)\|_{\elle^2(X,\mu)}^2.
\end{align*}
Let $q=\frac{p+1}2\in(1,p)$. From Lemma \ref{lemma:L2-Lp-norm} there exists a positive constant $a_p$, which only depends on $p$, such that
\begin{align*}
\|I_1(e^{-U}f)\|_{\elle^2(X,\mu)}^2\leq a_p\|e^{-U}f\|^2_{\elle^q(X,\mu)}.
\end{align*}
By applying the H\"older's inequality with $r=\frac{p}{q}$ and $r'=\frac{p}{p-q}$ we get
\begin{align}
\|e^{-U}f\|_{\elle^q(X,\mu)}^2
= & \left\|\left(e^{-\frac qpU}|f|^q\right)e^{\left(-q+\frac qp\right)U}\right\|_{\elle^1(X,\mu)}^{2/q}
\leq  \|fe^{-\frac1pU}\|_{\elle^p(X,\mu)}^{2}
\left(\int_Xe^{-r'(q-\frac qp)U}d\mu\right)^{2/(qr')} \notag \\
=: & C_p\|f\|_{\elle^p(X,\nu)}^2,
\label{stima_2_D_D^2}
\end{align}
where $C_p$ is a constant which only depends on $p$ and $U$. Further,
\begin{align}
\notag \sum_{n\in\N}\left|\int_Xf[D_H^\mu U,e_n]_{H}d\nu\right|^2
=& \left|\int_Xf\sum_{n\in\N}[D_H^\mu U,e_n]_He_nd\nu\right|_H^2 
= \left|\int_X f D_H^\mu U d\nu\right|_H^2 \\
\leq & \left(\int_X|f| |D_H^\mu U|_Hd\nu\right)^2
\leq \|f\|_{\elle^p(X,\nu)}^2\||D_H^\mu U|_H\|_{\elle^{p'}(X,\nu)}^2.
\label{stima_3_D_D^2}
\end{align}
From Remark \ref{rmk:U_sob_space_nu} we infer that $\||D_H^\mu U|_H\|_{\elle^{p'}(X,\nu)}^2<+\infty$. 

Collecting together \eqref{stima_1_stima_D_D^2}-\eqref{stima_3_D_D^2} we infer that
\begin{align*}
 \left(\sum_{n\in\N}\left|\int_X[D_{H}f,e_n]_{H}d\nu\right|^2\right)^{1/2}
 \leq & \sqrt2\left(C_p+\||D_H^\mu U|_{H}\|_{\elle^{p'}(X,\nu)}^2\right)^{1/2}\|f\|_{\elle^p(X,\nu)} \\
 =:& \widetilde C_p\|f\|_{\elle^p(X,\nu)},
\end{align*}
which gives the thesis.
\end{proof}

\begin{corollary}
\label{coro:eq_norme_2}
For every $p\in(1,+\infty)$ the operator $D_H^2:\fcon_b^\infty(X)\subseteq\elle^p(X,\nu)\to \elle^p(X,\nu;\mathcal H_2(H))$ is closable in $\elle^p(X,\nu)$ and the domain of its closure, endowed with the graph norm, coincides with $W^{2,p}(X,\nu)$.
\end{corollary}
\begin{proof}
Let us consider a sequence of functions $(f_n)\subseteq \fcon_b^\infty(X)$ such that $f_n\to 0$ in $\elle^p(X,\nu)$ and $D^2_Hf_n\rightarrow \Phi$ in $\elle^p(X,\nu;\mathcal H_2(H))$ as $n\rightarrow+\infty$. From \eqref{controllo_1_derivative} it follows that $(D_Hf_n)$ is a Cauchy sequence in $\elle^p(X,\nu;H)$ and so $f_n)$ converges in $W^{2,p}(X,\nu)$. Since $f_n\rightarrow 0$ in $\elle^p(X,\nu)$ it follows that $\Phi=0$, which gives the closability of $D_H^2$. From Theorem \eqref{thm:equivalece_space_2} the graph norm of $D^2_H$ and the $W^{2,p}(X,\nu)$-norm are equivalent, and the thesis follows.
\end{proof}

\subsection{The case \texorpdfstring{$k\geq3$}{k>3}}
\label{sub:eq_norm_>2}
Without any additional assumption on $U$ we are able to extend \eqref{controllo_1_derivative} to $D^k_H$ for every $k\in\N$, $k\geq2$.
\begin{proposition}
\label{prop:stima_norma_k-1}
Let $p\in(1,+\infty)$ and let $k\in\N$, $k\geq2$. Then, there exist a positive constant $\widetilde C_p$, which only depends on $p$ and $U$, such that for every $f\in W^{k+1,p}(X,\nu)$ we have
\begin{align}
\label{controllo_k_derivative}
\|D_{H}^{k}f\|_{\elle^p(X,\nu;\mathcal H_{k}(H))}
\leq k_{p}\|D_{H}^{k+1}f\|_{\elle^p(X,\nu;\mathcal H_{k+1}(H))}
+\widetilde C_p\|D_{H}^{k-1}f\|_{\elle^q(X,\nu;\mathcal H_{k-1}(H))},
\end{align}
where $k_p$ is the constant in \eqref{poincare_vett_derivate}. 
\end{proposition}
\begin{proof}
Let $k\in\N$, $k\geq 2$, and let $f\in\fcon_b^\infty(X)$. From \eqref{poincare_vett_derivate} we have  
\begin{align}
\label{stima_vett_totale}
\|D_{H}^{k}f\|_{\elle^p(X,\nu;\mathcal H_{k}(H))}
\leq & \|D_{H}^{k}f-\nu(D_{H}f^{k})\|_{\elle^p(X,\nu;\mathcal H_{k}(H))}
+\|\nu(D_{H}^{k}f)\|_{\mathcal H_{k}(H)}\notag  \\
\leq & k_p\|D^{k+1}_{H}f\|_{\elle^p(X,\nu;\mathcal H_{k+1}(H))} \notag \\
& + \left(\sum_{i_1,\ldots,i_{k}\in\N}\left|\int_XD^{k}_{H}f(e_{i_1},\ldots,e_{i_{k}})d\nu\right|^2\right)^{1/2},
\end{align}
where $\{e_n:n\in\N\}$ is any orthonormal basis of $H$. We notice that
\begin{align*}
D^{k}_{H}f(e_{i_1},\ldots,e_{i_{k}})=[D^{k}_{H}f(e_{i_1},\ldots,e_{i_{k-1}}),e_{i_{k}}]_{H}, \quad i_1,\ldots,i_{k}\in\N,
\end{align*}
where $D^k_Hf(e_{i_1},\ldots,e_{i_{k-1}})$ is understood as an element of $H$.
By applying formula \eqref{int_parti_k} to the last addend in \eqref{stima_vett_totale} we get
\begin{align}
\sum_{i_1,\ldots,i_{k}\in\N} & \left|\int_XD^{k}_{H}f(e_{i_1},\ldots,e_{i_{k}})d\nu\right|^2 \notag \\
=& \sum_{i_1,\ldots,i_{k}\in\N}\left|\int_XD^{k-1}_{H}f(e_{i_1},\ldots,e_{i_{k-1}})(\hat {e_{i_{k}}}+[D_H^\mu U,e_{i_{k}}]_{H})d\nu\right|^2 \notag \\
\leq & 2\sum_{i_1,\ldots,i_{k-1}\in\N}\sum_{i_k\in\N}\left|\int_XD^{k-1}_{H}f(e_{i_1},\ldots,e_{i_{k-1}})\hat {e_{i_{k}}}d\nu\right|^2\notag  \\
& +2\sum_{i_1,\ldots,i_{k-1}\in\N}\sum_{i_k\in\N}\left|\int_XD^{k-1}_{H}f(e_{i_1},\ldots,e_{i_{k-1}})[D_H^\mu U,e_{i_{k}}]_{H}d\nu\right|^2.
\label{stima_vett_totale_2}
\end{align}
By applying Lemma \ref{lemma:L2-Lp-norm} for every $i_1,\ldots,i_{k-1}\in\N$  we get
\begin{align*}
\sum_{i_k\in\N}\left|\int_XD^{k-1}_{H}f(e_{i_1},\ldots,e_{i_{k-1}})\hat {e_{i_{k}}}d\nu\right|^2
= & \|I_1(D^{k-1}_{H}f(e_{i_1},\ldots,e_{i_{k-1}})e^{-U})\|^2_{\elle^2(X,\mu)}\\
\leq & a_p \|D^{k-1}_{H}f(e_{i_1},\ldots,e_{i_{k-1}})e^{-U}\|^2_{\elle^q(X,\mu)},
\end{align*}
where $q=\min\left\{\frac{p+1}{2},\frac32\right\}$ and $a_p$ is a positive constant which only depends on $p$. Summing up $i_1,\ldots,i_{k-1}$ over $\N$ we get
\begin{align*}
& \sum_{i_1,\ldots,i_{k-1}\in\N}\sum_{i_k\in\N} \left|\int_XD^{k-1}_{H}f(e_{i_1},\ldots,e_{i_{k-1}})\hat {e_{i_{k}}}d\nu\right|^2 \\
\leq & a_p \sum_{i_1,\ldots,i_{k-1}\in\N}  \|D^{k-1}_{H}f(e_{i_1},\ldots,e_{i_{k-1}})e^{-U}\|^2_{\elle^q(X,\mu)} \\
= & a_p \sum_{i_1,\ldots,i_{k-1}\in\N} \left(\int_X|D_H^{k-1}f(e_{i_1},\ldots, e_{i_{k-1}})e^{-U}|^qd\mu\right)^{2/q}. 
\end{align*}
Let us apply the Minkowski's integral inequality with $\mu_1=\mu$, $\mu_2$ being the product of $k-1$ counting measures on $\N$ and $p=2/q>1$. It follows that
\begin{align}
\notag
& \sum_{i_1,\ldots,i_{k-1}\in\N} \left(\int_X|D_H^{k-1}f(e_{i_1},\ldots, e_{i_{k-1}})e^{-U}|^qd\mu\right)^{2/q} \\
= & \bigg(\Big(\sum_{i_1,\ldots,i_{k-1}\in\N} \Big(\int_X|D_H^{k-1}f(e_{i_1},\ldots, e_{i_{k-1}})e^{-U}|^qd\mu \Big)^{2/q}\Big)^{q/2}\bigg)^{2/q} \notag\\
\leq & \bigg(\int_X\Big(\sum_{i_1,\ldots,i_{k-1}\in\N}|D_H^{k-1}f(e_{i_1},\ldots,e_{i_{k-1}})e^{-U}|^2\Big)^{q/2} d\mu\bigg)^{2/q} \notag\\
= & \left\|\|D_H^{k-1}f\|_{\mathcal H_{k-1}(H)}e^{-U}\right\|_{\elle^q(X,\mu)}^2.
\label{stima_vett_minkowski}
\end{align}
By repeating the computations in \eqref{stima_2_D_D^2} it follows that
\begin{align}
\label{stima_vettoriale_totale_3}
\left\|\|D^{k-1}_{H}f\|_{\mathcal H_{k-1}(H)}e^{-U}\right\|^2_{\elle^q(X,\mu)}
\leq  & C_p\|D^{k-1}_{H}f\|_{\elle^p(X,\nu;\mathcal H_{k-1}(H))}^2,
\end{align}
for some positive constant $C_p$ which only depends on $p$ and $U$. Moreover,
\begin{align}
\notag
\sum_{i_1,\ldots,i_{k-1}\in\N}\sum_{i_k\in\N}& \left|\int_XD^{k-1}_{H}f(e_{i_1},\ldots,e_{i_{k-1}})[D_H^\mu U,e_{i_{k}}]_{H}d\nu\right|^2 \\
\leq &\left(\int_X\|D^{k-1}_{H}f\|_{\mathcal H_{k-1}(H)}|D_H^\mu U|_Hd\nu\right)^2 \notag \\
\leq & \|D_H^{k-1}f\|_{\elle^p(X,\nu;\mathcal H_{k-1}(H))}^2\||D_H^\mu U|_H\|_{\elle^{p'}(X,\nu)}^2,
\label{stima_finale_2_k>1}
\end{align}
for every $p\in(1,+\infty)$, and from Remark \ref{rmk:U_sob_space_nu} it follows that $\||D_H^\mu U|_H\|_{\elle^{p'}(X,\nu)}<+\infty$. Formulae \eqref{stima_vett_totale}-\eqref{stima_finale_2_k>1} imply that there exists a positive constant $\widetilde C_p$, which only depends on $p$ and $U$, such that
\begin{align*}
\left(\sum_{i_1,\ldots,i_{k}\in\N}\left|\int_XD^{k}_{H}f(e_{i_1},\ldots,e_{i_{k}})d\nu\right|^2\right)^{1/2}\leq \widetilde C_p \|D_H^{k-1}f\|_{\elle^p(X,\nu;\mathcal H_{k-1}(H))}.
\end{align*}
By density we conclude.
\end{proof}

A generalization of Theorem \ref{thm:equivalece_space_2} is available for $k\geq 3$ under additional assumptions on $U$. Indeed, the idea is going on with integration by parts of the second addend in \eqref{stima_vett_totale}
\begin{align*}
\sum_{i_1,\ldots,i_{k}\in\N} & \left|\int_XD^{k}_{H}f(e_{i_1},\ldots,e_{i_{k}})d\nu\right|^2,
\end{align*}
in order to estimate this term by means of $\|f\|_{\elle^p(X,\nu)}$. By applying this procedure the derivatives $D_H^{\mu,j}U$ arise, with $j=1,\ldots,k-1$, and so we need that the function $U$ belongs to $W^{k-1,p}(X,\mu)$.

\begin{remark}
\label{rmk:U_sob_der_k_nu}
If $U\in W^{k,p}(X,\mu)$ for every $p\in(1,+\infty)$, then $U\in W^{k,p}(X,\nu)$ for every $p\in(1,+\infty)$. It is enough to repeat the computations in Remark \ref{rmk:U_sob_space_nu}.
\end{remark}

\begin{thm}
\label{thm:equivalence_space_k}
Let $p\in(1,+\infty)$ and let $k\in\N$, $k\geq 3$. If $U\in W^{k-1,q}(X,\mu)$ for every $q\in[1,+\infty)$, then there exists positive constants $K_{p}$ and $\widetilde K_{p,k}$ such that
$K_{p}\|f\|_{p,D_{H}^k}\leq \|f\|_{k,p}\leq \widetilde K_{p,k}\|f\|_{p,D^k_H}$ for every $f\in \fcon_b^\infty(X)$.
\end{thm}
\begin{proof}
From the definition we get $\|f\|_{p,D_{H}^k}\leq 2^{1-1/p}\|f\|_{k,p}$ for every $f\in\fcon_b^\infty(X)$.





We limit ourselves to prove second inequality of the statement when $k=3$, the other cases following by analogous computations. 

Let $k=3$ and let $p\in(1,+\infty)$. We notice that from Theorem \ref{thm:equivalece_space_2} there exists a positive constant $\widetilde K_p$ such that \begin{align*}
\|f\|_{2,p}\leq \widetilde K_p\left(\|f\|_{\elle^p(X,\nu)}+\|D_H^2f\|_{\elle^p(X,\nu;\mathcal H_2(H))}\right), 
\end{align*}
for every $f\in\fcon_b^\infty(X)$.
This gives
\begin{align*}
\|f\|^p_{3,p}
\leq  & \|f\|^p_{2,p}+\|D_H^3f\|_{\elle^p(X,\nu;\mathcal H_3(H))}^p \\
\leq & 2^{p-1}(\widetilde K_p)^p\left(\|f\|^p_{\elle^p(X,\nu)}+\|D_H^2f\|_{\elle^p(X,\nu;\mathcal H_2(H))}^p\right) +\|D_H^3f\|^p_{\elle^p(X,\nu;\mathcal H_3(H))},
\end{align*}
for every $f\in\fcon_b^\infty(X)$.
Hence, it is enough to find a positive constant $c$, independent of $f$, such that 
\begin{align*}
\|D_H^2f\|_{\elle^p(X,\nu;\mathcal H_2(H))}\leq c\left(\|f\|_{\elle^p(X,\nu)}+\|D_H^3f\|_{\elle^p(X,\nu;\mathcal H_3(H))}\right),
\end{align*} 
for every $f\in\fcon_b^\infty(X)$.
Arguing as in \eqref{stima_vett_totale} and taking into account \eqref{poincare_vett_derivate}, we get
\begin{align}
\|D_{H}^{2}f\|_{\elle^p(X,\nu;\mathcal H_{2}(H))}
\leq k_p\|D_{H}^{3}f\|_{\elle^p(X,\nu;\mathcal H_{3}(H))}
+\left(\sum_{i,j\in\N}\left|\int_X[D^{2}_{H}fe_{i},e_{j}]_{H}d\nu\right|^2\right)^{1/2},
\label{eq_vettoriale_totale}
\end{align}
for every $f\in\fcon_b^\infty(X)$, where $\{e_n:n\in\N\}$ is any orthonormal basis of $H$. Integrating by parts twice in the second addend of the right-hand side of \eqref{eq_vettoriale_totale} it follows that
\begin{align}
& \int_X[D^{2}_{H}fe_{i},e_{j}]_{H}d\nu \notag\\
&= \int_X[D_{H}f,e_{i}]_{H}(\hat{e_{j}}+[D_H^\mu U,e_{j}]_{H})d\nu \notag \\
& = \int_Xf\big((\hat{e_{i}}+[D_H^\mu U,e_{i}]_{H})(\hat{e_{j}}+[D_H^\mu U,e_{j}]_{H})-\delta_{ij}-[D^{\mu,2}_{H}Ue_{i},e_{j}]_{H}\big) d\nu  \notag \\
& =  \int_Xf(\hat{e_{i}}\hat{e_{j}}-\delta_{ij})d\nu
-\int_Xf[D^{\mu,2}_{H}Ue_{i},e_{j}]_{H}d\nu+\int_Xf[D_H^\mu U,e_{i}]_{H}[D_H^\mu U,e_{j}]_{H}d\nu \notag  \\
& \ \ \ +\int_Xf(\hat{e_{i}}[D_H^\mu U,e_{j}]_{H}+\hat{e_{j}}[D_H^\mu U,e_{i}]_{H})d\nu \notag \\
& =:J_1^{i,j}+J_2^{i,j}+J_3^{i,j}+J_4^{i,j},
\label{somma_da_stimare_equivalenza}
\end{align}
for every $i,j\in\N$, where $\delta_{ij}$ is the Kronecker symbol and $[D_H^\mu\hat {e_i},e_j]_H=\delta_{ij}$ follows from \cite[Lemma 2.10.5]{Bog98}. Let us separately estimate the four terms in the last line of \eqref{somma_da_stimare_equivalenza}. From Remarks \ref{rmk:chaos_ind_base} and  \ref{rmk:dec_wiener_chaos} it follows that $J^{i,j}_1$ is related to the second Wiener chaos $E_2$: indeed,
\begin{align*}
J^{i,j}_1=\int_X\Phi_{\alpha(i,j)}fe^{-U}d\mu, \quad i,j\in\N, \ i\neq j,
\end{align*}
where $\alpha(i,j)$ is the multiindex which satisfies $|\alpha|=2$ and $\alpha_i=\alpha_j=1$, and
\begin{align*}
J^{i,i}_1=\sqrt2\int_X\Phi_{\alpha(i)}fe^{-U}d\mu, \quad i\in\N, 
\end{align*}
where $\alpha(i)$ is the multiindex which satisfies $|\alpha|=2$ and $\alpha_i=2$. Then,
\begin{align}
\label{prima_sommatoria_k=3_p>2}
\sum_{i,j\in\N}\left|J_1^{i,j}\right|^2\!\!
= & \!\!\sum_{i,j\in\N}\left|\int_Xfe^{-U}(\hat{e_{i}}\hat{e_{j}}-\delta_{ij})d\mu\right|^2
=2 \|I_2(fe^{-U})\|_{\elle^2(X,\mu)}^2
\leq \!a_p\|fe^{-U}\|_{\elle^p(X,\mu)}^2
\leq  \! C_p\|f\|_{\elle^p(X,\nu)}^2,
\end{align}
where the last inequality can be obtained arguing as in \eqref{stima_2_D_D^2} and $C_p$ is a positive constant which only depends on $p$ and $U$. 

As far as $J_2^{i,j}$ is considered, we get
\begin{align}
\sum_{i,j\in\N}\left|J_2^{i,j}\right|^2
= & \sum_{i,j\in\N}\left|\int_Xf[D^{\mu,2}_{H}Ue_{i},e_{j}]_{H}d\nu\right|^2
= \left\|\int_XD^{\mu,2}_HU f d\nu\right\|_{\mathcal H_2(H)}^2  \notag \\
\leq & \left( \int_X\|D^{\mu,2}_HU\|_{\mathcal H_2(H)}| f | d\nu\right)^2 
\leq \left \|\|D^{\mu,2}_{H}U\|_{\mathcal H_2(H)}\right\|^2_{\elle^{p'}(X,\nu)}\|f\|_{\elle^p(X,\nu)}^2,
\label{seconda_sommatoria_k=3_p>2}
\end{align}
where in the last inequality we have applied the H\"older's inequality with exponents $p$ and $p'$. From Remark \ref{rmk:U_sob_der_k_nu} we infer that $\left\|\|D^{\mu,2}_{H}U\|_{\mathcal H_2(H)}\right\|^2_{\elle^{p'}(X,\nu)}<+\infty$.

$J_3^{i,j}$ can be estimated as follows:
\begin{align}
\sum_{i,j\in\N}\left|J_3^{i,j}\right|^2
= & \sum_{i,j\in\N}\left|\int_Xf[D_H^\mu U,e_{i}]_{H}[D_H^\mu U,e_{j}]_{H}d\nu\right|^2 \notag 
\leq \left( \int_X|f||D^\mu_HU|_H^2d\nu\right)^2 \\
\leq & \|f\|_{\elle^p(X,\nu)}^2\||D^\mu_HU|_H^{2}\|_{\elle^{p'}(X,\nu)}^2,
\label{stima_3_eq_vettoriale}
\end{align}
and the last inequality follows from the H\"older's inequality with exponents $p$ and $p'$. 
Again, Remark \ref{rmk:U_sob_der_k_nu} gives
$\||D^\mu_{H}U|_H^2\|^2_{\elle^{p'}(X,\nu)}<+\infty$.


Finally, we take $J^{i,j}_4$ into account. We get
\begin{align}
\sum_{i,j\in\N}\left|J_4^{i,j}\right|^2
= & \sum_{i,j\in\N}\left|\int_Xf(\hat{e_{i}}[D_H^\mu U,e_{j}]_{H}+\hat{e_{j}}[D_H^\mu U,e_{i}]_{H})d\nu\right|^2 \notag \\
\leq & 2\sum_{i,j\in\N}\left|\int_Xf\hat{e_{i}}[D_H^\mu U,e_{j}]_{H}d\nu\right|^2+2\sum_{i,j\in\N}\left|\int_X\hat{e_{j}}[D_H^\mu U,e_{i}]_{H}d\nu\right|^2 \notag \\
= & 4\sum_{i,j\in\N}\left|\int_Xf\hat{e_{i}}[D_H^\mu U,e_{j}]_{H}e^{-U}d\mu\right|^2.
\label{eq_vettorial_4}
\end{align}
We recall that for every $j\in\N$ the element
\begin{align*}
\int_Xf\hat{e_{i}}[D_H^\mu U,e_{j}]_{H}e^{-U}d\mu,
\end{align*}
is the projection of $f[D_H^\mu U,e_{j}]_{H}e^{-U}$ on the subspace of $E_1$ generated by $\hat {e_i}$, for every $i\in\N$. Hence, from Lemma \ref{lemma:L2-Lp-norm} we have
\begin{align}
\sum_{i\in\N}\left|\int_Xf\hat{e_{i}}[D_H^\mu U,e_{j}]_{H}e^{-U}d\mu\right|^2
= & \|I_1(f[D_H^\mu U,e_{j}]_{H}e^{-U})\|_{\elle^2(X,\nu)}^2
\leq a_p\|f[D^\mu_HU,e_j]_He^{-U}\|_{\elle^q(X,\mu)}^2,
\label{eq_vettorial_4_1}
\end{align}
for every $j\in\N$, where $q=\min\left\{\frac{1+p}{2},\frac32\right\}$ and $a_p$ is a positive constant which only depends on $p$. By applying the Minkowski integral inequality, as in \eqref{stima_vett_minkowski}, with indeces $i,j\in\N$, we infer that
\begin{align}
\label{eq_vettorial_4_2}
\sum_{j\in\N}\|f[D^\mu_HU,e_j]_He^{-U}\|_{\elle^q(X,\mu)}^2
\leq \|f |D^\mu_HU|_H e^{-U}\|_{\elle^q(X,\mu)}^2.
\end{align}
We apply the H\"older's inequality with $r=\frac{p}{q}$ and $r'=\frac{p}{p-q}$ and we get
\begin{align}
\label{eq_vettorial_4_3}
\|f |D^\mu_HU|_H e^{-U}\|_{\elle^q(X,\mu)}^2
\leq \|f\|_{\elle^p(X,\nu)}^2\left(\int_X|D^\mu_HU|_H^{qr'}e^{-r'(q-\frac qp)U}d\mu\right)^{2/{qr'}}
\leq C_p \|f\|_{\elle^p(X,\nu)}^2,
\end{align}
where $C_p$ is a positive constant which only depends on $p$ and $U$. By collecting \eqref{eq_vettoriale_totale}-\eqref{eq_vettorial_4_3} and by noticing that
\begin{align*}
\left(\sum_{i,j\in\N}\left|\int_X[D^{2}_{H}fe_{i},e_{i_{2}}]_{H}d\nu\right|^2\right)^{1/2}
\leq 2\left(\sum_{i,j\in\N}\left(\left|J_1^{i,j}\right|^2+\left|J_2^{i,j}\right|^2+\left|J_3^{i,j}\right|^2+\left|J_4^{i,j}\right|^2\right)\right)^{1/2},
\end{align*}
we conclude that there exists a positive constant $\widetilde C_p$, which only depends on $p$ and $U$, such that
\begin{align*}
\left(\sum_{i,j\in\N}\left|\int_X[D^{2}_{H}fe_{i},e_{i_{2}}]_{H}d\nu\right|^2\right)^{1/2}
\leq \widetilde C_p\|f\|_{\elle^p(X,\nu)}.
\end{align*}
This gives the thesis for $k=3$.

The cases $k\geq4$ can be obtained arguing as for $k=3$, simply iterating the integration by parts as in \eqref{somma_da_stimare_equivalenza} $k-1$-times and estimating the terms which arise as for $k=3$. Computations are long but straightforward, and we left them to the reader.
\end{proof}

The following result follows from Theorem \ref{thm:equivalence_space_k}, arguing as in the proof of Corollary \ref{coro:eq_norme_2}. 
\begin{corollary}
\label{coro:eq_norme_k}
Let $k\in\N$, $k\geq3$, and assume that $U\in W^{k-1,p}(X,\mu)$ for every $p\in(1,+\infty)$. Then, the operator $D_H^3:\fcon_b^\infty(X)\subseteq\elle^p(X,\nu)\to \elle^p(X,\nu;\mathcal H_k(H))$ is closable in $\elle^p(X,\nu)$ and the domain of its closure, endowed with the graph norm, coincides with $W^{k,p}(X,\nu)$.
\end{corollary}

We conclude this subsection by providing a sufficient condition for the belonging of a function $f$ in $W^{k,p}(X,\nu)$.
\begin{lemma}
Let $p\in(1,+\infty)$, $k\geq 2$ and assume that $U\in W^{k,q}(X,\mu)$ for every $q\in(1,+\infty)$. If $(f_n)\subseteq W^{k,p}(X,\nu)$ converges to $f$ $\nu$-a.e. in $X$ and $\sup_{n\in\N}\|f_n\|_{p,D_H^k}<+\infty$, then $f\in W^{k,p}(X,\nu)$.    
\end{lemma}
\begin{proof}
From Corollary \ref{coro:eq_norme_k} it follows that $\sup_{n\in\N}\|f_n\|_{k,p}<+\infty$. The thesis follows by adapting the argument of the proof of \cite[Lemma 5.4.4]{Bog98}.   
\end{proof}

\subsection{Exponential decay}
\begin{proposition}
Let $V$ be a separable Hilbert space. Then:
\begin{enumerate}[(i)]
\item For every $p\in[2,+\infty)$, every $t>1$ and every $F\in\elle^p(X,\nu;V)$, we get
\begin{align*}
\|T^V(t)F-\nu(F)\|_{\elle^p(X,\nu;V)}\leq d_pe^{-t}\|F\|_{\elle^p(X,\nu;V)},
\end{align*}
where $d_p$ is the positive constant given by
\begin{align}
\label{vett_ex_decay_costante_p>2}
d_p:=e\sqrt{\frac{p-1}{2}}.
\end{align}
\item For every $p\in(1,2)$, every $t>2$ and every $F\in\elle^p(X,\nu;V)$ we have
\begin{align*}
\|T^V(t)F-\nu(F)\|_{\elle^p(X,\nu;V)}\leq \frac{(c_p)^{1/p}e^2}{\sqrt{2}}e^{-t}\|F\|_{\elle^p(X,\nu;V)},
\end{align*}
where $c_p$ is the positive constant in \eqref{stima_grad_funz_prima_parte_finale}.
\end{enumerate}
\end{proposition}
\begin{proof}
Let $V$ be a separable Hilbert space. By density, we can limit ourselves to prove the statement for $F\in\fcon_b^\infty(X;V)$.

\vspace{2mm}
{{\bf $(i)$}}. Let $p\in[2,+\infty)$ and let $t>1$. We fix $F\in\fcon_b^\infty(X;V)$, and we set $G:=T^V(t)F-\nu(F)\in\fcon_b^\infty(X)$. From \eqref{der_G_vett_poincare}, \eqref{poincare_vett_1} and \eqref{poincare_vett_2} it follows that
\begin{align}
\label{vett_exp_decay_1}
\int_X|G|^p_Vd\nu
=\int_t^\infty\left(\int_X[\overline{D_H}T^V(s)F,\overline{D_H}G^*]_{H\otimes V}d\nu\right)ds,
\end{align}
where $G^*:=|G|^{p-2}G$. Since $T^V(t)F\in W^{1,p}(X,\nu;V)$, from \eqref{der_duale} we infer that $G^*\in W^{1,p'}(X,\nu)$ and
\begin{align}
\label{vett_exp_decay_2}
\overline{D_H}G^*= (p-2)|G|^{p-4}_V[\overline {D_{H}}T^V(t)F,G]_V\otimes G+|G|^{p-2}_V\overline {D_{H}}T^V(t)F.
\end{align}
Arguing as in \eqref{poincare_vett_3}, from \eqref{vett_exp_decay_1} and \eqref{vett_exp_decay_2} we deduce that
\begin{align*}
\int_X|G|^p_Vd\nu
\leq & (p-1)\int_t^\infty\left( \int_X|\overline {D_H}T^V(s)F|_{H\otimes V}|\overline {D_H}T^V(t)F|_{H\otimes V}|G|^{p-2}_Vd\nu \right)ds.
\end{align*}
Let $p>2$. We apply the generalized H\"older's inequality with exponents $p,p$ and $\frac{p}{p-2}$, which gives
\begin{align*}
\int_X|G|^p_Vd\nu
\leq & (p-1)\left(\int_t^\infty\|\overline {D_H}T^V(s)F\|_{\elle^p(X,\nu;H\otimes V)}ds\right)\|\overline {D_H}T^V(t)F\|_{\elle^p(X,\nu;H\otimes V)}\|G\|_{\elle^p(X,\nu;V)}^{p-2}.
\end{align*}
Recalling \eqref{stima_grad_funz_mista_2+infty} and $t>1$, it follows that
\begin{align}
\int_X|G|^p_Vd\nu
\leq & \frac{(p-1)e^{-t+1}}{2}\int_t^\infty e^{-s+1}ds\|F\|_{\elle^p(X,\nu;V)}^2\|G\|_{\elle^p(X,\nu;V)}^{p-2} \notag \\
= & \frac{(p-1)e^{2}}{2}e^{-2t}\|F\|_{\elle^p(X,\nu;V)}^2\|G\|_{\elle^p(X,\nu;V)}^{p-2}.
\label{vett_exp_decay_3}
\end{align}
Dividing both the sides of \eqref{vett_exp_decay_3} by $\|G\|_{\elle^p(X,\nu;V)}^{p-2}$ we get
\begin{align*}
\|T(t)F-\nu(F)\|_{\elle^p(X,\nu;V)}\leq d_pe^{-t}\|F\|_{\elle^p(X,\nu;V)},
\end{align*}
with $d_p$ is the constant in \eqref{vett_ex_decay_costante_p>2}. The case $p=2$ follows noticing that $G=G^*$ and $\overline{D_H}G^*=\overline{D_H}F$, applying the H\"older's inequality with $p=2$ in \eqref{vett_exp_decay_1} and concluding as for $p>2$ by means of \eqref{stima_grad_funz_mista}.

\vspace{2mm}
{{\bf $(ii)$}}. Let $p\in(1,2)$, let $t>2$ and let $F\in\fcon_b^\infty(X;V)$. We have
\begin{align}
\label{vett_ex_decay_p>2_1}
\|T^V(t)F-\nu(F)\|_{\elle^p(X,\nu;V)}
= & \sup_{
\begin{array}{l}
G\in\fcon_b^\infty(X;V),\\ 
\|G\|_{\elle^{p'}(X,\nu;V)}\leq 1
\end{array}}\int_X[F-\nu(F),G]_Vd\nu.
\end{align}
Let $G\in \fcon_b^\infty(X;V)$ with $\|G\|_{\elle^{p'}(X,\nu;V)}\leq 1$. From the semigroup property of $(T^V(t))_{t\geq0}$, Lemma \ref{lem:prop_vettoriali}$(i)$ and \eqref{convergenza_as_vett} it follows that
\begin{align}
\int_X[T^V(t)F-\nu(F),G]_Vd\nu
= & \lim_{r\rightarrow +\infty}\int_X[T^V(t)F-T^V(r)F,G]_Vd\nu \notag \\
= & \lim_{r\rightarrow +\infty}\int_X[T^V(t/2)F,T^V(t/2)G-T^V(r-t/2)G]_Vd\nu,
\label{vett_ex_decay_p>2_2}
\end{align}
Arguing as in \eqref{poincare_vett_stima_duale_3},  
for every $r>t$  we have
\begin{align}
& \int_X]T^V(t/2)F,T^V(t/2)G-T^V(r-t/2)G]_Vd\nu \notag \\
= & \int_{t/2}^{r-t/2}\left(\int_X[\overline {D_{H}}T^V(t/2)F,\overline{D_{H}}T^V(s)G]_{H\otimes V}d\nu\right)ds \notag \\
\leq &  \left(\int_X|\overline {D_{H}}T^V(t/2)F|_{H\otimes V}^pd\nu\right)^{1/p}\!\!
\int_{t/2}^{r-t/2}\left(\int_X|\overline{D_{H}}T^V(s)G|^{p'}_{H\otimes V}d\nu\right)^{1/p'}\!\!ds.
\label{vett_ex_decay_p>2_3}
\end{align}

{By applying \eqref{stima_grad_funz_mista_2+infty} to $\int_X|\overline{D_{H}}T^V(s)G|^{p'}_{H\otimes V}d\nu$ and recalling that $\|G\|_{\elle^{p'}(X,\nu;V)}\leq 1$ and that $s>t/2>1$, we infer that
\begin{align*}
\left(\int_X|\overline{D_{H}}T^V(s)G|^{p'}_{H\otimes V}d\nu\right)^{1/p'}
\leq& \frac{e^{-s+1}}{\sqrt2}, \quad s\in(t/2,r-t/2),
\end{align*}
which gives
\begin{align}
\label{vett_ex_decay_p>2_4}
\int_{t/2}^{r-t/2}\left(\int_X|\overline{D_{H}}T^V(s)G|^{p'}_{H\otimes V}d\nu\right)^{1/p'}\!\!ds
\leq &\frac{e}{\sqrt2} \int_{t/2}^{r-t/2}e^{-s}ds
= \frac{e}{\sqrt2}(e^{-t/2}-e^{-r+t/2}),
\end{align}
for every $r>t$. Further, since $t>2$ estimate \eqref{stima_grad_funz_mista} gives
\begin{align}
\label{vett_ex_decay_p>2_5}
\left(\int_X|\overline {D_{H}}T^V(t/2)F|_{H\otimes V}^pd\nu\right)^{1/p}
\leq (c_p)^{1/p}{e^{-t/2+1}}\|F\|_{\elle^p(X,\nu;V)}.
\end{align}}

{Collecting \eqref{vett_ex_decay_p>2_1}-\eqref{vett_ex_decay_p>2_5} we get
\begin{align}
\|T^V(t)F-\nu(F)\|_{\elle^p(X,\nu;V)}
\leq \frac{(c_p)^{1/p}e^2}{\sqrt2}e^{-t}\|F\|_{\elle^p(X,\nu; V)},
\label{vett_exp_decaystima_fina_p<2}
\end{align}
which gives the thesis for $p\in(1,2)$.}
\end{proof}

\appendix
\section{Smooth approximations}
\label{app_A}
This section is devoted to the approximation procedure which have been exploited to obtain gradient estimates for the semigroup $(T(t))_{t\geq0}$. The results are quite technical and follows the lines of the proofs of similar statements in \cite{AngFerPal18}, but we deal with Banach spaces instead of Hilbert spaces and we weaken the assumptions on $U$. Hence, we provide the details for reader's convenience. 

We recall the definition of $H$-Lipschitz function: let $Y$ be a separable Banach space with norm $\|\cdot\|_Y$, we say that $F:X\rightarrow Y$ is a $H$-Lipschitz continuous function if there exists a positive constant $M$ such that
\begin{align}
\label{def_H_lips}
\|F(x+h)-F(x)\|_Y\leq M|h|_{H}, \quad \forall h\in H, \ \mu\textup{-a.e. }x\in X.
\end{align}
We denote by $[F]_{H-{\rm Lip}}$ the smaller constant $M$ which realizes \eqref{def_H_lips}.

\begin{remark}
\label{rmk:def_proj}
Hereafter, we fix an orthonormal basis $\Phi:\{e_m=i^*(e_m^*):m\in\N\}$ of $H$ (the existence of this basis follows from the density of $j(X^*)$ in $\elle^2(X,\mu)$ and the isometry between $\elle^2(X,\mu)$ and $H$), and for every $k\in\N\cup\{\infty\}$ we set $\fcon_{b,\Phi}^k(X):=\{f\in \fcon_b^k(X):f(x)=\varphi(\langle x, e_1^*	\rangle, \ldots,\langle x,e_m^*\rangle), \ m\in\N, \ \varphi\in C_b^k(\R^m)\}$. Clearly, $\fcon_b^\infty(X)$ is dense in $\elle^P(X,\nu)$ for every $p\in[1,+\infty)$. Further, for every $p\in[1,+\infty)$ the closure of the operator $D_{H}:\fcon_{b,\Phi}^\infty(X)\rightarrow \elle^p(X,\nu;H)$ and the domain of its closure coincide with $D_{H}$ and with $W^{1,p}(X,\nu)$, respectively. Indeed, let $p\in[1,+\infty)$ and let $f\in\fcon_b^\infty(X)$ be such that $f(x)=\varphi(\langle x,x^*_1\rangle, \ldots, \langle x,x^*_n\rangle)$ for some $n\in\N$, $\varphi\in C_b^\infty(\R^n)$ and $x_1^*,\ldots,x_n^*\in X^*$. For every $m\in\N$ let us denote by $P_m:X\rightarrow H$ the projection
\begin{align*}
P_mx:=\sum_{i=1}^m\langle x, e_i^*\rangle e_i, \quad x\in X,
\end{align*}
on ${\rm span}\{e_1,\ldots,e_m\}$, and let us set $f_m(x):=f(P_mx)$ for every $x\in X$. Clearly, $f_m\in\fcon_{b,\Phi}^\infty(X)$ for every $m\in\N$. Further, from \cite[Corollary 3.5.8]{Bog98} and from the dominated convergence theorem we infer that $f_m\rightarrow f$ in $\elle^p(X,\nu)$ as $n\rightarrow+\infty$. Let us compute $D_{H}f_m$. We have
\begin{align*}
D_{H}f_m(x)=P_mD_{H}f(P_m x), \quad x\in X, \ m\in\N.
\end{align*}
Again, \cite[Corollary 3.5.8]{Bog98} and the dominated convergence theorem give $D_{H}f_m\rightarrow D_{H}f$ in $\elle^p(X,\nu;H)$ as $n\rightarrow+\infty$. As a byproduct, $\fcon_{b,\Phi}^\infty(X)$ is dense in $W^{1,p}(X,\nu)$ for every $p\in[1,+\infty)$. Analogously, it is possible to prove that $\fcon_{b,\Phi}^\infty(X)$ is dense in $W^{k,p}(X,\nu)$ for every $k\in\N$ and every $p\in[1,+\infty)$. Finally, we set $\Pi_n:X\rightarrow \R^n$ the projection defined by
\begin{align*}
\Pi_nx:=(\langle x,e_1^*\rangle,\ldots,\langle x, e^*_n\rangle)\in\R^n, \quad x\in X,
\end{align*}
and for every $n\in\N$  we define $\Sigma_n:\R^n\rightarrow H$ as
\begin{align*}
\Sigma_n\xi:=\sum_{i=1}^n\xi_ie_i, \quad \xi=(\xi_1,\ldots,\xi_n)\in\R^n.
\end{align*}
\end{remark}

For every $f\in \elle^1(X,\mu)$ we introduce the conditional expectation $\mathbb E_nf$ as
\begin{align}
\label{conditional_ex_n}
\mathbb E_nf(x):=\int_Xf(P_nx+(I-P_n)y)\mu(dy), \quad \mu\textup{-a.e. }x\in X,
\end{align}
where $P_n$ has been defined in Remark \ref{rmk:def_proj}. We recall the following results (see \cite[Corollary 3.5.2 \& Proposition 5.4.5]{Bog98}).
\begin{proposition}
\label{prop:cond_exp}
Let $p\in[1,+\infty)$ and let $f\in \elle^p(X,\mu)$. Then, $\mathbb E_nf\rightarrow f$ in $\elle^p(X,\mu)$ and for every $n\in\N$ we have
\begin{align*}
\|\mathbb E_nf\|_{\elle^p(X,\mu)}\leq \|f\|_{\elle^p(X,\mu)}.
\end{align*} 
Further, if $f\in W^{1,p}(X,\mu)$ then $\mathbb E_nf\rightarrow f$ in $W^{1,p}(X,\mu)$ and for every $n\in\N$ we have
\begin{align*}
[D_{H}(\mathbb E_nf),e_j]_{H}
= \begin{cases}
\mathbb E_n([D_{H}f,e_j]_{H}), & 1\leq j\leq n, \\
0, & j\geq n+1,
\end{cases}
\end{align*}
and $\|D_H(\mathbb E_nf)\|_{\elle^p(X,\mu)}\leq \|D_Hf\|_{\elle^p(X,\mu)}$. Finally, the same results, with the obvious modifications, are true for $f\in W^{2,p}(X,\mu)$.
\end{proposition}

The final tool which we need are the {\it Moreau-Yosida approximants} of $U$ along $H$. Below we state the main results we use in the following, and we refer to \cite[Section 12.4]{BC11} for the classical theory, and to \cite{ACF20,AngFerPal18,CF16} for the case here considered.

\begin{proposition}
\label{prop:mor_yos_appr}
Let $f:X\rightarrow \R\cup\{+\infty\}$ be a proper convex and $|\cdot|_X$-lower semicontinuous function and denote by ${{\rm dom}}(f)=\{x\in X:f(x)<+\infty\}$. For every $\varepsilon>0$ and every $x\in X$ we set
\begin{align*}
f_\varepsilon(x):=\inf\left\{f(x+h)+\frac1{2\varepsilon}|h|_{H}^2\Big|h\in H\right\}.
\end{align*}
Then:
\begin{enumerate}[(i)]
\item $f_\varepsilon(x)\leq f(x)$ for every $\varepsilon>0$ and every $x\in X$. Moreover, $f_\varepsilon(x)$ monotonically converges to $f(x)$ as $\varepsilon\rightarrow 0^+$ for every $x\in X$.
\item $f_\varepsilon$ is $H$-differentiable and $D_{H}f_\varepsilon$ is a $H$-valued $H$-Lipschitz continuous function in $X$ with constant less or equal to $1/\varepsilon$. Hence, from \cite[Theorem 5.11.2]{Bog98} it follows that $D_Hf_\varepsilon\in W^{1,\infty}(X,\mu;H)$ and $|D_H^2 f_\varepsilon(x)|_{\mathcal H_2(H)}\leq \varepsilon^{-1}$ for $\mu$-a.e. $x\in X$.
\item $f_\varepsilon\in W^{2,p}(X,\mu)$ whenever $f\in\elle^p(X,\mu)$ with $p\in[1,+\infty)$.
\item If $x\in {\rm dom}(f)$ and $f\in W^{1,p}(X,\mu)$ for some $p\in[1,+\infty)$ then $D_{H}f_\varepsilon(x)$ converges to $D_{H}f(x)$ as $\varepsilon\rightarrow0^+$.
\end{enumerate}
\end{proposition}

Let us introduce smooth approximations of $U$. For every $\varepsilon>0$, let $U_\varepsilon$ be the Moreau-Yosida approximants of $U$. Further, we set
\begin{align*}
\psi_{\varepsilon,n}(\xi):=\mathbb E_n(U_\varepsilon)(\Sigma_n\xi), \quad 
\psi_{\varepsilon,n,\eta}(\xi):=(\psi_{\varepsilon,n}*\theta_\eta)(\xi), \quad \xi\in \R^n, 
\end{align*}
where the second term is the convolution of $\psi_{\varepsilon,n}$ with the family of mollifiers $\theta_\eta$, $\eta\in\R^+$. Here, $\theta_\eta(\xi)=\theta(\xi\eta)$ and $\theta\in C_c^\infty(\R^n)$ satisfies $0\leq \theta\leq 1$,  has support contained in the unit ball and $\int_{\R^n}\theta(\xi)d\xi=1$. For every $\varepsilon>0$ we set
\begin{align*}
\nu_\varepsilon:=e^{-U_\varepsilon}\nu.
\end{align*}
Arguing as in \cite[Proposition 5.12]{CF16} it follows that there exists $x^*\in X^*$ and $r\in \R$ such that
\begin{align*}
U_\varepsilon(x)\geq \langle x,x^*\rangle+r, \quad x\in X,
\end{align*}
for every $\varepsilon\in(0,1]$. This fact and Fernique Theorem (see \cite[Theorem 2.8.5]{Bog98}) imply that for every $\varepsilon\in(0,1]$ we have $e^{-U_\varepsilon}\in\elle^p(X,\mu)$ for every $p\in[1,+\infty)$, and so the measure $\nu_\varepsilon$ is well defined. We also deduce the following useful lemma.
\begin{lemma}
\label{lem:conv_U_epsilon}
Let $U_\varepsilon$ be as above. Then, for every $p\in[1,+\infty)$ and every decreasing and vanishing sequence $(\varepsilon_n)_{n\in\N}$ we have
\begin{align*}
\lim_{n\rightarrow+\infty}\|e^{-U_{\varepsilon_n}}-e^{-U}\|_{\elle^p(X,\mu)}=0.
\end{align*}
\end{lemma}
\begin{proof}
Thank to Proposition \ref{prop:mor_yos_appr}$(i)$, the thesis follows from the monotone convergence theorem.
\end{proof}
\begin{proposition}
\label{prop:conv_Ueps}
For every $\varepsilon\in(0,1]$, every $n\in\N$ and every $\eta\in\R^+$, the function $\psi_{\varepsilon,n,\eta}\in C_b^\infty(\R^n)$. Further, if we set $U_{\varepsilon,n,\eta}(x):=\psi_{\varepsilon,n,\eta}(\Sigma_n\Pi_nx)$, $x\in X$, then for every $\varepsilon\in(0,1]$ there exists a decreasing vanishing sequence $(\eta_n)\subseteq \R^+$ such that
\begin{align}
\label{convergenza_Uepsneta1}
& \elle^2(X,\nu_\varepsilon;H)-\lim_{n\rightarrow+\infty}D_H U_{\varepsilon,n,\eta_n}=D_HU_\varepsilon, \\
\label{convergenza_Uepsneta2}
& \elle^2(X,\nu_\varepsilon;\mathcal H_2(H))-\lim_{n\rightarrow+\infty}D^2_{H}U_{\varepsilon,n,\eta_n}=D^2_HU_\varepsilon.
\end{align}
We set $\widetilde U_{\varepsilon,n}:=U_{\varepsilon,n,\eta_n}$ for every $\varepsilon\in(0,1]$ and every $n\in\N$.
\end{proposition}
\begin{proof}
The first part of the statement follows arguing as in \cite[Lemma 2.3]{AngFerPal18}. Further, the belonging of $U_\varepsilon$ to $W^{2,p}(X,\nu_\varepsilon)$ follows from Proposition \ref{prop:mor_yos_appr}.  To show \eqref{convergenza_Uepsneta1}, at first we prove that for every vanishing sequence $(\eta_n)$ we have
\begin{align*}
\elle^2(X,\nu_\varepsilon;H)-\lim_{n\rightarrow+\infty}D_H  U_{\varepsilon,n,\eta_n}=D_HU_\varepsilon.
\end{align*}
Let $U_{\varepsilon,n}:=\mathbb E_n(U_\varepsilon)$ for every $ \varepsilon\in(0,1]$ and every $n\in\N$. Then, we have
\begin{align*}
& \|D_H U_{\varepsilon,n,\eta_n}-D_HU_\varepsilon\|_{\elle^2(X,\nu_\varepsilon;H)}^2 \\
\leq & 2\|D_H U_{\varepsilon,n}-D_HU_\varepsilon\|_{\elle^2(X,\nu_\varepsilon;H)}^2
+2\|D_H U_{\varepsilon,n,\eta_n}-D_HU_{\varepsilon,n}\|_{\elle^2(X,\nu_\varepsilon;H)}^2
=:I_1^n+I_2^n.
\end{align*}
As far as $I_1^n$ is considered, for every $p\in(1,+\infty)$ from Proposition \ref{prop:cond_exp} we get
\begin{align*}
I_1^n\leq \left(\int_Xe^{-p'U_\varepsilon}d\mu\right)^{1/p'}\|D_H U_{\varepsilon,n}-D_HU_\varepsilon\|_{\elle^{2p}(X,\mu;H)}^{2/p}\rightarrow 0, \quad n\rightarrow+\infty,
\end{align*}
where $p'$ is the conjugate exponent of $p$. Let us estimate $I_2^n$. From the definition of $\psi_{\varepsilon,n}$, $U_{\varepsilon,n}$, $\psi_{\varepsilon,n,\eta_n}$ and $U_{\varepsilon,n,\eta_n}$, and from Proposition \ref{prop:mor_yos_appr}$(ii)$ it follows that
\begin{align*}
I_2^n
= & \int_X\left|\int_X\left(D_HU_\varepsilon(P_nx+(I-P_n)y)-\int_{\R^n}D_HU_\varepsilon\left(P_nx+(I-P_n)y-\eta_n\big(\Sigma_n\xi\big)\right)\theta(\xi)d\xi\right)d\mu\right|_H^2d\nu_\varepsilon \\
\leq & [D_HU_\varepsilon]_{\rm Lip}^2\nu_\varepsilon(X)\eta_n\int_{\R^n}|\xi|^2\theta(\xi)d\xi \rightarrow0, \quad n\rightarrow+\infty.
\end{align*}
Let us consider the convergence of the second order derivatives. We claim that for every $\varepsilon\in(0,1]$ and every $n\in\N$, $D^2_HU_{\varepsilon,n,\eta}\rightarrow D^2_HU_{\varepsilon,n}$ in $\elle^2(X,\nu_\varepsilon;\mathcal H_2(H))$ as $\eta\rightarrow0^+$. If the claim is true, a diagonal argument as in the proof of \cite[Lemma 2.4]{AngFerPal18} allows us to conclude. It remains to prove the claim. We have
\begin{align}
& \|D^2_{H}U_{\varepsilon,n,\eta}-D^2_HU_{\varepsilon,n}\|_{\elle^2(X,\nu_\varepsilon;\mathcal H_2(H))}^2\notag \\
& \leq \int_X\left(\int_{\R^n}\left|D_H^2U_{\varepsilon,n}(x)-D_H^2U_{\varepsilon,n}(x-\eta\big(\Sigma_n \xi\big))\right|_{\mathcal H_2(H)}^2\theta(\xi)d\xi\!\right)e^{-U_\varepsilon(x)}\mu(dx) \notag\\
& =\int_{\R^n}\left(\int_{X}\left|D_H^2U_{\varepsilon,n}(x)-D_H^2U_{\varepsilon,n}(x-\eta\big(\Sigma_n \xi\big))\right|_{\mathcal H_2(H)}^2e^{-U_\varepsilon(x)}\mu(dx)\right)\theta(\xi)d\xi \notag \\
& \leq \|e^{-U_1}\|_{\elle^{q'}(X,\mu)}\int_{\R^n}\|D^2U_{\varepsilon,n}-D^2U_{\varepsilon,n}(\cdot-\eta\Sigma_n\xi)\|_{\elle^{2q}(X,\mu;{\mathcal H_2(H)})}^{2}\theta(\xi)d\xi,
\label{conv_der_seconde_n_1}
\end{align}
for every $q\in(1,+\infty)$.  From Propositions \ref{prop:cond_exp} and \ref{prop:mor_yos_appr}$(ii)$ it follows that $D^2_HU_{\varepsilon,n}\in \elle^q(X,\mu;\mathcal H_2(H))$ for every $q\in(1,+\infty)$, hence from \cite[Theorem 2.4.8]{Bog98} the function $H\ni h\mapsto[D^2_HU_{\varepsilon,n}(\cdot+h)e_i,e_j]_H$ is continuous from $H$ into $\elle^p(X,\mu)$, for every $i,j\in\N$. This implies that for every $h\in H$ the function
\begin{align*}
X\ni x\mapsto D^2_HU_{\varepsilon,n}(x+h),
\end{align*}
is well defined for $\mu$-a.e. $x\in X$, and it is $\mu$-measurable in $X$. In particular, since for every $\xi\in \R^n$ we have $\Sigma_n\xi\in H$, the function $x\mapsto D^2U_{\varepsilon,n}(x-\eta\Sigma_n\xi)$ is well defined for $\mu$-a.e. $x\in X$ and $\mu$-measurable. 
From Proposition \ref{prop:cond_exp} it follows that
\begin{align}
\|D^2U_{\varepsilon,n}-D^2U_{\varepsilon,n}(\cdot-\eta\Sigma_n\xi)\|_{\elle^{2q}(X,\mu;{\mathcal H_2(H)})}\leq \|D^2U_{\varepsilon}-D^2U_{\varepsilon}(\cdot-\eta\Sigma_n\xi)\|_{\elle^{2q}(X,\mu;{\mathcal H_2(H)})}.
\label{conv_der_seconde_n_2}
\end{align}
From \eqref{conv_der_seconde_n_1} and \eqref{conv_der_seconde_n_2} it follows that
\begin{align}
\notag
& \|D^2_{H}U_{\varepsilon,n,\eta}-D^2_HU_{\varepsilon,n}\|_{\elle^2(X,\nu_\varepsilon;\mathcal H_2(H))}^2 \\
 & \leq \|e^{-U_1}\|_{\elle^{q'}(X,\mu)}\int_{\R^n}\|D^2U_{\varepsilon}-D^2U_{\varepsilon}(\cdot-\eta\Sigma_n\xi)\|_{\elle^{2q}(X,\mu;{\mathcal H_2(H)})}^{2}\theta(\xi)d\xi, \quad \eta>0.
\label{conv_der_seconde_n_3}
\end{align}
Finally, we notice that
\begin{align*}
\|D^2U_{\varepsilon}-D^2U_{\varepsilon}(\cdot-\eta\Sigma_n\xi)\|_{\elle^{2q}(X,\mu;{\mathcal H_2(H)})}
\leq \||D^2U_{\varepsilon}|^2_{\mathcal H_2(H)}-|D^2U_{\varepsilon}(\cdot-\eta\Sigma_n\xi)|^2_{\mathcal H_2(H)}\|_{\elle^{q}(X,\mu)}.
\end{align*}
Again, from \cite[Theorem 2.4.8]{Bog98} the function $h\mapsto |D^2U_{\varepsilon}(\cdot+h)|^2_{\mathcal H_2(H)}$ is continuous from $H$ into $\elle^q(X,\mu)$ for every $q\in(1,\infty)$. Since $\eta\Sigma_n\xi\rightarrow0$ in $H$ for every $\xi\in\R^n$, $D^2U_{\varepsilon,n}$ is essentially bounded on $X$ and $\|\theta\|_\infty\leq 1$, the dominated convergence theorem applied to the right-hand side of \eqref{conv_der_seconde_n_3} gives
\begin{align*}
 \|D^2_{H}U_{\varepsilon,n,\eta}-D^2_HU_{\varepsilon,n}\|_{\elle^2(X,\nu_\varepsilon;\mathcal H_2(H))}^2\rightarrow0, \quad \eta\rightarrow 0^+,
\end{align*}
and the claim is so proved.
\end{proof}

By means of the family of measures $\nu_\varepsilon$, $\varepsilon\in(0,1]$, we introduce a family of operators $L_\varepsilon$. We set
\begin{align}
D(L^\varepsilon_2) & :=\left\{u\in W^{1,2}(X,\nu_\varepsilon):\exists g\in \elle^2(X,\nu_\varepsilon), \ \mathcal \int_X\!\![D_Hu,D_Hv]d\nu_\varepsilon=-\int_X\!\!gvd\nu_\varepsilon, \ \forall v\in\fcon_b^\infty(X)\right\} ,\vspace{2mm}\notag \\
L^\varepsilon_2u& :=g,
\label{def_op_L_epsilon}
\end{align}
where $W^{1,2}(X,\nu_\varepsilon)$ is the domain of the closure of $D_H:\fcon_b^\infty(X)\rightarrow \elle^2(X,\nu_\varepsilon;H)$ in $\elle^2(X,\nu_\varepsilon)$. We denote by $(T_\varepsilon)_{t\geq0}$ the analytic symmetric strongly continuous semigroup of contractions generated by $L_\varepsilon$ on $\elle^2(X,\nu_\varepsilon)$.

On smooth functions $u$ we have an explicit expression of $L_2^\varepsilon u$ and of its $H$-derivative.
\begin{proposition}
\label{propo:explicit_L_2_varepsilon}
We have $\fcon_b^2(X)\subseteq D(L^\varepsilon_2)$ and for every $u\in \fcon_b^2(X)$ it follows that
\begin{align}
\label{esplicit_L2_varepsilon}
(L^\varepsilon_2u)(x)={\rm Tr}[D^2_{H}u(x)]_{H}-\langle x,Du(x)\rangle-[D_H U_\varepsilon(x),D_H u(x)]_{H}, \quad x\in X.
\end{align}
Further, if $u\in\fcon_b^3(X)$ then $L^\varepsilon_2u\in W^{1,2}(X,\nu)$ and for every $h\in H$ we have
\begin{align}
\label{esplicit_DL2_varepsilon}
[(D_{H}L^\varepsilon_2u)(x),h]_{H}=\left(L^\varepsilon_2[D_H u(\cdot),h]_{H}\right)(x)-[D_H u(x),h]_{H}-[D^{2}_{H}U_\varepsilon(x)h,D_H u(x)]_{H},
\end{align}
for $\mu$-a.e. $x\in X$.
\end{proposition}
\begin{proof}
Formula \eqref{esplicit_L2_varepsilon} is well known and it is a direct consequence of integration by parts \eqref{int_by_parts}. Further, By differentiating \eqref{esplicit_L2_varepsilon}, long but straightforward computations give \eqref{esplicit_DL2_varepsilon}.
\end{proof}

\begin{remark}
Let $\varepsilon\in(0,1]$. If we denote by $(T_\varepsilon(t))_{t\geq0}$ the analytic $C_0$-semigroup generates from $L_2^\varepsilon$ in $\elle^2(X,\nu_\varepsilon)$, then Proposition \ref{prop:quasi_regular} also holds for $(T_\varepsilon(t))_{t\geq0}$.
\end{remark}
Let $f\in \fcon_{b,\Phi}^{\infty}(X)$ and let $\lambda>0$. From \cite{CF16}, for every $\varepsilon\in(0,1]$ there exists a sequence $u_{\varepsilon,n}\in\fcon_{b,\Phi}^{3}(X)$ such that
\begin{align}
\label{elliptic_equation_epsilon}
\lambda u_{\varepsilon,n}-L_2^\varepsilon u_{\varepsilon,n}=f+[D_H u_{\varepsilon,n},D_H U_\varepsilon-D_H \widetilde U_{\varepsilon,n}]_{H}=:f_{\varepsilon,n},
\end{align}
where $\widetilde U_{\varepsilon,n}$ has been defined in Proposition \ref{prop:conv_Ueps}. For every $\lambda>0$ and every $\varepsilon\in(0,1]$, we denote by $R(\lambda,L_2^\varepsilon)$ and by $R(\lambda,L_2)$ the resolvent of $L_2^\varepsilon$ and of $L_2$, respectively.
\begin{proposition}
\label{prop:collezione1}
Let $\varepsilon\in(0,1]$. For every $f\in\elle^2(X,\nu_\varepsilon)$ there exists a sequence $(f_{\varepsilon,n})\subseteq W^{1,2}(X,\nu_\varepsilon)$ such that $f_{\varepsilon,n}\rightarrow f$ in $\elle^2(X,\nu_\varepsilon)$, the sequence $(R(\lambda,L_2^\varepsilon)f_n)\subseteq \fcon_{b,\Phi}^3(X)$ converges to $R(\lambda,L_2^\varepsilon)f$ in $W^{2,2}(X,\nu_\varepsilon)$ as $n\rightarrow+\infty$, and
\begin{align}
\label{stima_risolvente}
\|R(\lambda,L_2^\varepsilon)f\|_{W^{2,2}(X,\nu_\varepsilon)}\leq \max\{\sqrt 2,\lambda^{-1},\lambda^{-1/2}\}\|f\|_{\elle^2(X,\nu_\varepsilon)}.
\end{align}
If $f\in W^{1,2}(X,\nu_\varepsilon)$, then $D_{H}f_{\varepsilon,n}\rightarrow D_{H}f$ in $\elle^1(X,\nu_\varepsilon;H)$ as $n\rightarrow+\infty$. Further, if $f\in\fcon_{b,\Phi}^\infty(X)$ then we can choose as $(f_{\varepsilon,n})$ the sequence defined in \eqref{elliptic_equation_epsilon}.
\end{proposition}
\begin{proof}
The proof is contained in the proof of \cite[Propositions 5.6 \& 5.10]{CF16}. The unique point which we have to show is that $D_Hf_{\varepsilon,n}\rightarrow D_Hf$ in $\elle^1(X,\nu_\varepsilon;H)$ as $n\rightarrow +\infty$. As usual, by density it is enough to prove the statement for $f\in \fcon_{b,\Phi}^\infty(X)$. In this case, $f_{\varepsilon,n}$ is the function defined in \eqref{elliptic_equation_epsilon}. By differentiating \eqref{elliptic_equation_epsilon} along $H$ in the direction of $e_i$, $i=1,\ldots,n$, by multiplying by $[D_Hu_{\varepsilon,n},e_i]_H$ and by summing up $i$ from $1$ to $n$, taking \eqref{esplicit_DL2_varepsilon} into account we get
 \begin{align*}
 & (\lambda+1) |D_Hu_{\varepsilon,n}|^2_H
 -\sum_{i=1}^n \left(L^\varepsilon_2[D_H u_{\varepsilon,n}(\cdot),e_i]_{H}\right)[D_Hu_{\varepsilon,n},e_i]_H+[D^2_{H}U_\varepsilon D_Hu_{\varepsilon,n},D_H u_{\varepsilon,n}]_{H} \\
 = & [D_Hf,D_Hu_{\varepsilon,n}]_H
 +[D^2_{H}u_{\varepsilon,n}D_Hu_{\varepsilon,n},D_H U_\varepsilon-D_H \widetilde U_{\varepsilon,n}]_{H}
 +[D_Hu_{\varepsilon,n},(D^2_{H} U_\varepsilon-D^2_{H}\widetilde U_{\varepsilon,n})D_Hu_{\varepsilon,n}]_{H}.
 \end{align*}
We notice that the above equality holds $\mu$-a.e. in $X$. From \cite[Proposition 4.4]{CF16}, there exists a positive constant $K$, independent of $n$, such that
\begin{align}
\label{stima_D_H_uepsilonn}
\|D_Hu_{\varepsilon,n}\|_\infty\leq K\|D_Hf\|_\infty, \quad \mu\textup{-a.e. in }X, \ n\in\N.
\end{align}
The convexity of $U_\varepsilon$, the definition of $L_2^\varepsilon$ (and an integration by parts) and \eqref{stima_D_H_uepsilonn} imply that
\begin{align*}
\int_X|D^2_Hu_{\varepsilon,n}|_{\mathcal H_2(X)}d\nu_\varepsilon
\leq & C\bigg(\|D_Hf\|_\infty+\sigma\int_X|D^2_Hu_{\varepsilon,n}|^2_{\mathcal H_2(X)}d\nu_\varepsilon \\
&  +\frac{1}{\sigma}\int_X|D_H U_\varepsilon-D_H \widetilde U_{\varepsilon,n}|^2_Hd\nu_\varepsilon
+ \int_X|D^2_{H} U_\varepsilon-D^2_{H}\widetilde U_{\varepsilon,n}|^2_{\mathcal H_2(H)}d\nu_\varepsilon\bigg),
\end{align*}
for some positive constant $C$ independent on $n$, for every $\sigma>0$. By choosing $\sigma=(2C)^{-1}$ and taking into account \eqref{convergenza_Uepsneta1} and \eqref{convergenza_Uepsneta2}, it follows that there exists a positive constant $M$, independent of $n$, such that 
\begin{align}
\label{stima_D2_H_u_epsilonn}
 \|D_H^2u_{\varepsilon,n}\|_{\elle^2(X,\nu_\varepsilon;\mathcal H_2(H))}\leq M, \quad n\in\N.
\end{align}
We are ready to prove that $D_Hf_{\varepsilon,n}\rightarrow D_Hf$ in $\elle^1(X,\nu_\varepsilon;H)$ as $n\rightarrow +\infty$. Indeed, we have
\begin{align*}
\|D_Hf_{\varepsilon,n}- D_Hf\|_{\elle^1(X,\nu_\varepsilon;H)}
\leq & \int_X|[D^2_{H}u_{\varepsilon,n}D_Hu_{\varepsilon,n},D_H U_\varepsilon-D_H \widetilde U_{\varepsilon,n}]_{H}|d\nu_\varepsilon \\
&  +\int_X|[D_H u_{\varepsilon,n},(D^2_{H} U_\varepsilon-D^2_{H}\widetilde U_{\varepsilon,n})D_Hu_{\varepsilon,n}]_{H}|d\nu_\varepsilon \\
& \rightarrow 0, \quad n\rightarrow+\infty,
\end{align*}
from \eqref{convergenza_Uepsneta1}, \eqref{convergenza_Uepsneta2}, \eqref{stima_D_H_uepsilonn} and \eqref{stima_D2_H_u_epsilonn}.
\end{proof}

\begin{proposition}
\label{prop:collezione2}
\begin{enumerate}[(i)]
\item  Let $\varepsilon\in(0,1]$, let $t>0$, let $f\in \elle^2(X,\nu_\varepsilon)$ and let $(f_{\varepsilon,n})$ be as in Proposition \ref{prop:collezione1}. The sequence $(T_\varepsilon(t)f_{\varepsilon,n})\subseteq \fcon_{b,\Phi}^3(X)$ converges to $T_\varepsilon(t)f$ in $W^{2,2}(X,\nu_\varepsilon)$ as $n\rightarrow+\infty$. 
\item For every $f\in C_b(X)$ and every $t>0$, we have $T_\varepsilon(t)f\rightarrow T(t)f$ as $\varepsilon\rightarrow0^+$ weakly in $W^{2,2}(X,\nu)$.
\end{enumerate} 
\end{proposition}
\begin{proof}
{$\bf (i)$}. The proof is identical to that of \cite[Proposition 2.8(i)]{AngFerPal18}.

\vspace{2mm}
{$\bf(ii)$}. Let $f\in C_b(X)$, and let $\varepsilon\in(0,1]$. Since $U(x)\geq U_\varepsilon(x)$ for every $x\in X$, from \eqref{stima_risolvente} it follows that the family $\{R(\lambda,L_2^\varepsilon)f:\varepsilon\in(0,1]\}$ is bounded in $W^{2,2}(X,\nu)$. Hence, there exists a vanishing sequence $(\varepsilon_n)\subseteq (0,1]$ such that $R(\lambda,L_{\varepsilon_n})f$ weakly converges to $g\in W^{2,2}(X,\nu)$ in $W^{2,2}(X,\nu)$. Same arguments as in the proof of \cite[Theorem 1.2]{CF16} imply that $g=R(\lambda,L_2)f$. Since every vanishing sequence $(\varepsilon_n)\subseteq (0,1]$ admits a subsequence $(\varepsilon_{k_n})$ such that $R(\lambda,L_2^{\varepsilon_{k_n}})f\rightarrow R(\lambda,L_2)f$ in $W^{2,2}(X,\nu)$ as $n\rightarrow+\infty$, it follows that $R(\lambda,L^\varepsilon_2)f\rightarrow R(\lambda,L_2)f$ weakly in $W^{2,2}(X,\nu)$ as $\varepsilon\rightarrow0^+$. We recall that $(T_\varepsilon(t))_{t\geq0}$ is an analytic semigroup for every $\varepsilon\in(0,1]$. It follows that
\begin{align*}
T_\varepsilon(t)f  =\frac{1}{2\pi i}\int_\sigma e^{\lambda t}R(\lambda,L_2^\varepsilon)fd\lambda, \quad t>0, 
\quad T(t)f =\frac{1}{2\pi i}\int_\sigma e^{\lambda t}R(\lambda,L_2)fd\lambda, \quad t>0, 
\end{align*}
where $\sigma$ is an unbounded curve in $\mathbb C$ which leaves on the left a sector containing the spectrum of $L^\varepsilon_2$. We remark that it is possible to choose $\sigma$ independent of $\varepsilon$. Therefore, for every $g\in C_b(X)$ we have
\begin{align*}
\int_XT_\varepsilon(t)fgd\nu
=\frac{1}{2\pi i}\int_X\left(\int_\sigma e^{\lambda t}R(\lambda,L^\varepsilon_2)fd\lambda\right) gd\nu_\varepsilon
=\frac{1}{2\pi i}\int_\sigma e^{\lambda t}\left(\int_XR(\lambda, L_2^\varepsilon)fgd\nu\right)d\lambda.
\end{align*}
By the dominated convergence theorem, letting $\varepsilon\rightarrow0^+$ we get
\begin{align*}
\lim_{\varepsilon\rightarrow0^+}\int_XT_\varepsilon(t)fgd\nu
= \int_\sigma e^{\lambda t}\left(\int_XR(\lambda,L_2)fgd\nu\right)d\lambda
= \int_XT(t)fgd\nu.
\end{align*}
This proves that $T_\varepsilon(t)f\rightarrow T(t)f$ as $\varepsilon\rightarrow0^+$ weakly in $\elle^2(X,\nu)$. The proof of the convergence of $D_HT_\varepsilon(t)f$ and of $D^2_HT_\varepsilon(t)f$ is analogous.
\end{proof}

\section{Proofs}
\label{app_B}
\subsection*{Proof of Proposition \ref{prop:mink_ineq}}
\begin{proof}
The measurability of $y\mapsto T(t)(F(\cdot,y))$ with respect to $\mathcal M$ can be proved by means of approximations with simple functions. Moreover, by density it is enough to prove the statement for functions $F$ of the form
\begin{align*}
F(x,y):=\sum_{j=1}^nf_j(y)\chi_{A_j}(x), 
\end{align*}
where $A_j\in B(X)$ and $f_j$ is a $\mathcal M$-measurable function for every $j=1,\ldots,n$. Then, for every $t\geq0$ and $\nu$-a.e. $x\in X$ we have
\begin{align*}
\left|T(t)(F(\cdot,y))(x)\right|^p
= \left|\sum_{j=1}^nT(t)(\chi_{A_j})(x)f_j(y)\right|^p, \quad y\in M.
\end{align*}
From the Minkowski inequality for $\gamma$ it follows that
\begin{align}
\left(\int_M\left|T(t)(F(\cdot,y))(x)\right|^p\gamma(dy)\right)^{1/p}
\leq & \sum_{j=1}^n\left(\int_M|T(t)(\chi_{A_j})(x)f_j(y)|^p\gamma(dy)\right)^{1/p} \notag \\
= & \sum_{j=1}^nT(t)(\chi_{A_j})(x)\left(\int_M|f_j(y)|^p\gamma(dy)\right)^{1/p} \notag \\
= & T(t)\left(\sum_{j=1}^m\chi_{A_j}\left(\int_M|f_j|^pd\gamma\right)^{1/p}\right)(x),
\label{minkowski_ineq_T_1}
\end{align}
for every $t\geq0$ and $\nu$-a.e. $x\in X$. 

Let us fix $x\in X$. If $x\in \cup_{i=1}^nA_i$ and $k\in\{1,\ldots,n\}$ satisfies $x\in A_k$, then, since $\{A_1,\ldots,A_n\}$ are pairwise disjoint, it follows that
\begin{align}
\sum_{j=1}^m\chi_{A_j}(x)\left(\int_M|f_j|^pd\gamma\right)^{1/p}
= & \chi_{A_k}(x)\left(\int_M|f_k|^pd\gamma\right)^{1/p}
=\left(\int_M|\chi_{A_k}(x)f_k|^pd\gamma\right)^{1/p} \notag \\
= & \left(\int_M\left|\sum_{j=1}^m\chi_{A_j}(x)f_j\right|^pd\gamma\right)^{1/p}
= \left(\int_M|F(x,\cdot)|^pd\gamma\right)^{1/p}.
\label{minkowski_ineq_T_2}
\end{align}
If $x\in \left(\cup_{i=1}^nA_i\right)^c$ then
\begin{align}
\sum_{j=1}^m\chi_{A_j}(x)\left(\int_M|f_j|^pd\gamma\right)^{1/p}=0=  \left(\int_M|F(x,\cdot)|^pd\gamma\right)^{1/p}. 
\label{minkowski_ineq_T_3}
\end{align}
Combining \eqref{minkowski_ineq_T_1}, \eqref{minkowski_ineq_T_2} and \eqref{minkowski_ineq_T_3} we get the thesis.
\end{proof}

\subsection*{Proof of Proposition \ref{prop:pointwise_stime_grad}}
\begin{proof}
We split the proof into three steps. In the former we prove that for every $\varepsilon\in(0,1]$ and every $f\in\fcon_{b,\Phi}^\infty(X)$ there exists a $\nu_\varepsilon$-measurable set $N_\varepsilon\subseteq X$ such that $\nu_\varepsilon(N_\varepsilon)=0$ and $|D_{H}T_\varepsilon(t)f(x)|_{H} \leq e^{-t} T_\varepsilon(t)|D_{H}f|_{H}(x)$ for every $x\in X\setminus N_\varepsilon$, in the second one we prove the statement for every $f\in\fcon_{b,\Phi}^\infty(X)$ and every $p\in[1,+\infty)$, in the latter we conclude.

\vspace{2mm}
{\bf STEP $1$}.
Let $\varepsilon\in(0,1]$, let $t>0$, let $f\in \fcon_{b,\Phi}^{\infty}(X)$, and let $g$ be a positive, bounded and continuous function on $X$. For every $\sigma>0$ we set $\eta_\sigma:[0,+\infty)\rightarrow[0,+\infty)$ defined by $\eta_\sigma(\xi):=\sqrt {\xi+\sigma}-\sqrt {\sigma}$, for every $\xi\in [0,+\infty)$. We notice that $\eta_\sigma$ satisfies the following:
\begin{align}
\label{eta_properties}
(i)\  \eta_\sigma(\xi)\leq \sqrt \xi, \quad (ii)\ \xi\eta_\sigma'(\xi)\geq \frac12\eta_\sigma(\xi), \quad (iii)\  \eta'_\sigma(\xi)+2\xi\eta_\sigma''(\xi)\geq0, \quad \xi\in[0,+\infty).
\end{align}
To lighten the notations we set $w^\varepsilon_n(s):=|D_{H}T_\varepsilon(s)f_{\varepsilon,n}|_{H}^2$, for every $n\in\N$ and every $s\geq0$, where the sequence $(f_{\varepsilon,n})$ is as in Proposition \ref{prop:collezione2}. We introduce the function
\begin{align*}
G(s):=\int_X\eta_\sigma(w_n^\varepsilon(t-s))T_\varepsilon(s)gd\nu_\varepsilon, \quad t\in[0,+\infty), s\in[0,t], \ n\in\N.
\end{align*}
The smoothness of $T_\varepsilon(t)f_{\varepsilon,n}$ (see Proposition\ref{prop:collezione2}$(i)$) implies that
\begin{align*}
\frac{d}{ds}\eta_\sigma(w^\varepsilon_n(t-s))
= & \eta_\sigma'(w^\varepsilon_n(t-s))\frac{d}{ds}[D_{H}T_\varepsilon(t-s)f_{\varepsilon,n},D_HT_\varepsilon(t-s)f_{\varepsilon,n}]_{H} \\
= & -2\eta'_\sigma(w^\varepsilon_n(t-s))[D_{H}L_2^\varepsilon [T_\varepsilon(t-s)f_{\varepsilon,n}],D_{H}T_\varepsilon(t-s)f_{\varepsilon,n}]_{H}.
\end{align*}
Then, we get
\begin{align}
G'(s)
= &-2\int_X\eta_\sigma'(w^\varepsilon_n(t-s)[D_{H}L_2^\varepsilon T_\varepsilon(t-s)f_{\varepsilon,n},D_{H}T_\varepsilon(t-s)f_{\varepsilon,n}]_{H}T_\varepsilon(s)g d\nu_\varepsilon \notag \\
& + \int_X\eta_\sigma(w^\varepsilon_n(t-s))L_2^\varepsilon T_\varepsilon(s)g d\nu_\varepsilon.
\label{stima_grad_grad_1}
\end{align}
Let us take into account the second addend in the right-hand side \eqref{stima_grad_grad_1}. From the definition of $L_2^\varepsilon$ we get
\begin{align}
\int_X\eta_\sigma(w^\varepsilon_n(t-s))L_2^\varepsilon T_\varepsilon(s)g d\nu_\varepsilon 
= & -\int_X\eta_\sigma'(w^\varepsilon_n(t-s))[(D_{H}w^\varepsilon_n(t-s)),D_{H}T_\varepsilon(s)g]_{H}d\nu_\varepsilon \notag\\
= &- \int_X[D_{H}w_n^\varepsilon(t-s),D_{H}(\eta_\sigma'(w^\varepsilon_n(t-s))T_\varepsilon(s)g)]_{H}d\nu_\varepsilon \notag \\
& +\int_X\eta_\sigma''(w^\varepsilon_n(t-s))|D_{H}w^\varepsilon_n(t-s)|_{H}^2T_\varepsilon(s)gd\nu_\varepsilon \notag \\
= & \int_X \eta_\sigma'(w^\varepsilon_n(t-s))T_\varepsilon(s)g L_2^\varepsilon [w^\varepsilon_n(t-s)]d\nu_\varepsilon \notag \\
& +\int_X\eta_\sigma''(w^\varepsilon_n(t-s))|D_{H}w^\varepsilon_n(t-s)|_{H}^2T_\varepsilon(s)gd\nu_\varepsilon.
\label{stima_grad_grad_2}
\end{align}
It follows that
\begin{align}
G'(s)
= & 2\int_X\eta_\sigma'(w^\varepsilon_n(t-s))T_\varepsilon(s)g
\left(\frac12 L^\varepsilon_2[w^\varepsilon_n(t-s)]-[D_{H}L_2^\varepsilon T_\varepsilon(t-s)f_{\varepsilon,n},D_{H}T_\varepsilon(t-s)f_{\varepsilon,n}]_{H}\right)d\nu_\varepsilon \notag \\
& + \int_X\eta_\sigma''(w^\varepsilon_n(t-s))|D_{H}w^\varepsilon_n(t-s)|_{H}^2T_\varepsilon(s)gd\nu_\varepsilon.
\label{stima_grad_grad_3}
\end{align}
Long but straightforward computations reveal that for every $x\in X$ we have
\begin{align}
\label{stima_grad_3_12}
\frac12L^\varepsilon_2[w^\varepsilon_n(t-s)](x)
= {\rm Tr}[(D^2_{H}T_\varepsilon(t-s)f_{\varepsilon,n})^2(x)]_{H}+\left[L^\varepsilon_2[D_{H}T_\varepsilon(t-s)f_{\varepsilon,n}](x),D_HT_\varepsilon(t-s)f_{\varepsilon,n}(x)\right]_H,
\end{align}
where in the second addend of the right-hand side of \eqref{stima_grad_3_12}, for every $x\in X$ the term $L_2^\varepsilon[D_HT_\varepsilon(t-s)f_{\varepsilon,n}](x)$ is seen as a fixed element of $H$. By combining \eqref{esplicit_DL2_varepsilon} and \eqref{stima_grad_3_12} we get
\begin{align}
\frac12 & L^\varepsilon_2w^\varepsilon_n(t-s)-[D_{H}L_2^\varepsilon T_\varepsilon(t-s)f_{\varepsilon,n},D_{H}T_\varepsilon(t-s)f_{\varepsilon,n}]_{H} \notag \\
= & {\rm Tr}[(D^2_{H}T_\varepsilon(t-s)f_{\varepsilon,n})^2]_{H}+|D_{H}T_\varepsilon(t-s)f_{\varepsilon,n}|_{H}^2
+[D^2_{H}U_\varepsilon D_{H}T_\varepsilon(t-s)f_{\varepsilon,n},D_{H}T_\varepsilon(t-s)f_{\varepsilon,n}]_{H}\notag  \\
\geq &  |D_{H}^2T_\varepsilon(t-s)f_{\varepsilon,n}|^2_{\mathcal H_2(H)}+|D_{H}T_\varepsilon(t-s)f_{\varepsilon,n}|_{H}^2,
\label{stima_grad_grad_4}
\end{align}
where in the last inequality we have used \eqref{traccia_hilbert_sc}, the symmetry of $D^2_H T_{\varepsilon}(t-s)f_{\varepsilon,n}$ as operator from $H\times H$ onto $\R$ and the convexity of $U_\varepsilon$. Further,
\begin{align}
|D_{H}w^\varepsilon_n(t-s)|^2_{H}
=& 4 |D_{H}^2T_\varepsilon(t-s)f_{\varepsilon,n}D_{H}T_\varepsilon(t-s)f_{\varepsilon,n}|^2_{H} \notag \\
\leq & 4|D_{H}^2T_\varepsilon(t-s)f_{\varepsilon,n}|^2_{\mathcal H_2(H)}|D_{H}T_\varepsilon(t-s)f_{\varepsilon,n}|_{H}^2.
\label{stima_grad_grad_5}
\end{align}
By collecting \eqref{stima_grad_grad_1}-\eqref{stima_grad_grad_5} and recalling \eqref{eta_properties} we infer that
\begin{align*}
G'(s)
\geq & 2 \int_X\left(\eta_\sigma'(w^\varepsilon_n(t-s))+2w^\varepsilon_n(t-s)\eta_\sigma''(w^\varepsilon_n(t-s))\right)T_\varepsilon(s)g{\rm Tr}[(D^2_{H}T_\varepsilon(t-s)f_{\varepsilon,n})^2]_{H}d\nu_\varepsilon \\
& + 2\int_X\eta_\sigma'(w^\varepsilon_n(t-s))w^\varepsilon_n(t-s)T_\varepsilon(s)gd\nu_\varepsilon \\
\geq & \int_X\eta_\sigma(w^\varepsilon_n(t-s))T_\varepsilon(s)gd\nu_\varepsilon \\
= & G(s).
\end{align*}
This gives
\begin{align*}
G(s)\geq G(0)e^{s}, \quad s\in[0,t]. 
\end{align*}
In particular, if we choose $s=t$ we get
\begin{align*}
\int_X\left(\sqrt {\sigma+|D_{H}T_\varepsilon(t)f_{\varepsilon,n}|_{H}^2}-\sqrt{\sigma}\right)gd\nu_\varepsilon
= & G(0)\leq e^{-t}G(t) \\
= &e^{-t}\int_X\left(\sqrt{\sigma +|D_{H}f_{\varepsilon,n}|^2}-\sqrt{\sigma}\right)T_\varepsilon(t)gd\nu_\varepsilon.
\end{align*}
By letting $\sigma\rightarrow0^+$ and by applying the symmetry of $T_\varepsilon(t)$ on $\elle^2(X,\nu_\varepsilon)$ we infer that
\begin{align*}
\int_X|D_{H}T_\varepsilon(t)f_{\varepsilon,n}|_{H} gd\nu_\varepsilon
\leq e^{-t}\int_XT_\varepsilon(t)|D_{H}f_{\varepsilon,n}|_{H} gd\nu_\varepsilon.
\end{align*} 
Letting $n\rightarrow+\infty$ and recalling Proposition \ref{prop:collezione2}$(i)$ we infer that
\begin{align*}
\int_X|D_{H}T_\varepsilon(t)f|_{H} gd\nu_\varepsilon
\leq e^{-t}\int_XT_\varepsilon(t)|D_{H}f|_{H} gd\nu_\varepsilon,
\end{align*} 
for every $\varepsilon\in(0,1]$. The arbitrariness of $g$ implies that for every $\varepsilon\in(0,1]$ there exists a $\nu_\varepsilon$-measurable set $N_\varepsilon\subseteq X$ such that $\nu_\varepsilon(N_\varepsilon)=0$ and $|D_{H}T_\varepsilon(t)f(x)|_{H} \leq e^{-t} T_\varepsilon(t)|D_{H}f|_{H}(x)$ for every $x\in X\setminus N_\varepsilon$. 

\vspace{2mm}
{\bf STEP $2$}.
Let us consider a decreasing and vanishing sequence $(\varepsilon_n)\subseteq (0,1]$ and let us set $N:=\cup_{n\in\N}N_{\varepsilon_n}$, where $N_{\varepsilon_n}$ has been defined in Step $1$. Since $\nu_\varepsilon$ and $\nu$ are equivalent measures on $X$ for every $\varepsilon>0$, it follows that $\nu(N)=0$ and 
\begin{align}
\label{stima_epsilon_n}
|D_{H}T_{\varepsilon_n}(t)f(x)|_{H} \leq e^{-t}T_{\varepsilon_n}(t)|D_{H}f|_{H}(x), \quad \forall x\in X\setminus N, \ \forall n\in\N. 
\end{align}
Let $g$ be a positive, bounded and continuous function. By multiplying both the sides of \eqref{stima_epsilon_n} by $g$ and integrating on $X$ with respect to $\nu$ we get
\begin{align}
\label{stima_pun_grad_grad_quasi_fin}
\int_X|D_{H}T_{\varepsilon_n}(t)f|_{H}gd\nu \leq e^{-t}\int_XT_{\varepsilon_n}(t)|D_{H}f|_{H}gd\nu, \quad \forall n\in\N. 
\end{align}
Let us consider the left-hand side of \eqref{stima_pun_grad_grad_quasi_fin}. Since $g$ is positive, it follows that
\begin{align}
\label{conv_epsilon_destra}
\int_X|D_{H}T_{\varepsilon_n}(t)f|_{H} gd\nu
= \int_X|gD_{H}T_{\varepsilon_n}(t)f|_{H} d\nu.
\end{align}
Let us set $V_n:=gD_{H}T_{\varepsilon_n}(t)f$. For every $\Phi\in \elle^\infty(X,\nu;H)$, from Proposition \ref{prop:collezione2}$(ii)$ we have
\begin{align*}
\int_X[V_n,\Phi]_Hd\nu=
\int_X[D_HT_{\varepsilon_n}(t)f,g\Phi]_Hd\nu\rightarrow 
\int_X[D_HT(t)f,g\Phi]_Hd\nu
= \int_X[gD_HT(t)f,\Phi]_Hd\nu,
\end{align*}
as $n\rightarrow+\infty$. We recall that $\elle^\infty(X,\nu;H)=\left(\elle^1(X,\nu;H)\right)^*$ (see \cite[Chapter IV, Theorem 1]{DU77}). Therefore, $V_n\rightarrow gD_HT(t)f$ weakly in $\elle^1(X,\nu;H)$ as $n\rightarrow+\infty$. This implies that
\begin{align}
\label{stima_epsilon_sinistra}
\int_X|D_HT(t)f|_Hgd\nu
= \int_X|gD_HT(t)f|_Hd\nu \leq \liminf_{n\rightarrow+\infty}\int_X|V_n|_Hd\nu
= \liminf_{n\rightarrow+\infty}\int_X|D_HT_{\varepsilon_n}(t)f|_Hgd\nu.
\end{align}
\eqref{stima_pun_grad_grad_quasi_fin}, \eqref{conv_epsilon_destra} and \eqref{stima_epsilon_sinistra} give
\begin{align}
\label{stima_grad_grad_prima_parte_quasi}
\int_X|D_HT(t)f|_Hgd\nu\leq 
e^{-t}\liminf_{n\rightarrow+\infty}\int_XT_{\varepsilon_n}(t)|D_Hf|_Hgd\nu.
\end{align}
From Proposition \ref{prop:collezione2}$(ii)$ we infer that
\begin{align*}
\lim_{n\rightarrow+\infty}\int_XT_{\varepsilon_n}(t)|D_Hf|_Hgd\nu=\int_XT(t)|D_Hf|_Hgd\nu
\end{align*}
From this formula and \eqref{stima_grad_grad_prima_parte_quasi} it follows that
\begin{align*}
\int_X|D_HT(t)f|_Hgd\nu\leq e^{-t}\int_XT(t)|D_Hf|_Hgd\nu.
\end{align*}
The arbitrariness of $g$ implies that
\begin{align}
\label{grad_est_p=1}
|D_{H}T(t)f(x)|_{H}\leq e^{-t}(T(t)|D_{H}f|_{H})(x), \quad t\geq0, \ \nu{\textup{ a.e. }}x\in X,
\end{align}
which gives the thesis for $f\in \fcon_{b,\Phi}^{\infty}(X)$ and $p=1$. 

Let $p>1$.
From \eqref{jensen_semigroup} (with $q=p$ and $p=1$) and \eqref{grad_est_p=1} we infer that
\begin{align*}
|D_{H}T(t)f(x)|_{H}^p
\leq &e^{-pt}(T(t)|D_{H}f|_{H}(x))^p
\leq  e^{-pt}(T(t)|D_{H}f|_{H}^p)(x)
, \quad t\geq0, \ \nu\textup{-a.e.} \  x\in X.
\end{align*}

\vspace{2mm}
{\bf STEP $3$}.
The general case follows by approximation. Let $p\in[1,+\infty)$, let $f\in W^{1,p}(X,\nu)$ and let $(g_n)\subseteq \fcon_{b,\Phi}^\infty(X)$ be such that $g_n\rightarrow f$ in $W^{1,p}(X,\nu)$. We get
\begin{align*}
\int_X |D_{H}T(t)(g_n-g_m)|_{H}^pd\nu
\leq \int_XT(t)|D_{H}(g_n-g_m)|_{H}^pd\nu=\int_X|D_{H}g_n-D_{H}g_m|_{H}^pd\nu,
\end{align*}
for every $n,m\in\N$, where in the last equality we have used the fact that $\nu$ is an invariant measure for $(T(t))_{t\geq0}$. This implies that $(D_{H}T(t)g_n)$ is a Cauchy sequence in $\elle^p(X,\nu;H)$ for every $t\geq0$, and we notice that $T(t)g_n\rightarrow T(t)f$ in $\elle^p(X,\nu)$ as $n\rightarrow+\infty$. The fact that $(D_{H},W^{1,p}(X,\nu))$ is a closed operator in $ \elle^p(X,\nu)$ implies that $T(t)f\in W^{1,p}(X,\nu)$ and 
\begin{align*}
D_{H}T(t)f=\elle^p-\lim_{n\rightarrow+\infty}D_{H}T(t)g_n, \quad t\geq0. 
\end{align*}
Let us fix $t\geq0$. From Step $2$ we know that for every $n\in\N$ there exists $N_n\in B(X)$, with $\nu(N_n)=0$, such that
\begin{align*}
|D_{H}T(t)g_{n}(x)|_{H}^p
\leq e^{-pt}(T(t)|D_{H}g_n|_{H}^p)(x), \quad x\in X\setminus N_n.
\end{align*}
We set $N:=\cup_n N_n$. Then, $\nu(N)=0$ and
\begin{align}
\label{stima_grad_grad_ultima}
|D_{H}T(t)g_{n}(x)|_{H}^p
\leq e^{-pt}(T(t)|D_{H}g_n|_{H}^p)(x), \quad x\in X\setminus N, \ n\in\N.
\end{align}
Finally, we recall that $|D_{H}g_n|_{H}\rightarrow |D_{H}f|_{H}$ in $\elle^1(X,\nu)$ as $n\rightarrow+\infty$, and therefore for every $t\geq0$ we have $T(t)|D_{H}g_n|_{H}\rightarrow T(t)|D_{H}f|_{H}$ in $\elle^1(X,\nu)$ as $n\rightarrow+\infty$. We consider a subsequence $(g_{k_n})$ such that $T(t)|D_{H}g_{k_n}|\rightarrow T(t)|D_{H}f|$ and $D_{H}T(t)g_{k_n}\rightarrow D_{H}T(t)f$ pointwise $\nu$-a.e. in $X$. Letting $n\rightarrow+\infty$ in \eqref{stima_grad_grad_ultima} with $g_n$ replaced by $g_{k_n}$ we get
\begin{align*}
|D_{H}T(t)f|_{H}^p
\leq e^{-pt}T(t)|D_{H}f|_{H}^p, \quad\nu{\textup{-a.e. in }}X.
\end{align*}
\end{proof}

\subsection*{Proof of Proposition \ref{prop:stime_gradiente_funzione}}
\begin{proof}
 
 For reader's convenience we split the proof into three steps. In the former we prove that for every  $p\in(1,2]$ and every  $f\in\fcon_{b,\Phi}^\infty(X)$ we have
\begin{align}
\label{stima_grad_funz_prima_parte}
|D_{H}T_\varepsilon(t)f(x)|^p_{H}\leq c_pt^{-p/2}T_\varepsilon(t)|f(x)|^p, \quad t>0, \ \varepsilon>0, \ x\in X.
\end{align}
in the second one we show \eqref{stima_grad_funz_prima_parte_finale} for every  $p\in(1+\infty)$ and every  $f\in\fcon_{b,\Phi}^\infty(X)$, in the latter we conclude.

\vspace{2mm}
{\bf STEP $1$}. Let $f\in\fcon_{b,\Phi}^\infty(X)$, let $(f_{\varepsilon,n})$ be the approximating sequence of $f$ defined in Proposition \ref{prop:collezione1}, let $p\in(1,2]$ and let us fix $t>0$. For every  $\delta>0$ and every  $n\in\N$ we set
\begin{align*}
G^\varepsilon_{\delta,n}(s):=T_\varepsilon(t-s)\left(\left(|T_\varepsilon(s)f_{\varepsilon,n}|^2+\delta\right)^{p/2}-\delta^{p/2} \right), \quad 0<s<t.
\end{align*}
$G^\varepsilon_{\delta,n}$ is differentiable in $(0,t)$, and differentiating it we get
\begin{align}
(G^\varepsilon_{\delta,n})'(s)
= & -L_2^\varepsilon T_\varepsilon(t-s)\left(\left(|T_\varepsilon(s)f_{\varepsilon,n}|^2+\delta\right)^{p/2}-\delta^{p/2} \right) \notag\\
& +T_\varepsilon(t-s)\left(p\left(|T_\varepsilon(s)f_{\varepsilon,n}|^2+\delta\right)^{p/2-1}(T_\varepsilon(s)f_{\varepsilon,n})( L_2^\varepsilon T_\varepsilon(s)f_{\varepsilon,n}) \right) \notag \\
= & T_\varepsilon(t-s) \notag \\
& \left[-L_2^\varepsilon \left(\left(|T_\varepsilon(s)f_{\varepsilon,n}|^2+\delta\right)^{p/2}-\delta^{p/2}\right)+p\left(|T_\varepsilon(s)f_{\varepsilon,n}|^2+\delta\right)^{p/2-1}(T_\varepsilon(s)f_{\varepsilon,n})( L_2^\varepsilon T_\varepsilon(s)f_{\varepsilon,n}) \right],
\label{stima_grad_funz_1}
\end{align}
where we have used the fact that $\left(|T_\varepsilon(s)f_{\varepsilon,n}|^2+\delta\right)^{p/2}\in \fcon_b^3(X)\subseteq D(L_2^\varepsilon)$. Further,
\begin{align}
L_2^\varepsilon\left(|T_\varepsilon(s)f_{\varepsilon,n}|^2+\delta\right)^{p/2}
= & p\left(|T_\varepsilon(s)f_{\varepsilon,n}|^2+\delta\right)^{p/2-1}(T_\varepsilon(s)f_{\varepsilon,n})(L_2^\varepsilon T_\varepsilon(s)f_{\varepsilon,n}) \notag \\
&+p\left(|T_\varepsilon(s)f_{\varepsilon,n}|^2+\delta\right)^{p/2-1}|D_{H}T_\varepsilon(s)f_{\varepsilon,n}|^2_{H} \notag \\
&+p(p-2)\left(|T_\varepsilon(s)f_{\varepsilon,n}|^2+\delta\right)^{p/2-2}(T_\varepsilon(s)f_{\varepsilon,n})^2|D_{H}T_\varepsilon(s)f_{\varepsilon,n}|^2_{H}.
\label{stima_grad_funz_2}
\end{align}
Combining \eqref{stima_grad_funz_1} and \eqref{stima_grad_funz_2} we get
\begin{align*}
(G^\varepsilon_{\delta,n})'(s)
= & -pT_\varepsilon(t-s)\left(\left(|T_\varepsilon(s)f_{\varepsilon,n}|^2+\delta\right)^{p/2-1}|D_{H}T_\varepsilon(s)f_{\varepsilon,n}|^2_{H}\right) \notag \\
&+p(2-p)T_\varepsilon(t-s)\left(\left(|T_\varepsilon(s)f_{\varepsilon,n}|^2+\delta\right)^{p/2-2}(T_\varepsilon(s)f_{\varepsilon,n})^2|D_{H}T_\varepsilon(s)f_{\varepsilon,n}|^2_{H}\right), \quad 0<s<t.
\end{align*}
From the positivity of $(T_\varepsilon(t))_{t\geq0}$ it follows that
\begin{align*}
(G^\varepsilon_{\delta,n})'(s)
\leq p(1-p)T_\varepsilon(t-s)\left(\left(|T_\varepsilon(s)f_{\varepsilon,n}|^2+\delta\right)^{p/2-1}|D_{H}T_\varepsilon(s)f_{\varepsilon,n}|^2_{H}\right), \quad 0<s<t.
\end{align*}
Integrating with respect to $s$ between $0$ and $t$ it follows that
\begin{align*}
& \left(|T_\varepsilon(t)f_{\varepsilon,n}|^2+\delta\right)^{p/2}-\delta^{p/2}-T_\varepsilon(t)\left(\left(|f_{\varepsilon,n}|^2+\delta\right)^{p/2}-\delta^{p/2}\right) \\
\leq &p(1-p) \int_0^tT_\varepsilon(t-s)\left(\left(|T_\varepsilon(s)f_{\varepsilon,n}|^2+\delta\right)^{p/2-1}|D_{H}T_\varepsilon(s)f_{\varepsilon,n}|^2_{H}\right)ds,
\end{align*}
which gives
\begin{align}
\notag
 p(p-1)\! \int_0^tT_\varepsilon(t-s) & \!\left(\left(|T_\varepsilon(s)f_{\varepsilon,n}|^2+\delta\right)^{p/2-1}|D_{H}T_\varepsilon(s)f_{\varepsilon,n}|^2_{H}\right)ds  \\
\leq &\delta^{p/2}+T_\varepsilon(t)\!\left(\left(|f_{\varepsilon,n}|^2+\delta\right)^{p/2}-\delta^{p/2}\right).
\label{stima_grad_funz_3}
\end{align}
Now we estimate $|D_{H}T_\varepsilon(t)f_{\varepsilon,n}|_{H}^p$. Let $s\in(0,t)$. From the semigroup property of $(T_\varepsilon(t))_{t\geq0}$, Proposition \ref{prop:collezione2}$(i)$ and \eqref{stima_grad_grad_tesi} we infer that
\begin{align}
\label{conti_stima_grad_funz_1}
|D_{H}T_\varepsilon(t)f_{\varepsilon,n}|_{H}^p
= & |D_{H}T_\varepsilon(t-s)T_\varepsilon(s)f_{\varepsilon,n}|_H^p
\leq e^{-p(t-s)}T_\varepsilon(t-s)\left(|D_{H}T_\varepsilon(s)f_{\varepsilon,n}|_{H}^p\right).
\end{align}
We multiply and divide the argument of $T_\varepsilon(t-s)$ by $\left(|T_\varepsilon(s)f_{\varepsilon,n}|^2+\delta\right)^{p(2-p)/4}$. By applying \eqref{holder_semigroup} with $q=\frac{2}{p}$ and $q'=\frac{2}{2-p}$ we infer that
\begin{align}
T_\varepsilon(t-s)& \!\left(|D_{H}T_\varepsilon(s)f_{\varepsilon,n}|_{H}^p\right) \notag \\
= & T_\varepsilon(t-s)\!\left(\!\!\left(|T_\varepsilon(s)f_{\varepsilon,n}|^2+\delta\right)^{-\frac{p(2-p)}4}\!|D_{H}\!T_\varepsilon(s)f_{\varepsilon,n}|_{H}^p\left(|T_\varepsilon(s)f_{\varepsilon,n}|^2+\delta\right)^{\frac{p(2-p)}4}\!\right) \notag\\
\leq & \left(T_\varepsilon(t-s)\left(\left(|T_\varepsilon(s)f_{\varepsilon,n}|^2+\delta\right)^{p/2-1}|D_{H}\!T_\varepsilon(s)f_{\varepsilon,n}|_{H}^2\right)\right)^{p/2} \notag \\
&\cdot \left(T_\varepsilon(t-s)\left(\left(|T_\varepsilon(s)f_{\varepsilon,n}|^2+\delta\right)^{p/2}\right)\right)^{(2-p)/2} \notag\\
\leq &\frac{p}{2}\eta^{2/p}T_\varepsilon(t-s)\left(\left(|T_\varepsilon(s)f_{\varepsilon,n}|^2+\delta\right)^{p/2-1}|D_{H}\!T_\varepsilon(s)f_{\varepsilon,n}|_{H}^2\right) \notag \\
& + \frac{2-p}{2}\eta^{2/(p-2)}T_\varepsilon(t-s)\left(\left(|T_\varepsilon(s)f_{\varepsilon,n}|^2+\delta\right)^{p/2}\right),
\label{conti_stima_grad_funz_2}
\end{align}
for every  $\eta>0$, where in the last inequality we have applied the Young's inequality. Further, the positivity of $(T_\varepsilon(t))_{t\geq0}$ and the fact that $p\in(1,2]$ give
\begin{align}
T_\varepsilon(t-s)\left(\left(|T_\varepsilon(s)f_{\varepsilon,n}|^2+\delta\right)^{p/2}\right)
\leq & T_\varepsilon(t-s)\left(|T_\varepsilon(s)f_{\varepsilon,n}|^p+\delta^{p/2}\right)
\leq T_\varepsilon(t-s)\left(T_\varepsilon(s)|f_{\varepsilon,n}|^p+\delta^{p/2}\right) \notag \\
= & T_\varepsilon(t)|f_{\varepsilon,n}|^p+\delta^{p/2}. 
\label{conti_stima_grad_funz_3}
\end{align}
Putting together \eqref{conti_stima_grad_funz_1}, \eqref{conti_stima_grad_funz_2} and \eqref{conti_stima_grad_funz_3} we infer that
\begin{align*}
e^{p(t-s)}|D_{H}T_\varepsilon(t)f_{\varepsilon,n}|_{H}^p
\leq & \frac{p}{2}\eta^{2/p}T_\varepsilon(t-s)\left(\left(|T_\varepsilon(s)f_{\varepsilon,n}|^2+\delta\right)^{p/2-1}|D_{H}\!T_\varepsilon(s)f_{\varepsilon,n}|_{H}^2\right) \\
+ &  \frac{2-p}{2}\eta^{2/(p-2)}\left(T_\varepsilon(t)|f_{\varepsilon,n}|^p+\delta^{p/2}\right).
\end{align*}
Integrating with respect to $s$ between $0$ and $t$ and recalling \eqref{stima_grad_funz_3} we infer that
\begin{align*}
\frac{e^{pt}-1}{p}|D_{H}T_\varepsilon(t)f_{\varepsilon,n}|_{H}^p
\leq & \frac{\eta^{2/p}}{2(p-1)}\left(\delta^{p/2}+T_\varepsilon(t)\!\left(\left(|f_{\varepsilon,n}|^2+\delta\right)^{p/2}-\delta^{p/2}\right)\right) \\
+ &  \frac{2-p}{2}\eta^{2/(p-2)}t\left(T_\varepsilon(t)|f_{\varepsilon,n}|^p+\delta^{p/2}\right).
\end{align*}
Letting $\delta\rightarrow0^+$ we get
\begin{align*}
\frac{e^{pt}-1}{p}|D_{H}T_\varepsilon(t)f_{\varepsilon,n}|_{H}^p
\leq & \frac{\eta^{2/p}}{2(p-1)}T_\varepsilon(t)|f_{\varepsilon,n}|^p+t\frac{2-p}{2}\eta^{2/(p-2)} T_\varepsilon(t)|f_{\varepsilon,n}|^p \\
= & \left(\frac{\eta^{2/p}}{2(p-1)}+\frac{2-p}{2}\eta^{2/(p-2)}t\right)T_\varepsilon(t)|f_{\varepsilon,n}|^p,
\end{align*}
for every  $\eta>0$. Therefore,
\begin{align}
\label{conti_grad_funz_4}
\frac{e^{pt}-1}{p}|D_{H}T_\varepsilon(t)f_{\varepsilon,n}|_{H}^p
\leq \min_{\eta>0}\left\{\frac{\eta^{2/p}}{2(p-1)}+\frac{2-p}{2}\eta^{2/(p-2)}t\right\}T_\varepsilon(t)|f_{\varepsilon,n}|^p
=  \widetilde c_pt^{-\frac p2+1}T_\varepsilon(t)|f_{\varepsilon,n}|^p,
\end{align}
for some positive constant $\widetilde c_p$ only depending on $p$. By dividing both the sides of \eqref{conti_grad_funz_4} by $(e^{pt}-1)p^{-1}$ we infer that
\begin{align*}
|D_{H}T_\varepsilon(t)f_{\varepsilon,n}|_{H}^p
\leq  \frac{pt}{e^{pt}-1} \widetilde c_pt^{-p/2}T_\varepsilon(t)|f_{\varepsilon,n}|^p\leq c_p t^{-p/2}T_\varepsilon(t)|f_{\varepsilon,n}|^p,
\end{align*}
since the function $t\mapsto pt(e^{pt}-1)^{-1}$ is bounded in $(0,+\infty)$. We notice that if $p=2$ computations simplify and we get
\begin{align*}
\frac{e^{2t}-1}{2}|D_{H}T_\varepsilon(t)f_{\varepsilon,n}|_{H}^2
\leq \frac12T_\varepsilon(t)|f_{\varepsilon,n}|^2.
\end{align*}
Hence,
\begin{align*}
|D_{H}T_\varepsilon(t)f_{\varepsilon,n}|_{H}^2
\leq \frac{2t}{2(e^{2t}-1)}t^{-1}T_\varepsilon(t)|f_{\varepsilon,n}|^2,
\end{align*}
which gives $c_2=\frac12$. In both the cases, we get
\begin{align}
\label{stima_finale_seconda_stima_pun_epsilon<2}
|D_{H}T_\varepsilon(t)f_{\varepsilon,n}|^p_{H}\leq c_pt^{-p/2}T_\varepsilon(t)|f_{\varepsilon,n}|^p, \quad t>0, \ p\in(1,2],
\end{align}
and the constant $c_p$ blows up as $p$ approaches $1$. From Proposition \ref{prop:collezione2}$(i)$, up to a subsequence, the left-hand side of \eqref{stima_finale_seconda_stima_pun_epsilon<2} converges to $|D\!_{H}T_\varepsilon(t)f(x)|^p_{H}$ as $n\rightarrow +\infty$ for $\nu_\varepsilon$-a.e. $x\in X$. Let us consider the right-had side of 	\eqref{stima_finale_seconda_stima_pun_epsilon<2}. Since $|f_{\varepsilon,n}|^p\rightarrow |f|^p$ in $\elle^1(X,\nu_\varepsilon)$ it follows that
\begin{align*}
T_\varepsilon(t)|f_{\varepsilon,n}|^p\rightarrow T_\varepsilon|f|^p, \quad n\rightarrow+\infty,  \ {\rm in} \ \elle^1(X,\nu_\varepsilon).
\end{align*}
Hence, up to a subsequence, $T_\varepsilon(t)|f_{\varepsilon,n}|^p(x)\rightarrow T_\varepsilon|f|^p(x)$ as $n\rightarrow+\infty$ for $\nu_\varepsilon$-a.e. $x\in X$. This gives \eqref{stima_grad_funz_prima_parte} for $p\in(1,2]$. 

\vspace{2mm}
{\bf STEP $2$}.
Let $f\in\fcon_{b,\Phi}^\infty(X)$ and let $p\in(1,2]$. Let us multiply both the sides of \eqref{stima_grad_funz_prima_parte} by a positive, bounded and continuous function $g$ and let us integrate on $X$ with respect to $\nu$. We get
\begin{align*}
\int_X|D_{H}T_\varepsilon(t)f|^p_{H}gd\nu
\leq  c_pt^{-p/2}\int_XT_\varepsilon(t)|f|^pgd\nu, \quad t>0.
\end{align*} 
Since $D_HT_\varepsilon(t)f\rightarrow D_HT(t)f$ weakly in $\elle^2(X,\nu)$ as $\varepsilon\rightarrow0^+$ (see Proposition \ref{prop:collezione2}$(ii)$), arguing as in Step $2$ in the proof of Proposition \ref{prop:pointwise_stime_grad} it is possible to prove that $g^{1/p}D_HT_\varepsilon(t)f$ weakly converges to $g^{1/p}D_HT(t)f$ in $\elle^p(X,\nu)$. Hence,
\begin{align*}
\int_X|D_{H}T(t)f|^p_{H}gd\nu
\leq \liminf_{\varepsilon\rightarrow0^+}
\int_X|D_{H}T_\varepsilon(t)f|^p_{H}gd\nu.
\end{align*}
From Proposition \ref{prop:collezione2}$(ii)$ we deduce that
\begin{align*}
\int_XT_\varepsilon(t)|f|^pgd\nu\rightarrow\int_XT(t)|f|^pgd\nu, \quad \varepsilon\rightarrow0^+.
\end{align*}
It follows that
\begin{align*}
\int_X|D_{H}T(t)f|^p_{H}gd\nu
\leq c _pt^{-p/2}\int_XT(t)|f|^pgd\nu, \quad t>0.
\end{align*}
The arbitrariness of $g$ implies that 
\begin{align}
\label{stima_grad_funz_finale_1}
|D_{H}T(t)f(x)|^p_{H}\leq c _pt^{-p/2}T(t)|f|^p(x), \quad \nu\textup{-a.e.} \ x\in X. 
\end{align}
If $p>2$, we apply \eqref{stima_grad_funz_prima_parte} with $p=2$ and \eqref{jensen_semigroup} with $p=1$ and $q=p/2$, and we get
\begin{align*}
|D_{H}T(t)f(x)|^p_{H}
=\left(|D_{H}T(t)f(x)|^2_{H}\right)^{p/2}
\leq \left(\frac12t^{-1}T(t)|f|^2(x)\right)^{p/2}\!\!\!\!
\leq c_pt^{-p/2}T(t)|f|^p(x), \  \nu\textup{-a.e. }x\in X,
\end{align*}
with $c_p=2^{-p/2}$.

\vspace{2mm}
{\bf STEP $3$}.
Let $f\in\elle^p(X,\nu)$ and let $(g_n)\subseteq \fcon_{b,\Phi}^\infty(X)$ converge to $f$ in $\elle^p(X,\nu)$ as $n\rightarrow+\infty$. Replacing $f$ with $g_n-g_m$ in \eqref{stima_grad_funz_finale_1} and integrating on $X$ with respect to $\nu$ we get
\begin{align}
\label{stima_grad_funz_seconda_parte}
\int_X|D_{H}T(t)(g_n-g_m)|^p_{H}d\nu\leq c_pt^{-p/2}\int_XT(t)|g_n-g_m|^pd\nu
= c_pt^{-p/2}\int_X|g_n-g_m|^pd\nu,
\end{align}
where in the last part we have used the fact that $\nu$ is an invariant measure for $(T(t))_{t\geq0}$. 
Conclusion follows by repeating the same computations as in Step $3$ in the proof of Proposition \ref{prop:pointwise_stime_grad}.
\end{proof}

\subsection*{Proof of Proposition \ref{prop:vett_smgr_prop}}
\begin{proof}
{$\boldsymbol {(i)}$}
Let $p\in[1,+\infty)$. Let us prove the result for $F\in \fcon_b^\infty(X;V)$. Let
\begin{align*}
T^V(t)F=\sum_{i=1}^nT(t)f_iv_i, \quad F=\sum_{i=1}^nf_iv_i, \quad f_i\in \fcon_b^\infty(X), \ v_i\in V, \ i=1,\ldots,n,
\end{align*}
Without loss of generality we can assume that $\{v_1,\ldots,v_n\}$ are orthonormal vectors in $V$. We have
\begin{align*}
\nu(F)=\sum_{i=1}^n\left(\int_Xf_id\nu\right) v_i.
\end{align*}
If $p> 2$, by applying the H\"older inequality with $q=p/(p-2)$ and $q'=p/2$ then
\begin{align*}
\|T^V(t)F-\nu(F)\|_{\elle^p(X,\nu;V)}^p
=&\int_X\left(\left| T^V(t)F-\int_XFd\nu\right|^2_V\right)^{p/2} d\nu \\
= & \int_X \left(\sum_{i=1}^n\left|T(t)f_i-\int_Xf_id\nu\right|^2\right)^{p/2}d\nu \\
= & \int_X \left(\sum_{i=1}^n\left|T(t)f_i-\int_Xf_id\nu\right|^2\right)^{p/2-1}\sum_{j=1}^n\left|T(t)f_j-\int_Xf_jd\nu\right|^2d\nu \\
= & \sum_{j=1}^n \int_X \left(\sum_{i=1}^n\left|T(t)f_i-\int_Xf_id\nu\right|^2\right)^{p/2-1}\left|T(t)f_j-\int_Xf_jd\nu\right|^2d\nu \\
\leq & \sum_{j=1}^n\|T^V(t)F-\nu(F)\|_{\elle^p(X,\nu;V)}^{p-2}
\|T(t)f_j-\nu(f_j)\|_{\elle^p(X,\nu)}^2 \\
= & \|T^V(t)F-\nu(F)\|_{\elle^p(X,\nu;V)}^{p-2}\sum_{j=1}^n\|T(t)f_j-\nu(f_j)\|_{\elle^p(X,\nu)}^2,
\end{align*}
\begin{align}
\label{stima_conp_as_vett_scal}
\|T^V(t)F-\nu(F)\|_{\elle^p(X,\nu;V)}^2
\leq \sum_{j=1}^n\|T(t)f_j-\nu(f_j)\|_{\elle^p(X,\nu)}^2,
\end{align}
for every $t\geq0$. If $p=2$ estimate \eqref{stima_conp_as_vett_scal} immediately follows, while if $p\in[1,2)$ then we get
\begin{align}
\|T^V(t)F-\nu(F)\|_{\elle^p(X,\nu;V)}^p
= & \int_X \Big(\sum_{i=1}^n\Big|T(t)f_i-\int_Xf_id\nu\Big|^2\Big)^{p/2}d\nu \notag \\
\leq & \int_X\sum_{i=1}^n|T(t)f_i-\nu(f_i)|^pd\nu \notag \\
= & \sum_{i=1}^n\|T(t)f_i-\nu(f_i)\|_{\elle^p(X,\nu)}^p.
\label{stima_conp_as_vett_scal_2}
\end{align}
From Proposition \ref{prop:comp_asintotico_fz}, \eqref{stima_conp_as_vett_scal} and \eqref{stima_conp_as_vett_scal_2} we conclude that
\begin{align}
\label{stima_conp_as_vett_cil}
\lim_{t\rightarrow+\infty}\|T^V(t)F-\nu(F)\|_{\elle^p(X,\nu;V)}=0.    
\end{align}

Let us consider a function $F\in \elle^p(X,\nu)$ and let $(F_m)\subseteq \fcon_b^\infty(X;V)$ be a sequence converging to $F$ in $\elle^p(X,\nu;V)$ as $n\rightarrow+\infty$. From the properties of Bochner integrals (see \cite[Chapter II, Thereom 4(ii)]{DU77}) and H\"older's inequality we infer that
\begin{align*}
\left\|\int_XF_md\nu -\int_XFd\nu\right\|_{\elle^p(X,\nu;V)}
\leq \|F-F_m\|_{\elle^p(X,\nu;V)}.
\end{align*}
Since $(T^V(t))_{t\geq0}$ is a semigroup of contractions in $\elle^p(X,\nu;V)$ it follows that
\begin{align}
\left\|T^V(t)F-\int_XFd\nu\right\|_{\elle^p(X,\nu;V)}
\leq & \left\|T^V(t)F-T^V(t)F_m\right\|_{\elle^p(X,\nu;V)} \notag \\
&+\left\|T^V(t)F_m-\int_XF_md\nu\right\|_{\elle^p(X,\nu;V)} \notag \\
&+\left\|\int_XF_md\nu -\int_XFd\nu\right\|_{\elle^p(X,\nu;V)} \notag \\
\leq&  2\|F-F_m\|_{\elle^p(X,\nu;V)}+\left\|T^V(t)F_m-\int_XF_md\nu\right\|_{\elle^p(X,\nu;V)}.
\label{stima_as_vett_finale}
\end{align}
Let $\varepsilon>0$. Then, there exists $\overline m\in\N$ such that $\|F-F_{\overline m}\|_{\elle^p(X,\nu;V)}\leq \varepsilon/4$. Further, from \eqref{stima_conp_as_vett_cil} there exists $\overline t$ such that for every $t>\overline t$ we have
\begin{align*}
\left\|T^V(t)F_{\overline m}-\int_XF_{\overline m}d\nu\right\|_{\elle^p(X,\nu;V)}\leq \varepsilon/2.
\end{align*}
Therefore, from \eqref{stima_as_vett_finale} with $m$ replaced by $\overline m$, for every $\varepsilon>0$ there exists $\overline t>0$ such that for every $t>\overline t$ we have
\begin{align*}
\left\|T^V(t)F-\int_XFd\nu\right\|_{\elle^p(X,\nu;V)}\leq \varepsilon.
\end{align*}
This gives the thesis.

\vspace{3mm}
{$\boldsymbol{(ii)}$}.
The further part follows by integrating \eqref{vector_stime_grad_grad} on $X$ with respect to $\nu$ and by recalling that $\nu$ is an invariant measure for $(T(t))_{t\geq0}$. Hence, let us prove the first part. By density we can limit ourselves to prove the result for $F\in\fcon_b^\infty(X;V)$.

Let $p\in[1,+\infty)$, let $t>0$ and let $F\in \fcon_b^\infty(X;V)$ of the form
\begin{align*}
F=\sum_{i=1}^nf_iv_i, \quad f_i\in \fcon_b^\infty(X), \ v_i\in V, \ i=1,\ldots,n.
\end{align*}
Without loss of generality we can assume that $v_i$ are orthonormal vectors in $V$. From \eqref{stima_grad_grad_tesi} with $p=1$, for every $x\in X$ we get
\begin{align*}
|\overline{D_{H}}T^V(t)F(x)|_{H\otimes V}
= & \left(\sum_{i=1}^n|D_{H}T(t)f_i(x)|_{H}^2\right)^{1/2}
\leq  e^{-t}\left(\sum_{i=1}^n(T(t)|D_{H}f_i|_H(x))^{2}\right)^{1/2}.
\end{align*}
Now we apply formula \eqref{minkowski_smgr} with $\gamma$ being the counting measure on $\N$ and with $q=2$. We get
\begin{align}
\label{stima_grad_smrp_p=1}
\left(\sum_{i=1}^n(T(t)|D_{H}f_i|_H(x))^{2}\right)^{1/2}
\leq T(t)\left(\sum_{i=1}^n|D_Hf_i|^2_H\right)^{1/2}\!\!\!(x), \quad x\in X.
\end{align}
If $p=1$ then we conclude. If $p\in(1,+\infty)$ then from \eqref{stima_grad_smrp_p=1} it follows that
\begin{align*}
|\overline{D_{H}}T^V(t)F(x)|_{H\otimes V}^p
\leq & e^{-pt}\left(T(t)\left(\sum_{i=1}^n|D_Hf_i|^2_H\right)^{1/2}\!\!\!(x)\right)^p
\leq e^{-pt}T(t)\left(\sum_{i=1}^n|D_Hf_i|^2_H\right)^{p/2}\!\!\!(x) \\
= & e^{-pt}T(t)|\overline {D_H}F(x)|_{H\otimes V}^p,
\end{align*}
for every $x\in X$, where the last inequality follows from \eqref{jensen_semigroup} with $p=1$ and $q=p/2$.

\vspace{3mm}
{$\boldsymbol{(iii)}$}. The further part follows by integrating \eqref{vector_stima_grad_funz} on $X$ with respect to $\nu$ and by recalling that $\nu$ is an invariant measure for $(T(t))_{t\geq0}$. Hence, let us prove the first part. By density we can limit ourselves to prove the result for $F\in\fcon_b^\infty(X;V)$. We split two cases: in the former we prove the statement for $p\in[2,+\infty)$, in the latter we consider the remaining values of $p$.

\vspace{2mm}
Let $p\geq2$, let $t>0$ and let $F\in \fcon_b^\infty(X;V)$ of the form
\begin{align*}
F=\sum_{i=1}^nf_iv_i, \quad f_i\in \fcon_b^\infty(X), \ v_i\in V, \ i=1,\ldots,n.
\end{align*}
Without loss of generality we can assume that $v_i$ are orthonormal vectors in $V$. From \eqref{stima_grad_funz_prima_parte_finale} with $p=2$ we get
\begin{align*}
|\overline{D_{H}}T^V(t)F|_{H\otimes V}^2
= & \sum_{i=1}^n|D_{H}T(t)f_i|_{H}^2 
\leq  \frac{1}{2t}\sum_{i=1}^nT(t)|f_i|^2 
=  \frac{1}{2t}T(t)\left(\sum_{i=1}^n|f_i|^2\right) = \frac{1}{2t}T(t)\left(|F|_{V}^2\right).
\end{align*}
Hence,
\begin{align*}
|\overline{D_{H}}T^V(t)F|_{H\otimes V}^p
= & \left(|\overline{D_{H}}T^V(t)F|_{H\otimes V}^2\right)^{p/2}
\leq  \frac{1}{(2t)^{p/2}}\left(T(t)\left(|F|_{V}^2\right)\right)^{p/2}
\leq  \frac{1}{(2t)^{p/2}}T(t)\left(|F|_{V}^p\right),
\end{align*}
where the last inequality follows from \eqref{jensen_semigroup} with $p=1$ and $q=p/2$. This gives the thesis since for $p\in[2,+\infty)$ we have $c_p=2^{-p/2}$.

\vspace{2mm}
Let $p\in(1,2)$, and let $t$ and $F$ as above. By applying \eqref{stima_grad_funz_prima_parte_finale} we have
\begin{align*}
|\overline{D_{H}}T^V(t)F(x)|_{H\otimes V}
= \left(\sum_{i=1}^n|D_{H}T(t)f_i(x)|_{H}^2 \right)^{1/2}
\leq c_p^{1/p}t^{-1/2}\left(\sum_{i=1}^n(T(t)|f_i|^p(x))^{2/p} \right)^{1/2}, \quad x\in X.
\end{align*}
Now we apply the Minkowski's integral inequality \eqref{minkowski_smgr} with $\gamma$ being the counting measure on $\N$ and with $q=2/p$. It follows that
\begin{align*}
\sum_{i=1}^n(T(t)|f_i(x)|^p)^{2/p}
\leq \left(T(t)\left(\sum_{i=1}^n|f_i|^2(x)\right)^{p/2}\right)^{2/p}
= (T(t)|F|_V^p(x))^{2/p}, \quad x\in X.
\end{align*}
Therefore, we conclude that
\begin{align*}
|\overline{D_{H}}T^V(t)F(x)|_{H\otimes V}^p
\leq & \frac{c_p}{t^{p/2}}T(t)|F|_V^p(x), \quad x\in X.
\end{align*}

\vspace{3mm}
{$\boldsymbol{(iv)}$}. By density, it is enough to show that \eqref{stima_grad_funz_mista} holds for $F\in\fcon_b^\infty(X;V)$. Let $p\in(1,+\infty)$, let $F\in\fcon_b^\infty(X;V)$ and let $t\in(0,1]$. By applying Proposition \ref{prop:vett_smgr_prop}$(iii)$ we get
\begin{align}
\label{stima_grad_funz_mista_1}
\int_X|\overline{D_{H}}T^V(t)F|_{H\otimes V}^pd\nu
\leq  & \frac{c_p}{t^{p/2}}\int_X|F|^p_V d\nu, \quad t\in(0,1].
\end{align}
If $t\geq1$, then from the semigroup property of $T^V(t)$ and by applying Proposition \ref{prop:vett_smgr_prop}$(ii)$ with $t-1$ and \eqref{int_vector_stima_grad_funz} with $t=1$ we have
\begin{align}
\int_X|\overline{D_{H}}T^V(t)F|_{H \otimes V}^pd\nu
= & \int_X|\overline{D_{H}}T^V(t-1)T^V(1)F|_{H\otimes V}^pd\nu\notag  \\
\leq & e^{-p(t-1)}\int_X|D_{H}T^V(1)F|^p_{H \otimes V}d\nu \notag \\
\leq  & c_p e^{-p(t-1)}\int_X|F|^p_Vd\nu, \quad t>1.
\label{stimagrad_funz_mista_1_bis}
\end{align}
By combining \eqref{stima_grad_funz_mista_1} and \eqref{stimagrad_funz_mista_1_bis} we get the thesis.
\end{proof}

\section*{Declarations}
\paragraph{Ethical approval.}
This declaration is ``not applicable''.

\paragraph{Competing interests.}
This declaration is ``not applicable''.

\paragraph{Authors' contribution.}
This declaration is ``not applicable''.

\paragraph{Funding.}
this research has financially been supported by the Programme ``FIL-Quota Incentivante'' of University of Parma and co-sponsored by Fondazione Cariparma.

\paragraph{Availability of data and materials.}
This declaration is ``not applicable''.

\bibliographystyle{plain}

\begin{thebibliography}{99}


\bibitem{ACF20}
D.~Addona, G.~Cappa, and S.~Ferrari.
\newblock Domains of elliptic operators on sets in {W}iener space.
\newblock {\em Infin. Dimens. Anal. Quantum Probab. Relat. Top.},
  23(1):2050004, 42, 2020.

\bibitem{AD20}
D.~Addona.
\newblock Analyticity of nonsymmetric {O}rnstein-{U}hlenbeck semigroup with
  respect to a weighted {G}aussian measure.
\newblock {\em Potential Analysis}, 54(1):95-117, 2021.

\bibitem{AdMuRo21}
D.~Addona, M.~Muratori, M.~Rossi.
\newblock Equivalence of Sobolev norms in Malliavin spaces.
\newblock{\em J. Funct. Anal.}, 283(7):2022.

\bibitem{AngBigFer23}
L.~Angiuli, D.A.~Bignamini, S.~Ferrari.
\newblock Harnack inequalities with power $p\in(1,+\infty)$ for transition semigroups in Hilbert spaces.
\newblock{\em NoDEA Nonlinear Differential Equations Appl.}, 30(1): 2023.

\bibitem{AngFerPal18}
L.~Angiuli, S.~Ferrari, and D.~Pallara.
\newblock Gradient estimates for perturbed Ornstein-{U}hlenbeck semigroups on
  infinite-dimensional convex domains.
\newblock {\em J. Evol. Equ.}, 19(3):677--715, 2019.


\bibitem{BaEm85}
D.~Bakry, and M.~\'{E}mery.
\newblock{Diffusions hypercontractives},
\newblock{\em {S\'{e}minaire de probabilit\'{e}s, {XIX}, 1983/84}},
\newblock{Lecture Notes in Math.}, {1123}, {177-206}, {Springer, Berlin}, {1985}.

\bibitem{BC11}
H.~H. Bauschke and P.~L. Combettes.
\newblock {\em Convex analysis and monotone operator theory in {H}ilbert
  spaces}.
\newblock CMS Books in Mathematics/Ouvrages de Math\'ematiques de la SMC.
  Springer, New York, 2011.

\bibitem{Big22}
D.A.~Bignamini.
\newblock $L^2$-theory for transition semigroups associated to dissipative systems.
\newblock{\em Stochastics and Partial Differential Equations: Analysis and Computations}, 2022.

\bibitem{BigFer21}
D.A.~Bignamini, S.~Ferrari.
\newblock Regularizing Properties of (Non-Gaussian) Transition Semigroups in Hilbert Spaces.
\newblock{\em 
Potential Analysis}, 2021.

\bibitem{BigFer22}
D.A.~Bignamini, S.~Ferrari.
\newblock On generators of transition semigroups associated to semilinear stochastic partial differential equations. 
\newblock{\em J. Math. Anal. Appl.}, 508(1): 40pp, 2022. 

\bibitem{Bog98}
V.~I. Bogachev.
\newblock {\em Gaussian measures}, volume~62 of {\em Mathematical Surveys and
  Monographs}.
\newblock American Mathematical Society, Providence, RI, 1998.

\bibitem{Bog07}
V.~I. Bogachev.
\newblock {\em Measure theory. {V}ol. {I}, {II}}.
\newblock Springer-Verlag, Berlin, 2007.

\bibitem{CF16}
G.~Cappa and S.~Ferrari.
\newblock Maximal {S}obolev regularity for solutions of elliptic equations in
  infinite dimensional {B}anach spaces endowed with a weighted {G}aussian
  measure.
\newblock {\em J. Differential Equations}, 261(12):7099--7131, 2016.

\bibitem{CF18}
G~Cappa, S.~Ferrari. \newblock Maximal Sobolev regularity for solutions of elliptic equations in Banach spaces endowed with a weighted Gaussian measure: the convex subset case. 
\newblock{\em J. Math. Anal. Appl.}, 458(1):330--331, 2018.

\bibitem{DPG01}
G.~Da~Prato and B.~Goldys.
\newblock Elliptic operators on {$\R^d$} with unbounded coefficients.
\newblock {\em J. Differential Equations}, 172(2):333--358, 2001.

\bibitem{DPL14}
G.~Da~Prato and A.~Lunardi.
\newblock Sobolev regularity for a class of second order elliptic {PDE}'s in
  infinite dimension.
\newblock {\em Ann. Probab.}, 42(5):2113--2160, 2014.

\bibitem{DPL15}
G.~Da Prato, A.~Lunardi.
\newblock Maximal Sobolev regularity in Neumann problems for gradient systems in infinite dimensional domains. 
\newblock{\em Ann. Inst. Henri Poincaré Probab. Stat.}, 51(3):1102--1123, 2015. 

\bibitem{Dav89}
E. B. Davies.
\newblock{Heat kernels and spectral theory}.
\newblock Cambridge Tracts in Mathematics, 1989.

\bibitem{DU77}
J.~Diestel and J.~J. Uhl, Jr.
\newblock {\em Vector measures}.
\newblock American Mathematical Society, Providence, R.I., 1977.
\newblock With a foreword by B. J. Pettis, Mathematical Surveys, No. 15.

\bibitem{Fer19}
S.~Ferrari.
\newblock Sobolev spaces with respect to a weighted {G}aussian measure in
  infinite dimensions.
\newblock {\em Infin. Dimens. Anal. Quantum Probab. Relat. Top.},
  22(4):1950026, 32, 2019.

\bibitem{FeUs00}
D.~Feyel, D. and A.~S.~{U}st\"{u}nel,
\newblock{The notion of convexity and concavity on {W}iener space},
\newblock{J. Funct. Anal.},176(2), 400-428, 2000.


\bibitem{Gr14}
L.~Grafakos.
\newblock {\em Classical {F}ourier analysis}, volume 249 of {\em Graduate Texts
  in Mathematics}.
\newblock Springer, New York, third edition, 2014.

\bibitem{Gr67}
L.~Gross.
\newblock Abstract {W}iener spaces.
\newblock In {\em Proc. {F}ifth {B}erkeley {S}ympos. {M}ath. {S}tatist. and
  {P}robability ({B}erkeley, {C}alif., 1965/66), {V}ol. {II}: {C}ontributions
  to {P}robability {T}heory, {P}art 1}, pages 31--42. Univ. California Press,
  Berkeley, Calif., 1967.

\bibitem{MR92}
Z.~M. Ma and M.~R\"ockner.
\newblock {\em Introduction to the theory of (nonsymmetric) {D}irichlet forms}.
\newblock Universitext. Springer-Verlag, Berlin, 1992.

\bibitem{Nu06}
D.~Nualart.
\newblock {\em The {M}alliavin calculus and related topics}.
\newblock Probability and its Applications (New York). Springer-Verlag, Berlin,
  second edition, 2006.

\bibitem{PrVe14}
M.~Pronk and M.~Veraar.
\newblock Tools for {M}alliavin calculus in {UMD} {B}anach spaces.
\newblock {\em Potential Anal.}, 40(4):307--344, 2014.

\bibitem{Ru91}
W.~Rudin.
\newblock {\em Functional analysis}.
\newblock International Series in Pure and Applied Mathematics. McGraw-Hill,
  Inc., New York, second edition, 1991.

\bibitem{Sh04}
I.~Shigekawa.
\newblock {\em Stochastic analysis}, volume 224 of {\em Translations of
  Mathematical Monographs}.
\newblock American Mathematical Society, Providence, RI, 2004.
\newblock Translated from the 1998 Japanese original by the author, Iwanami
  Series in Modern Mathematics.

\bibitem{Sh97}
Ichiro Shigekawa.
\newblock {$L^p$} contraction semigroups for vector valued functions.
\newblock {\em J. Funct. Anal.}, 147(1):69--108, 1997.

\bibitem{Us95}
A.~S. \"{U}st\"{u}nel.
\newblock {\em An introduction to analysis on {W}iener space}, volume 1610 of
  {\em Lecture Notes in Mathematics}.
\newblock Springer-Verlag, Berlin, 1995.

\bibitem{VN15}
J.~Van~Neerven.
\newblock The $L^p$-Poincar\'e inequality for analytic Ornstein-Uhlenbeck semigroups. 
\newblock Operator semigroups meet complex analysis, harmonic analysis and mathematical physics.
{\em Oper. Theory Adv. Appl.}, {\bf 250}, Birkh\"auser/Springer, Cham, 353--368, 2015.  
\end{thebibliography}

\end{document}